\documentclass{article}
\usepackage{amssymb, verbatim, multirow, xcolor, graphicx, caption, subcaption}
\usepackage{hyperref}
\usepackage{amsmath}

\def\dist{\mathop{\rm dist}}

\def\Riem{\mathop{\rm Rm}}

\def\be{\begin{eqnarray}}
\def\ee{\end{eqnarray}}
\def\beg{\begin{eqnarray*}}
\def\ees{\end{eqnarray*}}
\def\bel{\begin{aligned}}
\def\eel{\end{aligned}}

\def\XXint#1#2#3{{\setbox0=\hbox{$#1{#2#3}{\int}$ }
\vcenter{\hbox{$#2#3$ }}\kern-.6\wd0}}

\newcommand{\qed}{\hfill$\Box$}
\newtheorem{theorem}{Theorem}[section]
\newtheorem{proposition}[theorem]{Proposition}
\newtheorem{lemma}[theorem]{Lemma}

\newtheorem{corollary}[theorem]{Corollary}

\begin{document}

\title{Classification of polytope metrics and complete scalar-flat K\"ahler 4-Manifolds with two symmetries}
\date{September 2015}
\maketitle

\begin{abstract}
We study unbounded 2-dimensional metric polytopes such as those arising as K\"ahler quotients of complete K\"ahler 4-manifolds with two commuting symmetries and zero scalar curvature.
Under a mild closedness condition, we obtain a complete classification of metrics on such polytopes, and as a result classify all possible metrics on on the corresponding K\"ahler 4-manifolds.
If the polytope is the plane or half-plane then only flat metrics are possible, and if the polytope has one corner then the 2-parameter family of generalized Taub-NUTs (discovered by Donaldson) are indeed the only possible metrics.
Polytopes with $n\ge3$ edges admit an $(n+2)$-dimensional family of possible metrics.
\end{abstract}

\section{Introduction}

We study polytope metrics of the kind arising from the reduction of complete K\"ahler 4-manifolds with a pair of commuting holomorphic Killing fields $\mathcal{X}_1$, $\mathcal{X}_2$ and zero scalar curvature; we are particularly interested in the case that $M^4$ is complete.
In the simply-connected case, real-holomorphic fields are generated by potentials; here these are functions $\varphi^i$ with
\be
\omega(\mathcal{X}_i,\,\cdot)\;=\;-d\varphi^i.
\ee
This leads to the so-called moment map $\Phi:M^4\mapsto\mathbb{R}^2$ given by $\Phi(p)=(\varphi^1(p),\varphi^2(p))$, whose image is a polytope called the associated moment polytope \cite{De}\cite{KL}\cite{At}\cite{GS}, which we call $\Sigma^2$.
Since the $\mathcal{X}_i$ are Killing the map $\Phi:M^2\rightarrow\Sigma^2$, which is generically submersive, endows $\Sigma^2$ with a Riemannian metric $g_\Sigma$; the result is a metric polytope $(\Sigma^2,g_{\Sigma})$ called the K\"ahler reduction of $(M^4,J,\omega,\mathcal{X}_1,\mathcal{X}_2)$.
Setting
\be
\begin{aligned}
\mathcal{V}
&\;=\;|\mathcal{X}_1|^2|\mathcal{X}_2|^2-\left<\mathcal{X}_1,\,\mathcal{X}_2\right>^2 \;=\;|\nabla\varphi^1|^2|\nabla\varphi^2|^2-\left<\nabla\varphi^1,\,\nabla\varphi^2\right>^2
\end{aligned}
\ee
and letting $\triangle_\Sigma$ be the laplacian on $(\Sigma^2,g_\Sigma)$, we shall see that
\be
\triangle_\Sigma\sqrt{\mathcal{V}}\;+\;\frac12s\sqrt{\mathcal{V}}\;=\;0
\ee
where $s$ is the scalar curvature on $M^4$, expressed, with $\sqrt{\mathcal{V}}$ as a function on $\Sigma^2$.
This is a version of the Abreu equation from \cite{Ab1}; see the discussion in Section \ref{SectionReduction}.
When $s=0$ we have a Laplace equation and $\sqrt{\mathcal{V}}$ is a natural harmonic coordinate on $\Sigma^2$.

Our work involves only the metric polytopes themselves, whether or not they come from any K\"ahler reduction.
For this reason we take the time to state our hypotheses in a strictly analytic framework, apart from any K\"ahler geometry it may be reduced from.
Unless specifically stated otherwise, we shall always assume our polytopes $(\Sigma^2,g_\Sigma)$ satisfy
\begin{itemize}
\item[A)] (Closedness) $\Sigma^2$ is a closed subset of the $(\varphi^1,\varphi^2)$ coordinate plane.
\item[B)] (Boundary connectedness) The polytope boundary has just one component or no components.
\item[C)] (Natural polytope condition) $\Sigma^2$ is a convex polytope with finitely many edges.
The metric $g_\Sigma$ is $C^\infty$ up to boundary segments and Lipschitz at corners, and each boundary segment or ray is totally geodesic.
\end{itemize}
and our functions $\varphi^1$, $\varphi^2$ satisfy
\begin{itemize}
\item[D)] (Natural boundary conditions) The unit vector fields $\frac{\nabla\varphi^i}{|\nabla\varphi^i|}$ are well-defined and covariant-constant on boundary segments.
\item[E)] (Pseudo-toric condition) The functions $\varphi^1$, $\varphi^2$ are $C^\infty$ and $[\nabla\varphi^1,\nabla\varphi^2]=0$ on $(\Sigma^2,g_\Sigma)$.
\item[F)] (Pseudo-ZSC condition) With $\mathcal{V}=|\nabla\varphi^1|^2|\nabla\varphi^2|^2-\left<\nabla\varphi^1,\nabla\varphi^2\right>$, we have $\triangle_\Sigma\sqrt{\mathcal{V}}=0$.
\end{itemize}

Conditions (C), (D), and (E) are automatic when $(\Sigma^2,g_\Sigma)$ is the K\"ahler reduction of some $M^4$ (see \S\ref{SubSecPolyTopeEdges}).
Condition (F) is simply that $M^4$ also has zero scalar curvature.

Conditions (A) and (B), however, are genuine restrictions in the sense that there exist manifolds $(M^4,J,\omega,\mathcal{X}_1,\mathcal{X}_2)$ with associated polytopes that violate either or both; see Examples 7 and 8 below.
We will not discuss such polytopes, except to say that (A) is violated when some symmetry has a ``zero at infinity'' as in Example 7, and (B) is violated when the zero locus $\{\mathcal{X}_1=0\}\cup\{\mathcal{X}_2=0\}$ is disconnected, as in Example 8.
For a discussion of why these cases are more difficult, and remain unresolved, see the final remark in Section \ref{SectionQuarterPlaneAndGeneral}.

We completely classify polytope metrics under conditions (A)-(F).
Should $(\Sigma^2,g_\Sigma)$ be a K\"ahler reduction, it is well-known that all data on $M^4$ can be reconstructed from $g_\Sigma$ (see \S\ref{SectionReduction}), so obviously our theorems about $(\Sigma^2,g_\Sigma)$ have profound implications for K\"ahler manifolds $(M^4,J,\omega,\mathcal{X}_1,\mathcal{X}_2)$ with commuting symmetries.
\begin{theorem}[Compact Polytope] \label{TheoremFlatCompactPolytope}
If $(\Sigma^2,g_\Sigma)$ is compact, and even if (F) is relaxed to $\triangle_\Sigma\sqrt{\mathcal{V}}\ge0$, then $g_\Sigma$ is flat. (C.f. Corollary \ref{CorCompactFlatPolytope}.)
\end{theorem}
The following corollary is well-known:
\begin{corollary}[Compact $M^4$ with $s\le0$]\label{CorFlatCompactPolytope}
Assume $(M^4,J,\omega,\mathcal{X}_1,\mathcal{X}_2)$ has $s\le0$.
If the associated metric polytope $(\Sigma^2,g_\Sigma)$ is compact, then $M^4$ is flat. (C.f. Corollary \ref{CorCompactFlatFourManifold}.)
\end{corollary}
The simplest non-compact case is when $\Sigma^2$ is metrically complete; this is the subject of our first substantial results.
Note also the different sign required of $s$.
\begin{theorem}[Polytope is complete] \label{TheoremFlatRTwo}
Assume $(\Sigma^2,g_\Sigma)$ is a geodesically complete metric polytope with $\triangle_{\Sigma}\sqrt{\mathcal{V}}\le0$ (this is a relaxation of (A) and (F)).
Then $g_\Sigma$ is flat.
(C.f. Theorem \ref{ThmRTwoClassification}.)
\end{theorem}
\begin{corollary} \label{CorFlatRTwo}
Assume $(M^4,J,\omega,\mathcal{X}_1,\mathcal{X}_2)$ is simply connected and has $s\ge0$.
If the fields are no-where zero and no-where equal, then $M^4$ is flat. (C.f. Corollary \ref{CorRTwoClassification}.)
\end{corollary}
Our third theorem deals with the case that $\Sigma^2$ is a half-plane.
Unlike the previous results, its proof requires the full strength of (A)-(F).
\begin{theorem}[Polytope is a half-plane] \label{TheoremFlatHalfPlane}
If $(\Sigma^2,g_\Sigma)$ is a closed half-plane, then $g_\Sigma$ is flat.
(C.f. Theorem \ref{ThmHalfPlaneFlat}.)
\end{theorem}
\begin{corollary} \label{CorFlatHalfPlane}
Assume $(M^4,J,\omega,\mathcal{X}_1,\mathcal{X}_2)$ is simply connected, scalar-flat, and has associated polytope $(\Sigma^2,g_\Sigma)$ with just one edge.
Then $M^4$ is flat. (C.f. Corollary \ref{CorHalfPlaneFlat}.)
\end{corollary}

When $\mathcal{X}_1,\mathcal{X}_2$ have common zeros, the situation is more complicated.
Since at least the papers of Donaldson \cite{Do2} and Abreu--Sena-Dias \cite{AS}, it has been known that if $M^4$ is scalar-flat and has a polytope with one or more vertices, then $M^4$ admits at least a two-parameter family of scalar flat metrics (up to homothety).
As we show, when $\Sigma^2$ has just one vertex, these are {\it all} such metrics: one of these metrics is flat, one is the Taub-NUT metric, and the rest are the achiral generalized Taub-NUT metrics of Donaldson (from section 6 of \cite{Do2}).
\begin{theorem}[Polytope has one vertex] \label{TheoremQuarterPlane}
Assume $(\Sigma^2,g_\Sigma)$ has a single corner, so after affine recombination of $\varphi^1$, $\varphi^2$ we may assume $\Sigma^2$ is the closed first quadrant.
Up to homothety, the metric $g_\Sigma$ lies within a 2-parameter family of metrics (cf. Theorem \ref{ThmGeometricSingleVertex}).
\end{theorem}
\begin{corollary} \label{CorQuarterPlane}
Under the hypotheses of Theorem \ref{TheoremFlatHalfPlane}, suppose also that $\mathcal{X}_1$, $\mathcal{X}_2$ have a single common zero.
Up to homothety, $M^4$ has precisely a 2-parameter family of metrics which give it a metric of zero scalar curvature (cf. Corollary \ref{CorOneVertexMetricClassification}).
\end{corollary}

The situation of more than two edges is different still.

\begin{theorem}[Generic case] \label{TheoremGeneralCase}
Assume $\Sigma$ has $n\ge3$ many edges.
Under conditions (A)-(F), $\Sigma^2$ admits precisely an $(n+2)$-parameter family of metrics.
(C.f. Theorem \ref{ThmDegreesOfFreedomGeneralPolytope}.)
\end{theorem}
\begin{corollary} \label{TheoremGeneralCase}
Suppose $(\Sigma,g_\Sigma)$ is the K\"ahler reduction of a manifold $(M^4,J,\omega,\mathcal{X}_1,\mathcal{X}_2)$, and $\Sigma^2$ has $n\ge3$ edges, has no edges at infinity, and has connected boundary.
If $M^4$ is scalar flat, its metric is one in a $(n+2)$-parameter family of possible metrics.F
(C.f. Corollary \ref{CorDegreesOfFreedomGeneralScalarFlat}.)
\end{corollary}

We give a detailed construction of these polytope metrics in the proof of Corollary \ref{CorDegreesOfFreedomGeneralScalarFlat}, and a detailed recipe for constructing metrics on $M^4$ from the metric on $\Sigma^2$ in \S\ref{SubSectionReconstructionOfTheMetric}.
In Example 6 below, we explicitly construct all possible such metrics on those $M^4$ whose polytope has two corners; these include the metrics on $T\mathbb{C}P^2$ such as the Eguchi-Hanson metric, and the metrics of LeBrun \cite{LeB}.
Of course obtaining the $M^4$ metric from its polytope is nothing new, see \cite{G} \cite{Ab1}, but the method we find most helpful is a variation on what usually appears in the literature.
In particular we avoid the ``K\"ahler potential'' techniques of \cite{G} and \cite{Ab1}, which lead to a 4th order scalar PDE, and instead use essentially equivalent system of second order PDEs.

The proofs of Theorems \ref{TheoremFlatCompactPolytope} and \ref{TheoremFlatRTwo} are reasonably self-contained and easy.
Theorem \ref{TheoremFlatHalfPlane} requires the preparatory work of \S\ref{SubSecGlobalHarmonicCoordBehavior} and \S\ref{SubSecPolyTopeEdges} and the difficult analysis of \S\ref{SectionAnalysisHalfPlane}, which classifies non-negative solutions of $x(\varphi_{xx}+\varphi_{yy})-\varphi_x=0$ on the half-plane $H^2=\{x>0\}$ with zero boundary values.
The {\it geometric} classifications in Theorems \ref{TheoremQuarterPlane} and \ref{TheoremGeneralCase} are essentially corollaries of the analytic classification.
The equation $x(\varphi_{xx}+\varphi_{yy})-\varphi_x=0$ is in fact a geometric PDE that appeared first in a paper of Donaldson's \cite{Do1}; see Proposition \ref{PropDivOfVarphi} and the discussion at the beginning of Section \ref{SectionAnalysisHalfPlane}.

The main result of Section \S\ref{SectionAnalysisHalfPlane} is a Liouville-type result for the degenerate-elliptic equation $x(f_{xx}+f_{yy})+3f_x=0$ by combining the techniques of blow-ups, barriers, and Fourier analysis.
A certain duality relation exists between $x(f_{xx}+f_{yy})+3f_x=0$ with arbitrary boundary conditions on the half-plane $H^2$ and the PDE $x(\varphi_{xx}+\varphi_{yy})-\varphi_x=0$ with zero boundary conditions on $H^2$, so the Liouville-type theorem for the first gives the desired classification of all non-negative solutions for the second: these all have the form $Cx^2$ for $C>0$.

Degenerate-elliptic PDE have long been studied in their own right (eg. \cite{KN} \cite{Hor}).
The particular equations $x\left(g_{xx}+g_{yy}\right)+\nu{g}_x=0$ and the corresponding parabolic equation $g_t=x\left(g_{xx}+g_{yy}\right)+\nu{g}_x$ have been studied in a variety of contexts, almost exclusively for $\nu\ge0$ (although there is a duality between the cases $\nu>1$ and $\nu<1$: if $g$ solves the equation for $\nu$ then $\tilde{g}=x^{\nu-1}g$ solves it for $\tilde\nu=2-\nu$).

The operator $\triangle_\mu=x\left(\frac{\partial^2}{\partial{x}^2}+\frac{\partial^2}{\partial{y}^2}\right)$ can be considered a Laplace operator for the right half-plane with a hypoerbolic metric.
Thus from a naive point of view, $g_t=\triangle_\mu{g}+\nu{g}_x$ can be considered a heat flow through an isotropic but inhomogeneous medium of hyperbolic specific heat density, with a directional transfer bias indicated by $\nu$.
Less naively, this equation (or similar equations such as the Heston equation, equation (\ref{EqnHeston})) has been studied in connection with mathematical finance, stochastic PDE and Feynman-Kac theory \cite{Hes}, mathematical biology \cite{EM}, and the porous medium equation \cite{DH}.
The change of variables  $s=\frac12x^{2}$, $t=\sqrt{2}y$ turns $x(g_{xx}+g_{yy})-g_{x}=0$ into $s{f}_{ss}+f_{tt}=0$, which has been studied in connection with population dynamics \cite{EM}.

However, much of this study has been related to existence/uniqueness of flows or boundary value problems, H\"older estimates, admissibility of boundary values, and so on.
Much of this work has been local or confined to pre-compact domains, although there are global results such as global H\"older estimates.
To this author's knowledge, this paper contains the first Liouville-type result for these equations.

{\bf Remark.} Our 4-manifolds need not be ``toric'' exactly, but only have commuting holomorphic Killing fields.
Our work takes place mostly on the moment polytope where the vertex angles are irrelevant.
There is no particular need that the momentum construction rebuild any actual 4-manifold.

{\bf Remark}. Our conclusions will hold directly on certain orbifolds with commuting holomorphic Killing fields.
However we do not explore the significance of our results to the situation of ``marked polytopes'' \cite{LeTo} and the like.

{\bf Remark.} Notice that Theorem \ref{TheoremFlatCompactPolytope} requires $s\le0$ and Theorem \ref{TheoremFlatRTwo} requires $s\ge0$, but $s=0$ is necessary for the other three theorems.
We conjecture that Theorem \ref{TheoremFlatHalfPlane} (that the half-plane is flat) continues to hold if $s\ge0$.
Further, we conjecture that if $s\ge0$ is some pre-specified function, then Theorem \ref{TheoremQuarterPlane} holds with the ``up to homothety'' removed.
However the techniques of \S\ref{SectionAnalysisHalfPlane} will not work without $s\equiv0$, and a more sophisticated analysis is required.
See the discussion at the beginning of \S\ref{SectionAnalysisHalfPlane}.

{\bf Remark.} It is important to mention two issues not addressed in this paper.
The first is the question of how many of the metrics in our $(n+2)$-dimensional family are repeats.
Conceivably even the dimension of this family could be reduced by the action of some diffeomorphism group, although we do not believe this to be the case.

The second is that, although the polytopes $\Sigma^2$ with $n\ge3$ edges have an $(n+2)$-dimensional family of metrics, we do not address the issue of which of these may actually come from the reduction of some $M^4$.
This is a another way the family of scalar-flat metrics on $M^4$ may possibly be smaller than $(n+2)$-dimensional.

{\bf Acknowledgements}.
The author would like to thank Xiuxiong Chen, Ryan Hynd, Philip Gressman, and Camelia Pop for a number of helpful conversations that provided some valuable insights.

\section{K\"ahler Reduction} \label{SectionReduction}

Here we discuss the detailed relationship between the K\"ahler reduction $(\Sigma^2,g_\Sigma)$ and $(M^4,J,\omega,\mathcal{X}_1,\mathcal{X}_2$).
The philosophy is ``Guillemin's principle,'' to wit, {\it sympelctic coordinates are useful in K\"ahler geometry} \cite{G}, \cite{Ab2}.
This principle has been excellently developed by a number of authors \cite{G} \cite{Ab1} \cite{Do1} \cite{AS} \cite{CDG}, but begging the knowledgeable reader's forbearance, we redevelop some aspects as suits the present interest.
First we indicate this section's primary milestones.

In \S\ref{SecFundamentals} we perform the symplectic reduction itself, performing the Arnold-Liouville \cite{Ar} reduction process, and relating the symplectic coordinates $(\varphi^1,\varphi^2,\theta^1,\theta^2)$ to the complex-analytic coordinates $(z^1,z^2)$.
The resulting holomorphic volume form is a quarter of the parallelochoron volume, which is
\be
\mathcal{V}\;\triangleq\;|\mathcal{X}_1|^2|\mathcal{X}_2|^2-\left<\mathcal{X}_1,\mathcal{X}_2\right>^2
\;=\;|\nabla\varphi^1|^2|\nabla\varphi^2|^2-\left<\nabla\varphi^1,\nabla\varphi^2\right>^2, \label{EqnDefoOfV}
\ee
which implies that the scalar curvature of $(M^4,J,\omega)$ is $s=-\triangle\log\mathcal{V}$.
Our notational set-up this well-known construction will be useful later.

Even with this elliptic relation, and even with a sign on $s$, the behavior of $\log\mathcal{V}$ on $M^4$ may be hard to understand, and so in  \S\ref{SubSectionReductionToPolytope} we bring everything down to the metric polytope $(\Sigma^2,g_\Sigma)$ itself, which is naturally a Riemannian manifold that may or may not have boundary.
In the inherited metric we show $\triangle_\Sigma\mathcal{V}^{\frac12}+\frac12s\mathcal{V}^{\frac12}=0$, which is a version of the Abreu equation: see (10) of \cite{Ab1}.
Of course when $M^4$ is scalar-flat then $\sqrt{\mathcal{V}}$ is harmonic on $(\Sigma^2,g_\Sigma)$, and we have natural isothermal coordinates $(x,y)$ on $\Sigma^2$ where $x=\sqrt{\mathcal{V}}$ and $y$ is defined by $J_\Sigma{d}y=dx$.
We compute the Gauss curvature $K_\Sigma$ of $\Sigma^2$ in these coordinates.

In \S\ref{SubSectionReconstructionOfTheMetric} we show how to reconstruct the metrics $g_\Sigma$, $g$ simply by knowing $\varphi^1$, $\varphi^2$ as functions of $(x,y)$.
We also compute the intrinsic Gaussian curvature $K_\Sigma$ of $(\Sigma^2,g_\Sigma)$ in a second, independent way, and show that although $K_\Sigma$ is not necessarily signed, after a canonical conformal change it is indeed signed.

In \S\ref{SubSecGlobalHarmonicCoordBehavior} we prove the isothermal coordinate system $(x,y)$ actually has no critical points, and that the map $X=(x,y)^T$ sends $\Sigma^2$ onto the right half-plane $\{x\ge0\}$ in a one-to-one fashion, and the isothermal coordinates are therefore global.
This fundamentally relies on condition (B), connectedness of the polytope boundary.

Finally in \S\ref{SubSecPolyTopeEdges} we partially control the behaviors of $\varphi^1$, $\varphi^2$, $x$, and $y$ near the polytope edges.
This provides us with boundary conditions that are necessary to our analysis of \S\ref{SectionAnalysisHalfPlane} and also to our penultimate construction, in Theorem \ref{ThmDegreesOfFreedomGeneralPolytope}.

\subsection{Fundamentals} \label{SecFundamentals}

We have a simply connected K\"ahler 4-manifold with commuting Killing fields $(M^4,J,\omega,\mathcal{X}_1,\mathcal{X}_2)$.
The momentum construction consists of finding potentials $\varphi^1$, $\varphi^2$ satisfying $\omega(\mathcal{X}_i,\cdot)=-d\varphi^i$, traditionally called {\it momentum variables} or {\it action coordinates}.
Two additional coordinates $\theta^1$, $\theta^2$, called {\it cyclic variables} or {\it angle coordinates}, are defined by taking a transversal to the $\{\mathcal{X}_1,\mathcal{X}_2\}$ distribution and then pushing the natural $\mathbb{R}^2$ variables forward along the action of $\{\mathcal{X}_1,\mathcal{X}_2\}$.
Obviously the values of $\theta^1$, $\theta^2$ are not canonical, but the fields $\frac{\partial}{\partial\theta^1}$, $\frac{\partial}{\partial\theta^2}$ are canonical.
This construction gives a full coordinate system $\{\varphi^1,\varphi^2,\theta^1,\theta^2\}$ with fields
\be
\begin{array}{ll}
\frac{\partial}{\partial\varphi^1}\;=\;\frac{\left|\mathcal{X}_2\right|^2}{\mathcal{V}}\,\nabla\varphi^1 \,-\,\frac{\left<\mathcal{X}_1,\,\mathcal{X}_2\right>}{\mathcal{V}}\,\nabla\varphi^2, & \frac{\partial}{\partial\theta^1}\;=\;\mathcal{X}_1, \\
\\
\frac{\partial}{\partial\varphi^2}\;=\;-\frac{\left<\mathcal{X}_1,\,\mathcal{X}_2\right>}{\mathcal{V}}\,\nabla\varphi^1 \,+\,\frac{\left|\mathcal{X}_1\right|^2}{\mathcal{V}}\,\nabla\varphi^2, & \frac{\partial}{\partial\theta^2}\;=\;\mathcal{X}_2.
\end{array}
\ee
Letting $G^{-1}$ be the matrix
\be
G^{-1}\;=\;\left(\begin{array}{cc}
|\mathcal{X}_1|^2 & \left<\mathcal{X}_1,\,\mathcal{X}_2\right> \\
\left<\mathcal{X}_1,\,\mathcal{X}_2\right> & |\mathcal{X}_2|^2
\end{array}\right),
\ee
then the ordered frame $\frac{\partial}{\partial\varphi^1},\frac{\partial}{\partial\varphi^2},\frac{\partial}{\partial\theta^1},\frac{\partial}{\partial\theta^2}$ produces the metric, complex structure, and symplectic form
\be
g\;=\;\left(\begin{array}{c|c}
G & 0 \\
\hline
0 &\;\; G^{-1}\;\;
\end{array}\right), \quad
J\;=\;\left(\begin{array}{c|c}
0 & -G^{-1} \\
\hline
\;G & 0
\end{array}\right), \quad
\omega\;=\;\left(\begin{array}{c|c}
0 & \;-Id\; \\
\hline
Id\; & 0
\end{array}\right).
\label{EqnsGJOmegaM}
\ee

{\bf Remark}.
Compare (\ref{EqnsGJOmegaM}) to compare (4.8) of \cite{G} or (2.2), (2.3) of \cite{Ab2}.
The more usual construction sets $G=\nabla^2\,u$ for a function $u$ called the K\"ahler potential.
In the present work we do not find this formulation useful.
In our formulation, of course it may be objected that expressing the metric in terms of inner products is redundant.
Still, this formulation lends itself to the study of second order instead of fourth order PDEs, and is useful in relating the polytope metric $g_{\Sigma}$ with the original metric $g$.

\begin{lemma}[Symplectic and holomorphic coordinates on $M^4$] \label{LemmaHoloSymplCoordRelation}
The complex-valued functions
\be
\begin{aligned}
z^1&\;=\;f^1(\varphi^1,\,\varphi^2)\,+\,\sqrt{-1}\,\theta^1 \\
z^2&\;=\;f^2(\varphi^1,\,\varphi^2)\,+\,\sqrt{-1}\,\theta^2
\end{aligned} \label{EqnsHoloCoordForm}
\ee
form a holomorphic coordinate chart provided $f^1$, $f^2$ satisfy
\be
\begin{aligned}
&df^1\;=\;Jd\theta^1\;=\;\frac{|\mathcal{X}_2|^2}{\mathcal{V}}\,d\varphi^1\,-\,\frac{\left<\mathcal{X}_1,\,\mathcal{X}_2\right>}{\mathcal{V}}\,d\varphi^2 \\
&df^2\;=\;Jd\theta^2\;=\;-\frac{\left<\mathcal{X}_1,\,\mathcal{X}_2\right>}{\mathcal{V}}\,d\varphi^1\,+\,\frac{|\mathcal{X}_1|^2}{\mathcal{V}}\,d\varphi^2. \label{EqnsDFi}
\end{aligned}
\ee
Further, $dJd\theta^i=0$, so indeed such functions $f^1$, $f^2$ can be found locally, and therefore holomorphic coordinates of the form (\ref{EqnsHoloCoordForm}) exist on $(M^4,J,\omega)$.
\end{lemma}
{\it Proof}.
Using that $\bar\partial=\frac12\left(d+\sqrt{-1}Jd\right)$ on functions gives that $\bar\partial{z}^i=0$ if and only if (\ref{EqnsDFi}) holds.
That $dJd\theta^i=0$ follows from the preservation of $J$ under the action of $\mathcal{X}_1$, $\mathcal{X}_2$ (a lengthy though elementary calculation can give $dJd\theta^i=0$ directly).
\qed

One easily determines the holomorphic frame and coframe
\be
\begin{aligned}
\frac{\partial}{\partial{z}^1} \;=\;& \frac12\left(|\mathcal{X}_1|^2\frac{\partial}{\partial\varphi^1} \,+\, \left<\mathcal{X}_1,\,\mathcal{X}_2\right>\frac{\partial}{\partial\varphi^2} \,-\, \sqrt{-1}\,\frac{\partial}{\partial\theta^1}\right)
\;=\;\frac12\left(\nabla\varphi^1\,-\,\sqrt{-1}\,\mathcal{X}_1\right) \\
\frac{\partial}{\partial{z}^2} \;=\;& \frac12\left(\left<\mathcal{X}_1,\,\mathcal{X}_2\right>\frac{\partial}{\partial\varphi^1} \,+\, |\mathcal{X}_2|^2\frac{\partial}{\partial\varphi^2} \,-\,\sqrt{-1}\,\frac{\partial}{\partial\theta^2}\right)
\;=\; \frac12\left(\nabla\varphi^2\,-\,\sqrt{-1}\,\mathcal{X}_2\right) \\
dz^1 \;=\;&\frac{|\mathcal{X}_2|^2}{\mathcal{V}}dx^1\,-\,\frac{\left<\mathcal{X}_1,\,\mathcal{X}_2\right>}{\mathcal{V}}dx^2\,+\,\sqrt{-1}\,d\theta^1
\;=\; Jd\theta^1 \,+\,\sqrt{-1}\,d\theta^1 \\
dz^2 \;=\;&-\frac{\left<\mathcal{X}_1,\,\mathcal{X}_2\right>}{\mathcal{V}}dx^1\,+\,\frac{|\mathcal{X}_1|^2}{\mathcal{V}}dx^2\,+\,\sqrt{-1}\,d\theta^2
\;=\; Jd\theta^2 \,+\,\sqrt{-1}\,d\theta^2
\end{aligned} \label{EqnCxCoords}
\ee
so in the holomorphic frame the Hermitian metric and associated volume element are simply
\be
\begin{aligned}
&h_{i\bar\jmath}\;=\;
\frac12
\left(\begin{array}{cc}
|\mathcal{X}_1|^2 & \left<\mathcal{X}_1,\,\mathcal{X}_2\right> \\
\left<\mathcal{X}_1,\,\mathcal{X}_2\right> & |\mathcal{X}_2|^2
\end{array}\right)
\;=\;\frac12G^{-1}, \\
&\det\,h_{i\bar\jmath}\;=\;\frac14\mathcal{V}.
\end{aligned}
\ee
\begin{proposition} \label{PropRicciFormAndScalarOnM}
The Ricci form and scalar curvature of $(M^4,J,\omega)$ are
\be
\begin{aligned}
\rho&\;=\;-\sqrt{-1}\partial\bar\partial\log\,\mathcal{V} \\
s&\;=\;-\triangle\log\,\mathcal{V}.
\end{aligned} \label{EqnsRicScal}
\ee
\end{proposition}
{\it Proof}. Textbook formulas. \qed

{\bf Remark}. We have given, in Lemma \ref{LemmaHoloSymplCoordRelation} and (\ref{EqnCxCoords}), a very prosaic version of what is usually a expressed as a Legendre transform between the so-called K\"ahler and symplectic potentials, as developed in \cite{G} (see also \cite{Ab1}).
As we have chosen not to work with K\"ahler potentials, these Legendre transform methods won't apply.
In any case the workaday formulation above is probably more useful in the present context.

\subsection{Reduction of $M^4$ to a Metric Polytope} \label{SubSectionReductionToPolytope}

Let $(\Sigma^2,g_\Sigma)$ be the Riemannian quotient (leaf space) of $M^4$ by the Killing actions of $\mathcal{X}_1$, $\mathcal{X}_2$.
The commutativity of $\mathcal{X}_1$, $\mathcal{X}_2$ imply the functions $\varphi^1$, $\varphi^2$ pass to the quotient and there constitute a coordinate system.
As is well-known, the differentiable map $\Phi:\Sigma^2\rightarrow\mathbb{R}^2$, $\Phi(p)=(\varphi^1(p),\varphi^2(p))^T$ is 1-1 with differentiable inverse, so $\Sigma^2$ may be identified with its image in the $(\varphi^1,\varphi^2)$-coordinate plane.
This constitutes the moment polytope of $M^4$.

For clarity, objects on $\Sigma^2$ will be indicated with a subscript, so for instance $s_{\Sigma}$ and $s$ indicate the scalar curvatures on $(\Sigma^2,g_{\Sigma})$ and $(M^4,J,\omega)$, respectively.
The inherited metric and natural complex structure are
\be
&&g_\Sigma \;=\;\mathcal{V}^{-1}\left(\begin{array}{cc}
|\mathcal{X}_2|^2 & -\left<\mathcal{X}_1,\,\mathcal{X}_2\right> \\
-\left<\mathcal{X}_1,\,\mathcal{X}_2\right> & |\mathcal{X}_1|^2\\
\end{array}\right) \;=\; G, \label{EqnSigmaMetric} \\
&&J_\Sigma\;=\;\mathcal{V}^{-\frac12}\left(\begin{array}{cc}
\left<\mathcal{X}_1,\,\mathcal{X}_2\right> & -|\mathcal{X}_1|^2 \\
|\mathcal{X}_2|^2 & -\left<\mathcal{X}_1,\,\mathcal{X}_2\right>
\end{array}\right). \label{EqnJSigma}
\ee
Of course $J_\Sigma$ is not inherited from $J$, but is the Hodge star of $(\Sigma^2,g_\Sigma)$.
\begin{proposition}[$\Sigma^2$ and $M^4$ Laplacian relationship] \label{PropPojectionLaplacian}
If $f:M^4\rightarrow\mathbb{C}$ is any function on $M$ that is $\mathcal{X}_1$- and $\mathcal{X}_2$-invariant, then $\triangle{f}$ is $\mathcal{X}_1$, $\mathcal{X}_2$ invariant, so $f$, $\triangle{f}$ are naturally functions on $\Sigma^2$.
As functions on $\Sigma^2$, the two Laplacians are related by
\be
\triangle{}f
&=&\triangle_\Sigma{f}\,+\,\left<\nabla_{\Sigma}\log{\mathcal{V}}^{\frac12},\,\nabla_{\Sigma}f\right>_\Sigma
\ee
\end{proposition}
{\it Proof}. A function $f:M^4\rightarrow\mathbb{C}$ is invariant under $\mathcal{X}_1$, $\mathcal{X}_2$ if and only if it is a function of $\varphi^1$, $\varphi^2$ only.
Noting that $\det(g)=1$ and $\det(g_\Sigma)=\mathcal{V}^{-1}$ we have
\be
\begin{aligned}
\triangle{f}&\;=\;\frac{\partial}{\partial\varphi^i}\left(g^{ij}\frac{\partial{f}}{\partial\varphi^j}\right)
\;=\;\frac{\partial}{\partial\varphi^i}\left(g_\Sigma^{ij}\frac{\partial{f}}{\partial\varphi^j}\right) \\
&\;=\;\mathcal{V}^{\frac12}\frac{\partial}{\partial\varphi^i}\left(g_\Sigma^{ij}\mathcal{V}^{-\frac12}\frac{\partial{f}}{\partial\varphi^j}\right)
\;+\;\,g_\Sigma^{ij}\frac{\partial\log{\mathcal{V}}^{\frac12}}{\partial\varphi^i}\frac{\partial{f}}{\partial\varphi^j} \\
&\;=\;\triangle_\Sigma{f}\,+\,\frac12\left<\nabla_\Sigma\log\mathcal{V},\,\nabla_\Sigma{f}\right>_\Sigma.
\end{aligned}
\ee
\qed

\begin{corollary}[Second order Abreu equation] \label{CorReducedScalar}
The scalar curvature $s$ on $M^4$ passes to a function on $\Sigma^2$, and
\be
\triangle_\Sigma\mathcal{V}^{\frac12}\,+\,\frac12s\,\mathcal{V}^{\frac12}\;=\;0. \label{EllipticEqnForV}
\ee
\end{corollary}
{\it Proof}.
This follows from $s=-\triangle\log\mathcal{V}$ and Proposition \ref{PropPojectionLaplacian}.
\qed

\begin{proposition}[The $\varphi^i$ elliptic equations] \label{PropDivOfVarphi}
On $\Sigma^2$ we have $d\left(\mathcal{V}^{-\frac12}J_\Sigma{d}\varphi^i\right)=0$.
\end{proposition}
{\it Proof}. If $\Phi:M^4\rightarrow\Sigma^2$ is the moment map projection onto the polytope, then from (\ref{EqnsGJOmegaM}) and (\ref{EqnJSigma}) we obtain $df^i=\Phi^*(\mathcal{V}^{-\frac12}J_\Sigma{}d\varphi^i)$, where the $f^i$ are from Lemma \ref{LemmaHoloSymplCoordRelation}.
Finally $d\Phi^*=\Phi^*d$ and the fact that $\Phi$ is a submersion gives $d(\mathcal{V}^{-\frac12}J_\Sigma{}d\varphi^i)=0$
\qed

{\bf Remark.} When $s=0$, the equations $\mathcal{V}^{\frac12}d\left(\mathcal{V}^{-\frac12}Jd\varphi^i\right)=0$ are just the equations $\varphi^i_{xx}+\varphi^i_{yy}-x^{-1}\varphi^i_x=0$ on the right half-plane.
See Section \ref{SectionAnalysisHalfPlane}.

The following simple theorem illustrates a contrast between the compact case, where $s\le0$ is impossible, and the open case (eg. Theorem \ref{TheoremFlatRTwo}), where $s\le0$ is required, at least at one point.
\begin{corollary} \label{CorCompactFlatPolytope}
If $(\Sigma^2,g_\Sigma)$ is topologically compact, it is impossible that $\triangle_\Sigma\sqrt{\mathcal{V}}\ge0$.
\end{corollary}
{\it Proof}.
The maximum principle gives $\mathcal{V}=const$.
But $\mathcal{V}=0$ on $\partial\Sigma^2$, so $\mathcal{V}\equiv0$.
\qed

\begin{corollary} \label{CorCompactFlatFourManifold}
If $(M^4,J,\omega,\mathcal{X}_1,\mathcal{X}_2)$ is topologically compact, it is impossible that $s\le0$.
\end{corollary}
{\it Proof}.
If $s\le0$ then $\triangle_\Sigma\sqrt{V}\le0$, so the previous corollary applies.
But $\mathcal{V}\equiv0$ means $\mathcal{X}_1$ and $\mathcal{X}_2$ are co-linear throughout, an impossibility.
\qed

\subsection{Reconstruction of the metric on $M^4$ provided $s=0$} \label{SubSectionReconstructionOfTheMetric}

We show how to reconstruct the metrics on $(M^4,J,\omega)$ and $(\Sigma^2,g_\Sigma)$ from knowledge of the functions $\varphi^1$, $\varphi^2$, $\mathcal{V}$ on $\Sigma^2$.
Proposition \ref{PropDivOfVarphi} tells us the coordinates $\varphi^i$ are not harmonic on $(\Sigma^2,g_\Sigma)$; indeed if $\varphi^1$, $\varphi^2$ are harmonic then $g_\Sigma$ and $g$ are both flat metrics; see Lemma \ref{LemmaTConstFlat}.
Instead, when $s=0$, the equation $\triangle\mathcal{V}^{\frac12}=0$ shows $\mathcal{V}^{\frac12}$ and its dual harmonic function constitute a natural harmonic coordinate system.
We denote these by $(x^1,x^2)$, where $x^1$ and $x^2$ are defined as the solutions to
\be
x^1\;=\;\sqrt{\mathcal{V}},\quad\quad dx^2=J_{\Sigma}dx^1. \label{EqnsXYDefRules}
\ee
We shall use the notations $(x,y)$ and $(x^1,x^2)$ interchangeably.
Regarding $\varphi^1$, $\varphi^2$ as functions of $(x,y)$, the formula for $J_\Sigma$ in (\ref{EqnJSigma}) gives
\be
\begin{aligned}
&\frac{1}{x^1}\left<\nabla\varphi^i,\,\nabla\varphi^j\right>\,d\varphi^1\wedge{}d\varphi^2\;=\;-d\varphi^i\wedge{}J_\Sigma{}d\varphi^j.
\end{aligned}
\ee
Changing variables to express both sides in terms of $(x^1,x^2)$ we have
\be
\begin{aligned}
d\varphi^1\wedge{}d\varphi^2
&\;=\;\frac{d\varphi^1}{dx^i}\frac{d\varphi^2}{dx^j}dx^i\wedge{}dx^j
\;=\;\left(
\frac{d\varphi^1}{dx^1}\frac{d\varphi^2}{dx^2}\,-\,\frac{d\varphi^2}{dx^1}\frac{d\varphi^1}{dx^2}
\right)\,dx^1\wedge{}dx^2, \\
-d\varphi^i\wedge{}J_{\Sigma}d\varphi^j
&\;=\;\left(\frac{d\varphi^i}{dx^1}\frac{d\varphi^j}{dx^1}\,+\,
\frac{d\varphi^i}{dx^2}\frac{d\varphi^j}{dx^2}\right)dx^1\wedge{}dx^2,
\end{aligned}
\ee
and so
\be
\begin{aligned}
g_{\Sigma}^{ij}&\;=\;\left<\nabla\varphi^i,\,\nabla\varphi^j\right>\;=\;
x\,\cdot\,\frac{\frac{d\varphi^i}{dx}\frac{d\varphi^j}{dx}\,+\,
\frac{d\varphi^i}{dy}\frac{d\varphi^j}{dy}}
{\frac{d\varphi^1}{dx}\frac{d\varphi^2}{dy}\,-\,\frac{d\varphi^2}{dx}\frac{d\varphi^1}{dy}}.
\end{aligned} \label{EqnFirstSigmaMetricExpresssion}
\ee

We can succinctly express this by letting
$A=
\left(\frac{\partial\varphi^i}{\partial{x}^j}\right)$ and 
$B=
A^{-1}=
\left(\frac{\partial{x}^i}{\partial\varphi^j}\right)$ be the coordinate transition matrices, and writing
\be
\left(g_{\Sigma}^{ij}\right)\;=\;\frac{x}{\det(A)}\,A\,A^T, \quad\quad
\left(g_{\Sigma{}ij}\right)\;=\;\frac{1}{x\,\det(B)}\,B^T\,B.
\ee
Explicitly,
\be
\begin{aligned}
g_\Sigma
&\;=\;\frac{\det(A)}{x}\delta_{kl}\frac{\partial{x}^k}{\partial\varphi^i}\frac{\partial{x}^l}{\partial\varphi^j}\,d\varphi^i\otimes{}d\varphi^j
\;=\;\frac{\det(A)}{x}\delta_{ij}\,dx^i\otimes{}dx^j.
\end{aligned} \label{EqnMetricConstruction}
\ee
So the coordinate transition matrices fully determine $g$ and $g_\Sigma$.
A simple expression for the Gaussian curvature of $(\Sigma^2,g_\Sigma)$ is
\be
\begin{aligned}
K_\Sigma
&\;=\;-\frac{x}{\det(A)}\left(\left(\frac{\partial}{\partial{x}}\right)^2+\left(\frac{\partial}{\partial{y}}\right)^2\right)\log\frac{\det(A)}{x}.
\end{aligned}
\ee
Compare (\ref{EqnCoordSigmaScalarComp}).
Clearly $K_\Sigma$ is invariant under any $SL(2,\mathbb{R})$ recombination of $\varphi^1$, $\varphi^2$.
But transforming $\varphi^1$, $\varphi^2$ by a $GL(2,\mathbb{R})$ combination will usually involve a rescaling of the metric, given by the inverse of the determinant of the transformation.
This will be an important point in \S\ref{SectionQuarterPlaneAndGeneral}.

{\bf Remark}. It is important to point out the scaling properties of (\ref{EqnFirstSigmaMetricExpresssion}).
The right-hand side is scale-invariant in the $\{\varphi^i\}$ functions, whereas $g^{ij}_{\Sigma}$ clearly scales linearly in each; it might be objected that this appears nonsensical.
However the rules we have set up, in (\ref{EqnsXYDefRules}) and (\ref{EqnDefoOfV}), show that if the metric remains unchanged, the $x$, $y$ variables must actually scale quadratically with the $\{\varphi^i\}$, so indeed scale-invariance is retained.
That aside, certainly $(\varphi^1,\varphi^2)$ can be considered functions of $(x,y)$, without reference to any already-existing metric, and from this point of view $(\varphi^1,\varphi^2)$ can obviously be scaled without scaling $(x,y)$.
Via the metric-reconstruction procedure above, the effect is homothetic rescaling of the metric.

{\bf Remark}. For this section, particularly (\ref{EqnMetricConstruction}), see also \cite{Do1}, or \S3 of \cite{AS}.
One can find the K\"ahler potential $u$ from our functions $x$, $y$, $\varphi^1$, $\varphi^2$ by using the recipe of Theorem 1.1 of \cite{Do1}.

\subsection{A conformal change of the metric} \label{SubSecConfChangeOfMetric}

Here we make a second computation of the curvature of $(\Sigma^2,g_\Sigma)$, and show the effects of a certain conformal change of the metric.
In this section only, for convenience, we use both $g$ and $g_\Sigma$ to indicate the metric on $\Sigma$.
But the ``$s$'' in (\ref{EqnDivOfGamma}), (\ref{EqnConfChangeScalar}) and (\ref{EqnConformalScalar}) still refers to the $M^4$ scalar curvature.

The metric on $(\Sigma^2,g)$ has the ``pseudo-K\"ahler'' property
\be
\frac{\partial{g}_{ij}}{\partial\varphi^k}\;=\;\frac{\partial{g}_{ik}}{\partial\varphi^j}. \label{EqnPseudoKahler}
\ee
This was noted in \cite{AS}, and is equivalent to both $[\nabla\varphi^i,\nabla\varphi^j]=0$ and $d(\mathcal{V}^{-\frac12}J_\Sigma{d}\varphi^i)=0$.
The Christoffel symbols are
\be
\begin{aligned}
\Gamma_{ij}^k\;=\;\frac12\frac{\partial{g}_{ij}}{\partial\varphi^s}g^{sk}, \quad\quad
\Gamma^k
\;\triangleq\;g^{ij}\Gamma_{ij}^k
\;=\;-g^{ks}\frac{\partial}{\partial\varphi^s}\log\mathcal{V}^{\frac12}.
\end{aligned} \label{EqnChristoSymbs}
\ee
The divergence of $\Gamma^k$ has a nice interpretation: using $\sqrt{\det(g)}=\mathcal{V}^{-\frac12}$ we get
\be
\begin{aligned}
\frac{\partial\Gamma^k}{\partial\varphi^k}
&\;=\;-\mathcal{V}^{-\frac12}\frac{1}{\sqrt{\det(g)}}\frac{\partial}{\partial\varphi^k}\left(\sqrt{\det(g)}g^{sk}\frac{\partial}{\partial\varphi^s}\mathcal{V}^{\frac12}\right)
\;=\;-\mathcal{V}^{-\frac12}\triangle_\Sigma\mathcal{V}^{\frac12}
 \;=\;\frac12s.
\end{aligned} \label{EqnDivOfGamma}
\ee
The usual formula for scalar curvature in terms of Christoffel symbols is
\be
s_{\Sigma}\;=\;
g^{ij}\frac{\partial\Gamma_{ij}^s}{\partial\varphi^s}
-g^{ij}\frac{\partial\Gamma_{sj}^s}{\partial\varphi^i}
+g_{st}\Gamma^s\Gamma^t
\,-\,g^{is}g^{jt}g_{kl}\Gamma_{ij}^k\Gamma_{st}^l,
\ee
and using $\frac{\partial}{\partial\varphi^s}\left(g^{ij}\Gamma_{ij}^s\right)=\frac{\partial}{\partial\varphi^i}\left(g^{ij}\Gamma_{sj}^s\right)$ this transforms to
\be
\begin{aligned}
s_{\Sigma}
&\;=\;
g^{is}g^{jt}g_{kl}\Gamma_{ij}^k\Gamma_{st}^l-g_{st}\Gamma^s\Gamma^t
\;=\;
\left|\Gamma_{ij}^k\right|^2\,-\,\left|\Gamma^k\right|^2
\;=\;
\left|\Gamma_{ij}^k\right|^2\,-\,\left|\nabla\log\mathcal{V}^{\frac12}\right|^2.
\end{aligned} \label{EqnCoordSigmaScalarComp}
\ee
Now modify the metric by $\widetilde{g}_{\Sigma}=\mathcal{V}^{\frac12}{g}_\Sigma$.
The usual conformal-change formula gives
\be
\begin{aligned}
\widetilde{s}_\Sigma
&\;=\;\mathcal{V}^{-\frac12}\left(s_\Sigma\,-\,\triangle_\Sigma\log\mathcal{V}^{\frac12}\right) \\
&\;=\;\mathcal{V}^{-\frac12}\left(
\left|\Gamma_{ij}^k\right|^2\,-\,\left|\nabla\log\mathcal{V}^{\frac12}\right|^2
\,-\,\triangle_\Sigma\log\mathcal{V}^{\frac12}\right).
\end{aligned} \label{EqnConfChangeScalar}
\ee
But $\triangle_\Sigma\log\mathcal{V}^{\frac12}+\left|\nabla\log\mathcal{V}^{\frac12}\right|^2=\mathcal{V}^{-\frac12}\triangle_\Sigma\mathcal{V}^{\frac12}$ so the conformally related scalar is simply
\be
\begin{array}{ll}
\widetilde{s}_\Sigma
&\;=\;\mathcal{V}^{-\frac12}\left(\left|\Gamma_{ij}^k\right|^2\,-\,\mathcal{V}^{-\frac12}\triangle_\Sigma\mathcal{V}^{\frac12}\right) \\
&\;=\;\mathcal{V}^{-\frac12}\left(\left|\Gamma_{ij}^k\right|^2\,+\,\frac12s\right).
\end{array} \label{EqnConformalScalar}
\ee
Therefore if $s\ge0$, then $\widetilde{s}_\Sigma\ge0$.
This will be used decisively in the proof Theorem \ref{TheoremFlatRTwo}.

\subsection{Global behavior of the harmonic coordinate system} \label{SubSecGlobalHarmonicCoordBehavior}

Fundamental to our paper is the use of the harmonic coordinates $(x,y)$ as defined in (\ref{EqnsXYDefRules}).
But we need some a priori facts about this system, stemming from either the geometrical situation or from properties (A)-(E), in order to utilize our analytic theorems of \S\ref{SectionAnalysisHalfPlane}.
Specifically, we need to know that $(x,y)$ is actually a global coordinate system, and doesn't have any branch points.

\begin{proposition} \label{PropXISOneOneOnto}
Assume the polytope $(\Sigma^2,g_\Sigma)$ is closed and the boundary is connected.
Then map $X=(x,y)^T$ is $1$-to-$1$ and onto the right half-plane.
\end{proposition}
{\it Proof}.
The function $z=x+\sqrt{-1}y$ is a holomorphic function on $\Sigma^2$, and maps the right-half plane $\{\varphi^1\ge0\}$ into the right half-plane $\{x\ge0\}$.
From ({\it{ix}}) below, we have $y=\varphi^2+O(r^2)$, so that the line $\{x=0\}$ is indeed in the image of $z$.

\underline{\it Definition of pseudoboundary point}. 
We define the notion of a pseudo-boundary point: a point $X_0=(x_0,y_0)$ is a pseudoboundary point if the following holds: given any neighborhood $U\subset\mathbb{R}^2$ of $X_0$, there is some component $U_0$ of $X^{-1}(U)\subseteq\Sigma^2$ such that $X(U_0)\ne{U}$
(this essentially means $U$ is not evenly covered, excluding the phenomenon of ramification).
Letting $U\in\mathbb{R}^2$ be any neighborhood of a point $(0,y)$ of the actual boundary $\{x=0\}$, we see the boundary locus $\{x=0\}$ is contained in the pseudoboundary locus.

\underline{\it All pseudoboundary points lie on $\{x=0\}$}.
Let $X_0$ be a pseudoboundary point, and let $U$ be a neighborhood that is not evenly covered, and let $U_0\subseteq{X}^{-1}(U)$ be a component that does not completely cover $U$.
We can assume that $U$ itself is precompact, and that $X:U_0\rightarrow{U}$ is 1-1 but not onto.
We can assume this because the pre-image of $X_0$ in $U_0$ cannot be a branch point, or else $U_0\rightarrow{U}$ would certainly be onto.
By shrinking $U$ if necessary, we can assume $U_0$ it contains no branch points at all.

The pre-image $U_0\subset\Sigma^2$ is either pre-compact or not.
If it is precompact, then it must intersect the boundary of $\Sigma^2$.
This is because $U_0$ has no critical points, so either $U_0$ intersects the boundary of $\Sigma^2$ or else $X$ is a diffeomorphism in a neighborhood of $U_0$, meaning $U_0$ must evenly cover $U$, an impossibility.
Thus $X_0$ is actually in $\{x=0\}$, and so is actually part of the boundary.

Now assume the closure $\overline{U_0}$ is not compact.
All points on $\partial\left(X(U_0)\right)\setminus\partial{U}$ are in fact pseudoboundary points.
Thus $X_0$ is not an isolated pseudoboundary point.

Because $U_0$ is not compact, it extends to infinity in $(\varphi^1,\varphi^2)$ coordinates.
Then consider the pushforwards of $\left.\varphi^i\right|_{U_0}$ to $U$.
We see that $\varphi^1$ is infinite at all of the pseudoboundary points on $\partial{U}_0$.
But the $\varphi^i$ satisfy the equation $\varphi^i_{xx}-x^{-1}\varphi^i_x+\varphi^i_{yy}=0$, which is uniformly elliptic away from $\{x=0\}$, and therefore each has isolated singularities away from $\{x=0\}$.
Thus, again, we see that the pseudoboundary points $\partial\left(X(U_0)\right)\setminus\partial{U}$ must lie on the $\{x=0\}$ locus.

\underline{\it Pre-images of domains containing points of $\{x=0\}$ are pre-compact.}

Let $X_0\in\{x=0\}$, and let $U$ be a sufficiently small neighborhood around $(0,0)$.
Consider the ``barrier function''
\be
\varphi_A(x,y)\;=\;A\left(1-y^2+x^2\log(x) - x^2\right), \quad {\rm where} \quad A>0.
\ee
This solves the equation $x(x^{-1}\varphi_x)_x+\varphi_{yy}=0$.
Let $\Omega$ be the component of $\{\varphi>0\}\cap\{x\ge0\}$ that contains the point $(0,0)$.

Pulling $\varphi$ back along $X:\Sigma^2\rightarrow\mathbb{R}^2$, we obtain a function on $\Sigma^2$, still denoted $\varphi$.
On $\Sigma^2$ this function still satisfies $d\left(\mathcal{V}^{-\frac12}J_\Sigma{d}\varphi\right)=0$, which is an elliptic equation (see Proposition \ref{PropDivOfVarphi} above and Equation (\ref{EqnOrigUnmodifiedPDE}) from Section \ref{SectionAnalysisHalfPlane}).

Consider the domains $\Omega_\epsilon=\Omega\cap\{x\ge\epsilon\}$, and the images $X^{-1}(\Omega_\epsilon)$ of these domains on $\Sigma^2$.
We have that $\varphi=0$ on $X^{-1}(\Omega_\epsilon\setminus\{x=\epsilon\})$, so therefore $\varphi_A<\varphi_i$ on this boundary component.
As $\epsilon\rightarrow0$, we have that $\varphi^i\rightarrow\infty$, because the $\{x=0\}$ locus occurs at coordinate infinity.
Thus, given any $A$, there is some $\overline{\epsilon}$ so that $\epsilon<\overline\epsilon$ implies $\varphi_A<\varphi_i$ on $\Omega_\epsilon$.
Thus $\varphi_A<\varphi^i$ on all of $X^{-1}(\Omega)$.
Sending $A\rightarrow\infty$ implies that $\varphi^i=\infty$ on all of $X^{-1}(\Omega)$, an impossibility.

\underline{\it The pre-image of $\{x=0\}$ has a single component.}
The pre-image of $\{x=0\}$ is the boundary of the polytope: this follows trivially from the fact that $x\ne0$ on the interior of $\Sigma^2$, and from the fact, proved in the paragraph above, that points at infinite do not map to $\{x=0\}$.
Conversely, because it is assumed that the polytope contains all of its edges, the boundary of the polytope is precisely the locus on $\Sigma^2$ on which $\{x=0\}$.
By Proposition \ref{PropCTwoAndLipschitz} below (or by (D)), the function $X:\Sigma^2\rightarrow\mathbb{R}^2$ is Lipschitz, therefore continuous.
Thus, if $X$ is $k$-to-$1$ on $\{x=0\}$, the pre-image of $\{x=0\}$ must have $k$ components.

\underline{\it The map $(x,y):\Sigma^2\rightarrow\mathbb{R}^2$ is one-to-one.}
If $X=(x,y)^T$ is not one-to-one, then there must be pseudoboundary points, or else $x$ is not one-to-one as a map onto $\{x=0\}$.
Both of these possibilities were disproven above.
\qed

{\bf Remark}. For an example of how Proposition \ref{PropXISOneOneOnto} fails if the polytope's boundary is disconnected, see Example 8.

{\bf Remark.} The preceding theorem is far from a general statement on non-negative harmonic functions on manifolds.
Of course functions on virtually all manifolds must have critical points.
Absolutely critical to Proposition \ref{PropXISOneOneOnto} here is the interplay between the harmonic functions and the moment functions.

\subsection{Coordinate behavior near Polytope Edges} \label{SubSecPolyTopeEdges}

Of course all functions $\varphi^1$, $\varphi^2$, $x=\sqrt{\mathcal{V}}$, and $y$ are themselves generically $C^\infty$, but at $\{x=0\}$ it is possible that $\varphi^1,\varphi^2$ are not $C^\infty$ with respect to the $(x,y)$ system.

In this subsection, we prove the essential facts that near polytope edges, the functions $\varphi^1$, $\varphi^2$ are of class at least $C^2$ with respect to the variables $(x,y)$, and that near vertices, these functions are Lipschitz.
Of course the two coordinate systems $(\varphi^1,\varphi^2)$ and $(x,y)$ are $C^\infty$, but a priori their inverses might not be, so $\varphi^i$ might not even be $C^1$ as measured in $(x,y)$ coordinates.

Via the momentum reduction, we have the metric polytope $(\Sigma^2,g_\Sigma)$.
For the purposes of this section, define $r$ to be the $g_\Sigma$-distance from the $\{x=0\}$ locus.
This lifts to $M^4$, where $r$ is still a distance function.
Note that $\nabla{r}$ is well-defined on $\{x=0\}$ on $\Sigma^2$, but not at $\{x=0\}$ on $M^4$.

We focus our attention on a segment of the polytope, and explore the behavior of $x$, $y$, and $\varphi^1$, and $\varphi^2$ near that segment.
In what follows, we may either assume the polytope satisfies the conditions (A)-(E), or else that the polytope metric reduction of some $(M^4,J,\omega,\mathcal{X}_1,\mathcal{X}_2)$.

\begin{lemma}
Assume $(\Sigma^2,g_\Sigma)$ comes from a K\"ahler reduction of some $(M^4,J,\omega,\mathcal{X}_1,\mathcal{X}_2)$.
After picking any boundary segment and possibly re-combining the $\varphi^1$, $\varphi^2$, assume $\nabla\varphi^2$ is parallel to that segment.
Then $\frac{\nabla\varphi^2}{|\nabla\varphi^2|}$ is convariant-constant.
\end{lemma}
{\it Proof}.
We are trying to show $\nabla\frac{\nabla\varphi^2}{|\nabla\varphi^2|}=0$.
First note that $\left<\nabla{r},\nabla\varphi^2\right>=0$ on the boundary segment, and that
\be
\left<\nabla_{\nabla\varphi^2}\nabla\varphi^2,\,\nabla{r}\right>
=\left<\nabla\varphi^2,\,\nabla_{\nabla\varphi^2}\nabla{r}\right>\;=\;0
\ee
by total geodesy.
This easily implies $\nabla_{\nabla\varphi^2}\frac{\nabla\varphi^2}{|\nabla\varphi^2|}=0$.
To evaluate $\nabla_{\nabla{r}}\frac{\nabla\varphi^2}{|\nabla\varphi^2|}$, we lift the situation to $M^4$, and compute
\be
\left<\nabla_{\nabla{r}}\nabla\varphi^2,\,\nabla{r}\right>
=\left<\nabla_{\nabla{r}}\mathcal{X}_2,\,J\nabla{r}\right>
=\left<\nabla_{\mathcal{X}_2}\nabla{r},\,J\nabla{r}\right>=0
\ee
again by total geodesy.
Once again, this easily implies $\nabla_{\nabla\varphi^2}\frac{\nabla\varphi^2}{|\nabla\varphi^2|}=0$.
\qed

\begin{lemma}[Pseudo-Jacobi identity for the $\nabla\varphi^i$] \label{LemmaPseudoJacobi}
If $\gamma$ is a geodesic on $(\Sigma^2,g_\Sigma)$ then $\nabla\varphi^i$ obeys a second order equation along $\gamma$ given by $\left<\nabla_{\dot\gamma}\nabla_{\dot\gamma}\nabla\varphi^i,\nabla\varphi^j\right>=\frac12\dot\gamma\dot\gamma\left<\nabla\varphi^i,\nabla\varphi^j\right>-\left<\nabla_{\dot\gamma}\varphi^i,\nabla_{\dot\gamma}\varphi^j\right>$.
In the case that $(\Sigma^2,\gamma_\Sigma)$ is a K\"ahler reduction of some $(M^4,J,\omega,\mathcal{X}_1,\mathcal{X}_2)$ then
\be
\nabla_{\dot\gamma}\nabla_{\dot\gamma}\nabla\varphi^i\,+\,\Phi_*J\left(\Riem(J\nabla\varphi^i,\,\Phi_*^{-1}\dot\gamma)\Phi_*^{-1}\dot\gamma\right)\;=\;0
\ee
where $\Phi:M^4\rightarrow\Sigma^2$ is the reduction map, and $\Phi_*^{-1}$ indicates a choice of horizontal lift.
\end{lemma}
{\it Proof}.
The first claim is obvious.
The second claim follows from noting that, while $[\nabla\varphi^i,\dot\gamma]$ might not equal zero, certainly $[\mathcal{X}_i,\dot\gamma]=0$, so the derivation of the Jacobi equation can proceed as normal after multiplication of $\nabla_{\dot\gamma}\nabla_{\dot\gamma}\varphi^i$ by $J$.
\qed

Our condition (D) assumes also $\nabla\frac{\nabla\varphi^1}{|\nabla\varphi^1|}=0$.
This is actually a consequence of ({\it{vii}}) below, so we do not prove it separately.

We prove the following:
\begin{itemize}
\item[{\it{i}})] The $\{x=0\}$ locus is totally geodesic on both $\Sigma^2$ and $M^4$ (this is essentially obvious).
\item[{\it{ii}})] We may make an $SL(2,\mathbb{R})$-recombination of $\varphi^1$, $\varphi^2$, followed by possible addition of a constant, so that $\varphi^1=0$ on the segment under consideration, while leaving $x=\sqrt{\mathcal{V}}$ unchanged.
\item[{\it{iii}})] $[\nabla{r},\nabla\varphi^2]=0$ at $\{x=0\}$ on $\Sigma^2$ and $\left<\nabla{r},\nabla\varphi^2\right>=O(r^2)$.
\item[{\it{iv}})] In the $(r,\varphi^2)$-coordinate system on $\Sigma^2$, we have Taylor series $\varphi^1=C(\varphi^2)r^2\,+\,O(r^4)$.
\item[{\it{v}})] On each segment of the polytope, we have $C(\varphi^2)=const$; we take $C=\frac12$.
\item[{\it{vi}})] We have $\left<\nabla\varphi^1,\nabla\varphi^2\right>=O(r^4)$.
\item[{\it{vii}})] $|\nabla\varphi^1|^2=r^2+O(r^5)$ and $|\nabla\varphi^2|^2=C_1(\varphi^2)(1+2Kr^2)+O(r^3)$, where $K=K(\varphi^2)$ is the $g_\Sigma$-curvature of the polytope at the point $(r,\varphi^2)=(0,\varphi^2)$.
Here, $C_1(\varphi^2)$ is simply $|\nabla\varphi^2|^2$ evaluated at $(r,\varphi^2)=(0,\varphi^2)$.
\item[{\it{viii}})] $x=\sqrt{\mathcal{V}}=\sqrt{C_1}r\left(1\,+\,Kr^2\right)\,+\,O(r^4)$
\item[{\it{ix}})] $y=\varphi^2+O(r^2)$
\item[{\it{x}})] $\varphi^1$, $\varphi^2$ are at least $C^2$ on segments and are Lipschitz everywhere in the $(x,y)$-coordinate plane
\end{itemize}

\underline{\it Proof of (i)}:
$x=\sqrt{|\nabla\varphi^1|^2|\nabla\varphi^2|^2-\left<\nabla\varphi^1,\nabla\varphi^2\right>^2}$, so $\{x=0\}$ precisely when $\{\nabla\varphi^1,\nabla\varphi^2\}$ is a dependent set, so there are constants with $c_1\mathcal{X}_1+c_2\mathcal{X}_2=0$ on $M^4$.
But a combination of Killing fields is Killing, and has totally geodesic zero locus.
Thus $\{x=0\}$ is totally geodesic.

\underline{\it Proof of (ii)}:
An $SL(2,\mathbb{R})$-recombination of $\{\varphi^1,\varphi^2\}$ gives rise to an $SL(2,\mathbb{R})$-recombination of $\nabla\varphi^1,\nabla\varphi^2$, which obviously leaves $x=\sqrt{\mathcal{V}}$ unchanged.
If $\nabla\varphi^1,\nabla\varphi^2$ are co-linear at a point, we can certainly make such a transformation to obtain $\nabla\varphi^1=0$ at that point.
But the total geodesy of the segment and the fact that the $\nabla\varphi^i$ are Jacobi fields, we have that $\nabla\varphi^1=0$ on the entire segment.
Addition a constant, we can certainly make $\varphi^1=0$ on that segment as well.

\underline{\it Proof of (iii)}:
On $\Sigma^2$, the gradient $\nabla\varphi^2$ is parallel to the polytope boundary segment.
By total geodesy (vanishing of the second fundamental form), we have $\nabla_{\nabla\varphi^2}\nabla{r}=0$ so $[\nabla{r},\nabla\varphi^2]=-\nabla_{\nabla{r}}\nabla\varphi^2$.
In case we already know that condition (D) is satisfied, this is already zero.
In case $(\Sigma^2,g_\Sigma)$ is the reduction of some $M^4$, so we perhaps do know (D) apriori.
But we can lift back to $M^4$ (choosing a horizontal lift, of course), where we can see (using $J$-invariance and taking derivatives) that $\left<\nabla{r},\nabla\varphi^2\right>=0$ and $\nabla{r}\left<\nabla{r},\nabla\varphi^2\right>=\left<\nabla{r},\nabla_{\nabla{r}}\nabla\varphi^2\right>=0$ as we limit to $\{x=0\}$.

\underline{\it Proof of (iv)}:
On $M^4$ we have that $\mathcal{X}_1$ is a Jacobi field along trajectories of $\nabla{r}$ with an initial condition of $\mathcal{X}_1=0$.
Then $\nabla\varphi^1$ is also Jacobi, with initial conditions $\nabla\varphi^1=0$ and $\nabla_{\nabla{r}}\nabla\varphi^1\propto\nabla{r}$.
With $\varphi^1=0$ on $r=0$, this is sufficient to establish the claim.

\underline{\it Proof of (v)}: Using $\varphi^1=C(\varphi^2)r^2+O(r^3)$ we have
\be
\begin{aligned}
0&\;=\;[\nabla\varphi^1,\nabla\varphi^2] \\
&\;=\;[C'r^2\nabla\varphi^2+2Cr\nabla{r}+O(r^2)\nabla{r}+O(r^3)\nabla\varphi^2,\,\nabla\varphi^2] \\
&\;=\;C''r^2\nabla\varphi^2+2C'r\left<\nabla\varphi^2,\nabla{r}\right>\nabla\varphi^2
+2C'r\nabla{r}+2C\left<\nabla\varphi^2,\nabla{r}\right>\nabla{r}+3Cr[\nabla{r},\nabla\varphi^2]+O(r^2)
\end{aligned}
\ee
Using that $[\nabla{r},\nabla\varphi^2]=O(r)$ and $\left<\nabla\varphi^2,\nabla{r}\right>=O(r^2)$ we have the Taylor series
\be
0\;=\;2C'(\varphi^2)r\nabla{r}\,+\,O(r^2)
\ee
and so therefore $C'=0$.

\underline{\it Proof of (vi)}: Of course $\left<\nabla\varphi^1,\nabla\varphi^2\right>=0$ on the boundary so the inner product is at least $O(r)$.
We take one, two, and three derivatives along $\nabla{r}$, and show that each is at least $O(r)$.
Taking one derivative
\be
\nabla{r}\left<\nabla\varphi^1,\nabla\varphi^2\right>\;=\;\left<\nabla_{\nabla\varphi^1}\nabla\varphi^2,\,\nabla{r}\right>
\ee
which is $O(r)$ because $\nabla\varphi^1=Cr\nabla{r}+O(r^2)$, and using ({\it{iii}}) above.
Taking a second derivative,
\be
\begin{aligned}
&\frac12\nabla{r}\nabla{r}\left<\nabla\varphi^1,\nabla\varphi^2\right>
\;=\;
\left<\nabla_{\nabla{r}}\nabla_{\nabla\varphi^1}\nabla\varphi^2,\,\nabla{r}\right> \\
&\quad\quad\;=\;
\left<\nabla_{\nabla\varphi^1}\nabla_{\nabla{r}}\nabla\varphi^2,\,\nabla{r}\right>
\,+\,\left<\nabla_{[\nabla{r},\nabla\varphi^1]}\nabla\varphi^2,\,\nabla{r}\right>
\,+\,\left<\Riem(\nabla{r},\,\nabla\varphi^1)\nabla\varphi^2,\,\nabla{r}\right>
\end{aligned}
\ee
is zero because $\nabla\varphi^1=O(r)$ and $[\nabla{r},\nabla\varphi^1]=C\nabla{r}+O(r)$ and $\left<\nabla_{\nabla{r}}\nabla\varphi^2,\nabla{r}\right>=O(r)$.
Finally we take a third derivative to obtain
\be
\begin{aligned}
&\frac12\nabla{r}\nabla{r}\nabla{r}\left<\nabla\varphi^1,\nabla\varphi^2\right> \\
&\quad\;=\;
\nabla{r}\left<\nabla_{\nabla\varphi^1}\nabla_{\nabla{r}}\nabla\varphi^2,\,\nabla{r}\right>
\,+\,\nabla{r}\left<\nabla_{[\nabla{r},\nabla\varphi^1]}\nabla\varphi^2,\,\nabla{r}\right>
\,+\,\nabla{r}\left<\Riem(\nabla{r},\,\nabla\varphi^1)\nabla\varphi^2,\,\nabla{r}\right>
\end{aligned}
\ee
The first term can be re-written
\be
\begin{aligned}
&
\nabla{r}\left<\nabla_{\nabla\varphi^1}\nabla_{\nabla{r}}\nabla\varphi^2,\,\nabla{r}\right> \\
&\quad=\;
\nabla{r}\nabla\varphi^1\left<\nabla_{\nabla{r}}\nabla\varphi^2,\,\nabla{r}\right>
-\nabla{r}\left<\nabla_{\nabla{r}}\nabla\varphi^2,\,\nabla_{\nabla\varphi^1}\nabla{r}\right>  \\
&\quad=\;
\nabla\varphi^1\nabla{r}\left<\nabla_{\nabla{r}}\nabla\varphi^2,\,\nabla{r}\right>
+[\nabla{r},\nabla\varphi^1]\left<\nabla_{\nabla{r}}\nabla\varphi^2,\,\nabla{r}\right>
-\nabla{r}\left<\nabla_{\nabla{r}}\nabla\varphi^2,\,\nabla_{\nabla\varphi^1}\nabla{r}\right> 
\end{aligned}
\ee
Of these three terms, the first is zero because $\nabla\varphi^1=O(r)$.
The second is zero because $[\nabla{r},\nabla\varphi^1]=C\nabla{r}+O(r)$ and $\nabla\varphi^2$ satisfies the Jacobi equation $\nabla_{\nabla{r}}\nabla_{\nabla{r}}\nabla\varphi^2+\Riem(\nabla\varphi^2,\nabla{r})\nabla{r}=0$.
The third term is zero for this same reason, along with using a Jacobi equation for $\nabla\varphi^1$.

The second term is zero because $[\nabla{r},\nabla\varphi^1]=C\nabla{r}+O(r)$, and because $[\nabla{r},[\nabla{r},\nabla\varphi^1]]$ is proportional to $\nabla\varphi^2$ on $\{x=0\}$, and because $\nabla\varphi^2$ is totally geodesic on $\{x=0\}$ and satisfies a Jacobi equation.

The third term is zero simply by noting $\nabla_{\nabla{r}}\nabla{r}\equiv0$, by $\nabla\varphi^1=Cr\nabla{r}+O(r^2)$, and using the anti-symmetry of the Riemann tensor in the first two positions.

\underline{\it Proof of (vii)}: With $C=\frac12$ we have $\nabla\varphi^1=r\nabla{r}+O(r^2)$ so $|\nabla\varphi^1|^2=r^2+O(r^3)$.
To see that $O(r^3)$ can be improved to $O(r^5)$, take two more derivatives, and use the Jacobi equation for $\nabla\varphi^1$.
Then we look at the Taylor series for $|\nabla\varphi^2|^2$.
The $r^0$-term is, of course, given by some $C_1(\varphi^2)$.
To see there is no $r^1$-term, note that $\frac12\nabla{r}|\nabla\varphi^2|^2=\left<\nabla_{\nabla{r}}\nabla\varphi^2,\nabla\varphi^2\right>=0$ by the proof of ({\it{iii}}).
Taking two derivatives, we see that
\be
\begin{aligned}
\frac12\nabla{r}\nabla{r}|\nabla\varphi^2|^2
&\;=\;\left<\nabla_{\nabla{r}}\nabla_{\nabla{r}}\nabla\varphi^2,\,\nabla\varphi^2\right>\,+\,|\nabla_{\nabla{r}}\nabla\varphi^2|^2
\end{aligned}
\ee
We have that $\nabla_{\nabla{r}}\nabla\varphi^2=0$ on $\{x=0\}$ and that $\left<\nabla_{\nabla{r}}\nabla_{\nabla{r}}\nabla\varphi^2,\,\nabla\varphi^2\right>=\tilde{K}|\nabla\varphi^2|^2$, where $\widetilde{K}$ is a bisectional curvature on $M^4$ (by Lemma \ref{LemmaPseudoJacobi}), or else comes from second order data on the metric of $(\Sigma^2,g_\Sigma)$ if we don't know $\Sigma^2$ is a K\"ahler reduction.
We therefore obtain that the second-order term indeed has coefficient $2C_1(\varphi^2)\widetilde{K}$.

\underline{\it Proof of (viii)}:
Using our Taylor expansions for $|\nabla\varphi^1|^2$, $|\nabla\varphi^2|^2$ and using that $\left<\nabla\varphi^1,\nabla\varphi^2\right>=O(r^4)$, it is very simple to compute
\be
x\;=\;\sqrt{\mathcal{V}}\;=\;\sqrt{C_1}\cdot{}r\cdot\left(1\,+\,\widetilde{K}r^2\right)\,+\,O(r^4).
\ee

\underline{\it Proof of (ix)}: Using that $\nabla{x}=\sqrt{C_1}\nabla{r}+O(r^2)$, we have that $\nabla{y}=J\nabla{x}=\sqrt{C_1}\nabla{r}+O(r^2)$.
Of course $\sqrt{C_1(\varphi^2)}=\left.|\nabla\varphi^2|\right|_{(r,\varphi^2)=(0,\varphi^2)}$, so $|\nabla{y}|=|\nabla\varphi^2|$ along $\{x=0\}$.
Since $\nabla{y}$ is parallel to $\nabla\varphi^2$, we can take $y=\varphi^2$ along $\{x=0\}$.
This justifies $y=\varphi^2+O(r^2)$.

We re-express assertion ({\it{x}}) as a proposition.
\begin{proposition} \label{PropCTwoAndLipschitz}
Assume $(\Sigma^2,g_\Sigma)$ is a closed polytope that either obeys (A)-(F) or else is the reduction of $(M^4,J,\omega,\mathcal{X}_1,\mathcal{X}_2)$.
As measured in the $(x,y)$-coordinate system, the functions $\varphi^1$, $\varphi^2$ are at least $C^2$ along any polytope edge, and are Lipschitz everywhere.
\end{proposition}
{\it Proof}.
The series expansions of $x$ and $y$ from the proofs of ({\it{viii}}) and ({\it{ix}}) easily imply that
\be
\begin{aligned}
&\frac{\partial}{\partial{x}}=\frac{\partial}{\partial{r}}\,+\,O(r^2) \\
&\frac{\partial}{\partial{y}}=\frac{\partial}{\partial\varphi^2}\,+\,O(r)\frac{\partial}{\partial{r}}\,+\,O(r^2)
\end{aligned}
\ee
so that
\be
\begin{aligned}
&\frac{\partial^2}{\partial{x}^2}=\frac{\partial^2}{\partial{r}^2}\,+\,O(r) \\
&\frac{\partial^2}{\partial{y}^2}=\left(\frac{\partial}{\partial\varphi^2}\right)^2\,+\,O(1)\frac{\partial}{\partial{r}}\,+\,O(r) \\
\end{aligned}
\ee
This implies $\varphi^1$, $\varphi^2$ are at least $C^2$ on edges.
(Note that we cannot immediately assert $C^3$ due to the fact that the ``$O(1)$'' term might conceivably not have a well-defined derivative in the $\varphi^2$ direction.)

To consider the case on corners, note that the function $r$ is Lipschitz, so $x=\sqrt{C_1}r+O(r^2)$ is also Lipschitz, and so is $y$.
Thus the map ${{X}}:M^4\rightarrow\mathbb{R}^2$ given by ${{X}}(p)=(x(p),y(p))^T$ is a Lipschitz map.
Therefore pushing forward the invariant $C^\infty$ functions $\varphi^1,\varphi^2$ along ${{X}}$ does indeed give Lipschitz functions on the $(x,y)$-plane.
\qed

{\bf Remark}.
A key part of the above theorem is that the locus $\{x=0\}$ actually exists on $M^4$ or on $\Sigma^2$ and so no edges that are ``infintiely far away.''
Proposition \ref{PropCTwoAndLipschitz}, which asserts Lipschitzness of $\varphi^1$, $\varphi^2$ on the $\{x=0\}$ locus, will not apply to polytope edges that are infinitely far away; see Example 7.
The prototype is the standard matric on $M^4=\mathbb{S}^2\times\mathbb{H}^2$ where $\mathbb{H}^2$ is the pseudosphere with constant curvature $-1$.
This has the polytope given by a half-strip $\{-1\le\varphi^1\le1\}\cap\{\varphi^2>1\}$ in the $(\varphi^1,\varphi^2)$-plane, whose vertical edge rays {\it are} part of the polytope, but whose vertical edge segment {\it is not} part of the polytope.

\section{The case that the polytope is $\mathbb{R}^2$} \label{SectionRTwoCase}

Here we prove Theorem \ref{TheoremFlatRTwo}, namely that if $\triangle_\Sigma\sqrt{\mathcal{V}}\le0$ and the polytope $\Sigma^2$ is metrically complete (not necessarily coordinate-complete), then it is flat.
The outline of the proof is as follows.
First we prove that if $\mathcal{V}$ is a constant function, then $\widetilde{g}_\Sigma$ is flat.
Second, we show that if $(\Sigma^2,g_\Sigma)$ is complete, the ``canonical'' conformal change from \S\ref{SubSecConfChangeOfMetric} gives a metric $\widetilde{g}_\Sigma$ that in fact remains complete.
Using this along with $\triangle_\Sigma\sqrt{\mathcal{V}}\le0$ and some basic elliptic theory, we prove $\mathcal{V}$ is indeed constant, completing the proof.

\begin{lemma} \label{LemmaTConstFlat}
Assume $(\Sigma^2,g_\Sigma)$ is metrically metric.
If $\mathcal{V}$ is constant, or else if $\triangle_\Sigma\sqrt{\mathcal{V}}\le0$ and $(\Sigma^2,J_\Sigma)$ is biholomorphic to $\mathbb{C}$, then $(\Sigma^2,g_\Sigma)$ is flat; indeed $g_\Sigma$ has constant coefficients when expressed in $\varphi^1$, $\varphi^2$ coordinates.
\end{lemma}
{\it Proof}.
First assume $\mathcal{V}^{\frac12}$ is constant.
Then the expression $d(\mathcal{V}^{-\frac12}J_\Sigma{d}\varphi^k)=0$ from Proposition \ref{PropDivOfVarphi} reduces to $dJ_\Sigma{}d\varphi^k=0$.
Thus $\varphi^1$ and $\varphi^2$ are harmonic, so each determines a complex coordinate: either $z=\varphi^1+\sqrt{-1}\eta^1$ or $w=\varphi^2+\sqrt{-1}\eta^2$ (where $d\eta^k=-J_\Sigma{}d\varphi^k$).
Then
\be
\begin{aligned}
\frac{d}{dz}&\;=\;\frac{1}{2|\nabla\varphi^1|^2}\left(\nabla\varphi^1\,+\,\sqrt{-1}J_\Sigma\nabla\varphi^1\right) \\
\frac{d}{dw}&\;=\;\frac{1}{2|\nabla\varphi^2|^2}\left(\nabla\varphi^2\,+\,\sqrt{-1}J_\Sigma\nabla\varphi^2\right)
\end{aligned}
\ee
and we easily compute the transition function:
\be
\frac{dw}{dz}
\;=\;g\left(\frac{d}{dz},\,\overline{\nabla}{w}\right)
\;=\;\frac{\left<\nabla\varphi^1,\,\nabla\varphi^2\right>}{|\nabla\varphi^1|^2}\,-\,\sqrt{-1}\frac{\sqrt\mathcal{V}}{|\nabla\varphi^1|^2}.
\ee
But $\frac{dw}{dz}$ is holomorphic, so its real and imaginary parts are harmonic.
Therefore $|\nabla\varphi^1|^{-2}=-\mathcal{V}^{-\frac12}Im(\frac{dw}{dz})$ is harmonic.
In the $z$-coordinate the Hermitian metric is $h_\Sigma=\left|\frac{d}{dz}\right|^2=\frac12|\nabla\varphi^1|^{-2}$, and therefore the curvature
\be
K_\Sigma\;=\;-h^{-1}\triangle_\Sigma\log{h}_\Sigma\;=\;8|\nabla|\nabla\varphi^1||^2
\ee
is non-negative, forcing the complete manifold $(\Sigma^2,g_\Sigma)$ to be parabolic (this is due to the Cheng-Yau condition for parabolicity; see \cite{CY} or the remark below).
Equivalently the complete manifold $(\Sigma^2,J_{\Sigma})$ is biholomorphic to $\mathbb{C}$.
But then the Liouville theorem implies $|\nabla\varphi^1|^{-2}$ is actually constant, as it is an harmonic function on $\mathbb{C}$ bounded from below.
Similarly $|\nabla\varphi^1|^{-2}$ and $\left<\nabla\varphi^1,\nabla\varphi^2\right>$ are constant.
Expression (\ref{EqnSigmaMetric}) implies $g_\Sigma$ are constant in $(\varphi^1,\varphi^2)$ coordinates, and therefore have zero curvature.

For the second assertion, assume $(\Sigma^2,J_\Sigma)$ is biholomorphic to $\mathbb{C}$ and $\triangle_\Sigma\sqrt{\mathcal{V}}\le0$.
Since the function $\mathcal{V}^{\frac12}$ is superharmonic and bounded from below, it is constant, and we may use the argument from above.
\qed

\begin{lemma} \label{LemmaConformalComplete}
Assume $(\Sigma^2,g_\Sigma)$ is complete and $\triangle_\Sigma\sqrt{\mathcal{V}}\le0$.
Setting $\widetilde{g}_\Sigma=\mathcal{V}^{\frac12}g_\Sigma$, then $(\Sigma^2,\widetilde{g}_\Sigma)$ is also complete.
\end{lemma}
{\it Proof}.
See \S\ref{SubSecConfChangeOfMetric} for some basic facts on the conformal change $\widetilde{g}_\Sigma=\mathcal{V}^{\frac12}g_\Sigma$.
Assume $\gamma:[0,R]\rightarrow\Sigma^2$ is a geodeisc in the $\widetilde{g}_\Sigma$ metric that gives a shortest distance to $\partial\Sigma^2$.
For a contradiction, we show that if $\gamma$ has finite length in $\widetilde{g}_\Sigma$, it also has finite length in $g_\Sigma$.
Set $p=\gamma(0)$ and $r=\dist(p,\cdot)$.
We are free to re-choose $p$ along $\gamma$ to make $R$ smaller if convenient.

As a first step we verify that, for $R$ sufficiently small, we have $\mathcal{V}^{\frac12}$ is strictly monotone along $\gamma$; in particular $\left|\frac{\partial}{\partial{r}}\left(T^{\frac12}\circ\gamma\right)\right|\ge\delta$ for some $\delta>0$.
The standard method (eg. Lemma 3.4 of \cite{GT}) is the ``Hopf Lemma'' applied to the equation $\widetilde\triangle_\Sigma{\mathcal{V}}^{\frac12}\le0$, but in our situation there is an issue concerning how the elliptic operator itself might degenerate near the boundary.

To emulate the standard Hopf lemma proof we must first use Laplacian comparison \cite{SY}.
On the annulus $\mathcal{A}=B_p(R)\setminus{B}_p(R/2)$ we have $\mathcal{V}>0$ except possibly on $\partial{B}_p(R)$, so we have the existence of some $\epsilon>0$ with $\epsilon\left(r^{-1}-R^{-1}\right)\le\mathcal{V}^\frac12$ on $\partial\mathcal{A}$.
On the other hand, since $\tilde{s}_\Sigma>0$ by (\ref{EqnConformalScalar}), the standard Bochner technique along with non-negativity of $\widetilde{s}_\Sigma$ gives $\widetilde\triangle_\Sigma{r}\le{r}^{-1}$ so $\widetilde\triangle_\Sigma{r}^{-1}\ge{r}^{-3}$.
Thus
\be
\begin{aligned}
&\mathcal{V}^{\frac12}\,-\,\epsilon\left(r^{-1}\,-\,R^{-1}\right) \;\ge\;0 \quad\quad \text{on $\partial\mathcal{A}$} \\
&\widetilde{\triangle}_\Sigma\left(\mathcal{V}^{\frac12}\,-\,\epsilon\left(r^{-1}\,-\,R^{-1}\right)\right)\;\le\;-\epsilon{r}^{-3}\;<\;0.
\end{aligned}
\ee
We have super-harmonicity within $\mathcal{A}$ and non-negativity on $\partial\mathcal{A}$, so $\mathcal{V}^{\frac12}-\epsilon\left(r^{-1}-R^{-1}\right)>0$ in $\mathcal{A}$.
The using the fact that $\lim_{t\rightarrow{R}}\mathcal{V}^\frac12\circ\gamma(t)=0$ and $\lim_{r\rightarrow{R}}\left(r^{-1}-R^{-1}\right)\circ\gamma(t)=0$,  we must have
\be
\begin{aligned}
&\lim_{t\rightarrow{R}}\frac{\partial}{\partial{t}}\left[\left(\mathcal{V}^{\frac12}-\epsilon\left(r^{-1}-R^{-1}\right)\right)\circ\gamma(t)\right]\le0, \quad \text{so} \\
&\lim_{t\rightarrow{R}}\frac{\partial}{\partial{t}}\left(\mathcal{V}^{\frac12}\circ\gamma\right)(t)
\le\;\epsilon\lim_{t\rightarrow{R}}\frac{\partial}{\partial{t}}\left(\left(r^{-1}-R^{-1}\right)\circ\gamma\right)(t)=-\epsilon{R}^{-2}.
\end{aligned}
\ee
Therefore, for some small $\delta$, we must have $\left|\frac{\partial}{\partial{t}}\mathcal{V}^{\frac12}\circ\gamma\right|>\delta$ for $t\in(R-\delta,R)$.
Thus also $\left(\mathcal{V}^{\frac12}\circ\gamma\right)(t)\ge\delta(R-t)$.
Then to estimate the $g_\Sigma$-length of $\gamma$ as $t$ goes to its limiting value $R$, we use
\be
\begin{aligned}
|\dot\gamma|_{g_\Sigma}
&\;=\;\left(\mathcal{V}^{-\frac14}\circ\gamma\right)\cdot|\dot\gamma|_{\tilde{g}_\Sigma}
\;\le\;\delta^{-\frac12}(R-t)^{-\frac12}
\end{aligned}
\ee
for $t\in(R-\delta,R)$ to get
\be
\begin{aligned}
Length_{g_\Sigma}(\gamma)
&\;=\;\int_{R-\delta}^R|\dot\gamma|_{g_\Sigma}dt
\;\le\;\int_{R-\delta}^R\delta^{-\frac12}\cdot(R-t)^{-\frac12}dt
\;=\;2 \;<\;\infty.
\end{aligned}
\ee
This is impossible, and completes the proof.
\qed

\begin{theorem} [cf. Theorem \ref{TheoremFlatRTwo}] \label{ThmRTwoClassification}
Assume $(\Sigma^2,g_\Sigma)$ is metrically complete and $\triangle_\Sigma\sqrt{\mathcal{V}}\le0$.
Then $(\Sigma^2,\widetilde{g}_\Sigma)$ and $(\Sigma^2,g_\Sigma)$ are flat Riemannian manifolds; indeed $g_\Sigma$ has constant coefficients when expressed in $\varphi^1$, $\varphi^2$ coordinates.
\end{theorem}
{\it Proof}.
By Lemma \ref{LemmaConformalComplete}, $(\Sigma^2,\widetilde{g}_\Sigma)$ is also complete, and by (\ref{EqnConformalScalar}) also $\tilde{s}_\Sigma\ge0$.
Thus $(\Sigma^2,\widetilde{g}_\Sigma)$ is parabolic, therefore biholomorphically equivalent to $\mathbb{C}$.
Then $\mathcal{V}^{\frac12}$ is a positive superharmonic function on $\mathbb{C}$, so it is constant.
Lemma \ref{LemmaTConstFlat} gives the conclusion.
\qed

\begin{corollary} [cf. Corollary \ref{CorFlatRTwo}] \label{CorRTwoClassification}
Assume $(M^4,J,\omega,\mathcal{X}_1,\mathcal{X}_2)$ has $s\ge0$, and assume $\mathcal{X}_1$, $\mathcal{X}_2$ have no zeros and are nowhere parallel (in other words $\mathcal{V}$ is nowhere zero).
Then $M^4$ is a flat Riemannian manifold.
\end{corollary}
{\it Proof}.
The corresponding metric polytope $(\Sigma^2,g_\Sigma)$ is complete, and since $\mathcal{V}=0$ if and only if $\{\mathcal{X}_1,\mathcal{X}_2\}$ is a linearly dependent set, the polytope contains no edges.
By equation (\ref{EqnConformalScalar}) we have $\triangle_\Sigma\sqrt{\mathcal{V}}\le0$.
The theorem now implies $g_\Sigma$ is flat; in fact $g_\Sigma$ is a constant matrix.
In the context of the equations (\ref{EqnsGJOmegaM}), we have that $G$ is a constant matrix, so $g$, $J$ are constant.
Thus $M^4$ is flat.
\qed

{\bf Remark.}
In higher dimensions, the methods here fail at almost every stage.

{\bf Remark.}
A crucial part of the proofs of Lemma \ref{LemmaTConstFlat} and Theorem \ref{ThmRTwoClassification} is the fact that a complete $\Sigma^2$ with $K_\Sigma\ge0$ is parabolic.
This is a simple consequence of volume comparison and the Cheng-Yau criterion for parabolicity: that $\int_1^\infty\frac{t}{Vol\,B_t}dt=\infty$.
The assertion that parabolicity of a simply connected Riemann surface implies the surface is actually $\mathbb{C}$ is a simple consequence of uniformization, or even just the Riemann mapping theorem.
The subject of parabolicity has received a great deal of attention: see for example \cite{LT1} \cite{LT2} \cite{HK} \cite{GM} \cite{Web} and references therein.

\section{Global analysis of the coordinate functions} \label{SectionAnalysisHalfPlane}

The purpose of this section is proving that any $\varphi\ge0$ that solves $x(\varphi_{xx}+\varphi_{yy})-\varphi_x=0$ on $\{x>0\}$ with zero boundary conditions, then $\varphi(x,y)=Ax^2$ for some constant $A\ge0$.

\subsection{The elliptic system, coordinate transitions, separation of variables, and first properties}

Equation (\ref{PropDivOfVarphi}) is $dJ_\Sigma{d}\varphi+\mathcal{V}^{\frac12}d\mathcal{V}^{-\frac12}\wedge{J}_\Sigma{d}\varphi=0$, which, with $\triangle_\Sigma\mathcal{V}^{\frac12}=-\frac12s\mathcal{V}^{\frac12}$ in (\ref{EllipticEqnForV}), is the elliptic system
\be
\begin{aligned}
&\triangle_{\Sigma}\varphi\;-\;\left<\nabla_{\Sigma}\log\mathcal{V}^{\frac12},\,\nabla_\Sigma\varphi\right>_\Sigma\;=\;0, \\
&\triangle_{\Sigma}\mathcal{V}^{\frac12}\;+\;\frac12s\mathcal{V}^{\frac12}\;=\;0.
\end{aligned} \label{EqnsPDESystem}
\ee
If $s=0$ then $x=\mathcal{V}^{\frac12}$ is an harmonic coordinate.
Therefore with $H^2=\{x>0\}$ being the open half-plane in $(x,y)$ coordinates, (\ref{EqnsPDESystem}) reduces to
\be
\begin{aligned}
&\varphi_{xx}\,+\,\varphi_{yy}\,-\,x^{-1}\varphi_x\;=\;0, \quad or \quad x(x^{-1}\varphi_x)_x\,+\,\varphi_{yy}\;=\;0
\end{aligned} \label{EqnOrigUnmodifiedPDE}
\ee
on $H^2$ (see equation (3) of \cite{Do1} or equation (2) of \cite{AS}).
By Proposition \ref{PropXISOneOneOnto}, the map $X=(x,y)^T$ is a bijection from the polytope $\Sigma^2$ to the half-plane $\{x\ge0\}$, so therefore this equation is valid not just locally but globally on $\overline{H^2}$ (for an example of what happens if $X$ is not one-to-one, see Example 1).

We are considering the situation where one moment function, say $\varphi^1$, is zero at $x=0$; geometrically this is the situation that the $\Sigma^2=\{\varphi^1\ge0\}$ is a half-plane.
For convenience we will simply denote this variable by $\varphi$.
The fundamental result is that $\varphi(x,y)=Ax^2$ for some constant $A$.
Now the function $x^2$ is unbounded, so $\varphi$ is likely difficult to work with.
Setting $f=\varphi\cdot{x}^{-2}$ and noting that $f$ solves
\be
\begin{aligned}
&f_{xx}\,+\,f_{yy}\,+\,3x^{-1}f_x\;=\;0, \quad or \quad x^{-3}(x^{3}f_x)_x\,+\,f_{yy}\;=\;0, \label{EqnOrigModifiedPDE}
\end{aligned}
\ee
we expect $f$ to be bounded, and so, presumably, easier to work with.
Using the coordinate transformation $s=\frac12x^{-2}$, $t=\frac{1}{\sqrt{8}}y$, equation (\ref{EqnOrigModifiedPDE}) is
\be
\begin{aligned}
s^3f_{ss}\,+\,f_{tt}\;=\;0. \label{EqnNewModifiedPDE}
\end{aligned}
\ee
We shall use both $(x,y)$ and $(s,t)$ coordinates in our proofs below.
Note that the $\{x=\infty\}$ locus becomes the $\{s=0\}$ locus, and vice-versa, so we begin with no information whatever on the behavior of $f$ at $s=0$.
We prove, perhaps unexpectedly, that $f$ is bounded at $s=0$, that it must actually be constant there, and that this constant is the global minimum of $f$; see Corollary \ref{CorBoundedAtSEqZero}.
This is a crucial step toward our ultimate conclusion that $f$ is constant; see Theorem \ref{ThmFinalConstancy}.

\begin{lemma}[$f>0$ even on the boundary] \label{LemmaNewEqnBondaryFacts}
Assume $\varphi\ge0$, $\varphi\in{C}^2(\overline{H^2})\cap{C}^\infty(H^2)$, and $\varphi(0,y)=0$.
Then $f$ is $C^0(\overline{H^2})\cap{}C^\infty(H^2)$, and $f>0$ on $\overline{H^2}$.
\end{lemma}
{\it Proof}.
We can express $\varphi$ as a partial Taylor series with remainder: $\varphi(x,y)=C_0(y)+C_1(y)x+C_2(y)x^2+R$ where $R=o(x^2)$.
In other words, $R\in{C}^2(\overline{H^2})\cap{C}^\infty(H^2)$ and $x^{-2}R\in{C}^0(\overline{H^2})\cap{C}^\infty(H^2)$.

The boundary conditions force $C_0(y)=0$.
Plugging into the PDE gives, near $x=0$, $\lim_{x\rightarrow0}\left(-x^{-1}C_1(y)+R_{xx}-x^{-1}R_x\right)=0$.
Because $\lim_{x\rightarrow0}x^{-1}R_x=\lim_{x\rightarrow0}R_{xx}=0$, then necessarily $C_1(y)=0$.
Thus $f=x^{-2}\cdot\varphi=C_2(y)+x^{-2}{R}\in{C}^0(\overline{H^2})\cap{C}^\infty(H^2)$.

For the assertion that $f>0$ on $\overline{H^2}$, note the maximum principle implies $f>0$ on the interior of $H^2$, so we have to prove that $f$ is non-zero on $\{x=0\}=\partial{H}^2$.
We do this using a lower barrier function.
The barrier will be the function $\psi(x,y)=\epsilon(x^2-4(y-\bar{y})^2)+\epsilon'-\delta{x}^{-2}$ solves $x^{-3}(x^3\psi_x)_x+\psi_{yy}=0$.
On the strip $\{0<x<1\}$, the region on which $\psi$ is positive is pre-compact.
Thus to show that indeed $\psi<f$ on $\{0<x<1\}$ we must only show $\psi<f$ on $\{x=1\}$.
But since $f>0$ on the interior of $H^2$ we may choose $\epsilon>0$, $\epsilon'>0$ so that $\epsilon(1-4(y-\bar{y})^2)+\epsilon'<f$ along $\{x=1\}$ regardless of $\delta>0$.
By sending $\delta\rightarrow0$ we see that $f(0,\bar{y})>\epsilon'$.
Therefore $f$ is not zero anywhere on the boundary $\{x=0\}$.
\qed

{\bf Remark}. The barrier used in this lemma will be used again in Corollary \ref{CorAlmostIncreasingFunction} below.
See Figure (\ref{FigXBarrier}) for its depiction.

{\bf Remark}. The condition that $\varphi=0$ on the boundary is actually sufficient, and we do not need the additional assumption that $\varphi\in{C}^2(\overline{H^2})$.
However we shall not explore this at present, as our geometric situation implies actually $\varphi\in{C}^\infty(\overline{H^2})$.
In Proposition \ref{PropCTwoAndLipschitz}, we proved simply $\varphi\in{C}^2(\overline{H^2})$.
See Example 4 for a solution that is only $C^{1,\alpha}(\overline{H^2})\cap{C}^\infty(H^2)$, although it does not have zero boundary values.

\begin{proposition}[Unspecifiability of boundary values at $x=0$] \label{PropUnspecifiabilityOfBoundaryValues}
Assume $f$ is a bounded solution to (\ref{EqnOrigModifiedPDE}) on any precompact domain $\Omega\subset\overline{H^2}$ with zero boundary values on $\partial\Omega\setminus\{x=0\}$.
Then $f$ is zero.
\end{proposition}
{\it Proof}.
Suppose $f$ is a solution with $f=0$ on $\partial\Omega\setminus\{x=0\}$.
The PDE (\ref{EqnOrigModifiedPDE}) is uniformly elliptic on any open half-plane $0<\delta<x<\infty$ so the strong maximum principle is available on any $\Omega_\delta=\Omega\cap\{x\ge\delta\}$.
The the uniform boundedness of $f$ implies that given any $\epsilon$ there is a $\delta'$ so that for any $\delta\in(0,\delta']$, then $\epsilon{x}^{-2}$ is an upper barrier for $f$ on $\Omega_{\delta}$.
Thus for any $\epsilon>0$, $\epsilon{x}^{-2}$ is an upper barrier for $f$ on $\Omega$.
Sending $\epsilon$ to zero gives $f\le0$.
Replacing $f$ by $-f$ gives also $f\ge0$.
\qed

{\bf Remark.} This proposition is used in Proposition \ref{PropRegularityAtBoundary} to help show that boundary values of $f$ at $\partial\Omega\setminus\{x=0\}$ uniquely determine boundary values at $\{x=0\}$, under the condition of boundedness of $f$.
In particular, the non-uniformly elliptic equation (\ref{EqnOrigModifiedPDE}), in the form $x(f_{xx}+f_{yy})+3f_x=0$, is not even hypoelliptic at the boundary $\{x=0\}$, in that it fails H\"ormander's criterion there \cite{Hor}.
Non-hypoellipticity at the boundary also holds for (\ref{EqnOrigUnmodifiedPDE}); this also follows from H\"ormander's criterion, but one can see this directly from Example 4.

{\bf Remark.} The hypothesis that $\Omega$ is precompact in Proposition \ref{PropUnspecifiabilityOfBoundaryValues} can be relaxed, but we shall not explore this further.
The assumption that $f\in{L}^\infty(\Omega)$ is however necessary; see Example 3 for an example where this is violated.

{\bf Remark.}
The linear equations (\ref{EqnOrigUnmodifiedPDE}) and (\ref{EqnOrigModifiedPDE}) yield easily to separation of variables.
In the following tables the $\lambda$ are positive, and the separation is $f(x,y)=u(x)v(y)$ where either/both $u$, $v$ may be replaced by $\tilde{u}$, $\tilde{v}$.
\be
\begin{array}{|c|l|l|}
\hline
\multicolumn{3}{|c|}{Equation \;\; f_{xx}+f_{yy}-x^{-1}f_x=0} \\
\hline
Eigenvalue & \multicolumn{2}{|c|}{Eigenfunctions} \\
\hline
\multirow{2}{*}{$\lambda^2$}
&u_\lambda(x)\;=\;xJ_1\left({\lambda}\,{x}\right) & v_\lambda(y)\;=\;\sinh\left({\lambda}\,y\right) \\
&\tilde{u}_\lambda(x)=xY_1\left({\lambda}\,x\right) & \tilde{v}_\lambda(y)\;=\;\cosh\left({\lambda}\,y\right) \\
\hline
\multirow{2}{*}{$0$}
&u_0(x)\;=\;x^2 & v_0(y)\;=\;y \\
&\tilde{u}_0(x)\;=\;1 & \tilde{v}_0(y)\;=\;1 \\
\hline
\multirow{2}{*}{$-\lambda^2$}
&u_\lambda(x)\;=\;xI_1\left({\lambda}\,{x}\right) & v_\lambda(y)\;=\;\sin\left({\lambda}\,y\right) \\
&\tilde{u}_\lambda(y)\;=\;xK_1\left({\lambda}\,{x}\right) & \tilde{v}_\lambda(y)\;=\;\cos\left({\lambda}\,y\right) \\
\hline
\end{array} \label{TableUnmodEigenfunctions}
\ee
\be
\begin{array}{|c|l|l|}
\hline
\multicolumn{3}{|c|}{Equation \;\; f_{xx}+f_{yy}+3x^{-1}f_x=0} \\
\hline
Eigenvalue & \multicolumn{2}{|c|}{Eigenfunctions} \\
\hline
\multirow{2}{*}{$\lambda^2$}
&u_\lambda(x)\;=\;x^{-1}J_1\left({\lambda}\,{x}\right) & v_\lambda(y)\;=\;\sinh\left({\lambda}\,y\right) \\
&\tilde{u}_\lambda(x)=x^{-1}Y_1\left({\lambda}\,x\right) & \tilde{v}_\lambda(y)\;=\;\cosh\left({\lambda}\,y\right) \\
\hline
\multirow{2}{*}{$0$}
&u_0(x)\;=\;1 & v_0(y)\;=\;y \\
&\tilde{u}_0(x)\;=\;x^{-2} & \tilde{v}_0(y)\;=\;1 \\
\hline
\multirow{2}{*}{$-\lambda^2$}
&u_\lambda(x)\;=\;x^{-1}I_1\left({\lambda}\,{x}\right) & v_\lambda(y)\;=\;\sin\left({\lambda}\,y\right) \\
&\tilde{u}_\lambda(y)\;=\;x^{-1}K_1\left({\lambda}\,{x}\right) & \tilde{v}_\lambda(y)\;=\;\cos\left({\lambda}\,y\right) \\
\hline
\end{array} \label{TableModEigenfunctions}
\ee
The functions $J_1$, $Y_1$, $I_1$, $K_1$ are the usual Bessel functions.
These eigenfunctions will be used extensively in \S\ref{SubSectionFundSolMethods}.

\subsection{Barrier Function Methods} \label{SubSecBarrierFunctionMethods}

We shall switch freely between the $(x,y)$ and $(s,t)$ coordinate systems and the respective equations (\ref{EqnOrigModifiedPDE}) and (\ref{EqnNewModifiedPDE}).

We take a moment to outline the argument we are about to make.
Proposition \ref{PropIntGradientBounds} gives universal gradient control on any non-negative solution, bounded at $\{x=0\}$ or not, of (\ref{EqnOrigModifiedPDE}); these bounds instantly give polynomial growth/decay in $x$ and exponential growth/decay in $y$.
Then Corollary \ref{CorPolyGrowthInT} improves the exponential control in $y$ to polynomial control---this critical result allows the use of Fourier transform methods in \S\ref{SubSectionFundSolMethods}.
Corollary \ref{CorAlmostIncreasingFunction} proves an almost-monotonicity theorem, meaning $f$ is ``almost'' decreasing in $x$, or increasing in $s$.
Finally Corollary \ref{CorBoundedAtSEqZero} proves that st $s=0$ (which is $x=\infty$) $f$ is actually continuous and constant with respect to $t$, and this constant is actually the global infimum of $f$; this is necessary in the proof of Lemma \ref{LemmaFourierOnBigX}.

\begin{proposition}[Interior gradient bounds] \label{PropIntGradientBounds}
There is some universal constant $c$ so that if $f\in{C}^\infty(H^2)$ is a non-negative solution to $x^{-3}(x^3f_x)_x+f_{yy}=0$, then $x|\nabla\log{f}|<c$.
\end{proposition}
{\it Proof}. If not, there exists some sequence of points $(x_i,y_i)$ so that $x_i\left|\nabla\log{f}\right|_{(x_i,y_i)}\rightarrow\infty$.

We make a point-picking improvement argument to find a ball on which $|\nabla\log{f}|$ is controlled by $|\nabla\log{f}|_{(x_i,y_i)}$ on a sufficiently large ball around $(x_i,y_i)$.
Note that we can assume the ball of radius $2i|\nabla\log{f}|^{-1}$ actually remains in the half-plane, meaning $2i|\nabla\log{f}|^{-1}<<x_i$.
This is because $x_i|\nabla\log{f}|_{(x_i,y_i)}\rightarrow\infty$, so we can pass to a subsequence if necessary to ensure $i<<x_i|\nabla\log{f}|$.
Now if there is any point $(x'_i,y'_i)$ within the ball of radius $i\cdot\left(|\nabla\log{f}|_{(x_i,y_i)}\right)^{-1}$ with $|\nabla\log{f}|_{(x'_i,y'_i)}>2|\nabla\log{f}|_{(x_i,y_i)}$, then re-choose the point $(x_i,y_i)$ to be the point with this larger gradient.
There may be a nearby point with a larger gradient still, so we can repeat this process.
We may repeat this as many times as necessary and still remain all within a ball of radius $2i|\nabla\log{f}|^{-1}$ from the original point, and therefore still remain in the half-plane.
Thus the pointpicking process terminates, resulting in a point $(x_i,y_i)$ with the following two properties:  first the ball of radius $i\cdot\left(|\nabla\log{f}|_{(x_i,y_i)}\right)^{-1}$ is in the right half-plane, and second $|\nabla\log{f}|\le2|\nabla\log{f}|_{(x_i,y_i)}$ at all points in a ball of radius $i\cdot\left(|\nabla\log{f}|_{(x_i,y_i)}\right)^{-1}$.

Now for each $i$ we transition from the coordinates $(x,y)$ to coordinates $(\xi,\upsilon)$ by setting $x=\alpha_i\xi+x_i$, $y=\alpha_i\upsilon+y_i$ where $\alpha_i=\left(|\nabla{f}|_{(x_i,y_i)}\right)^{-1}$.
Using $\nabla_{(x,y)}$ and $\nabla_{(\xi,\upsilon)}$ to distinguish the gradients in the $(x,y)$ and the $(\xi,\upsilon)$ coordinates systems, we have the following properties:
\begin{itemize}
\item[{\it{i}})] $|\nabla_{(\xi,\upsilon)}\log{f}|_{(0,0)}=1$ and $|\nabla_{(\xi,\upsilon)}\log{f}|<2$ on the $(\xi,\upsilon)$-coordinate ball of radius $i$,
\item[{\it{ii}})] $x_i\alpha_i^{-1}\;=\;x_i\left|\nabla_{(x,y)}f\right|_{(x_i,y_i)}\;\longrightarrow\;\infty$,
\item[{\it{iii}})] $\frac{d^2f}{d\xi^2}\,+\,\frac{3}{\xi+x_i\alpha_i^{-1}}\frac{d}{d\xi}\,+\,\frac{d^2f}{d\upsilon^2}\;=\;0$ on the ball of radius $i$.
\end{itemize}
We can scale so that $f(0,0)=1$, and since $f$ and $|\nabla{f}|$ have at worst exponential growth by ({\it{i}}), we have at least $C^{0,\alpha}$ convergence.
By ({\it{ii}}) and ({\it{iii}}) $f$ satisfies an elliptic differential with bounded coefficients, so the convergence is actually $C^\infty$.
Then passing to the limit, we obtain a function $f\ge0$ defined on the entire $(u,v)$-plane so that $|\nabla{f}|=1$ at one point, but since $x_i\alpha_i\rightarrow\infty$ also $f$ satisfies the Laplace equation $f_{\xi\xi}+f_{\upsilon\upsilon}=0$.
However a non-negative harmonic function on $\mathbb{R}^2$ is constant by Liouville's theorem, so we have a contradiction.
\qed

With $x\frac{d}{dx}=-2s\frac{d}{ds}$ and $x\frac{d}{dy}=\frac{1}{4}s^{-\frac12}\frac{d}{dt}$, the previous lemma gives polynomial control over $s\mapsto{f}(s,t)$ and exponential control over $t\mapsto{f}(s,t)$.
\begin{corollary} \label{CorPolyAndExpBound}
There is a constant $c<\infty$ such that $s\left|\frac{d\log\,f}{ds}\right|<c$ and $s^{-1/2}\left|\frac{d\log{f}}{dt}\right|<c$.
In particular, $\left(\frac{s}{s'}\right)^{-c}\le\frac{f(s,t)}{f(s',t)}\le{}\left(\frac{s}{s'}\right)^c$ for $0<s'<s$ and $e^{-c\sqrt{s}(t-t')}\le\frac{f(s,t)}{f(s,t')}<e^{c\sqrt{s}(t-t')}$ for $0<t'<t$.
\end{corollary}
Next we improve the exponential constraint in $t$ to a polynomial constraint.
This is done using a choice of lower barrier; the barrier is related to the one from Lemma \ref{LemmaNewEqnBondaryFacts}.
\begin{proposition}[Quadratic Growth/Decay in $t$] \label{PropPolynomialTBounds}
Assume $f\in{C}^\infty(H^2)$ is a non-negative solution to (\ref{EqnNewModifiedPDE}).
Given any point $(s_0,t_0)$ in the right half-plane, and given any $M\in\mathbb{R}$, we have that
\be
\begin{aligned}
f(s_0,t)&\ge\;
f(s_0,t_0)\left[\frac{c^2-1}{2}\,\frac{\left(\frac{c}{c-1}\right)^c\left(\tau_0^2-\sqrt{\tau_0^4-\frac{c^2-1}{c^2}}\right)^c}{-\tau_0^2\;+\;c\sqrt{\tau_0^4-\frac{c^2-1}{c^2}}}\right]\cdot{}s_0(t-M)^2,  \\
& where \quad \tau_0=\sqrt{1+\frac12s_0(t_0-M)^2} \\
\end{aligned}
\ee
for $t\in[M,t_0]$ or $t\in[t_0,M]$ respectively as $M<t_0$ or $M>t_0$.
If some sharpness is sacrificed, this can be simplified to the estimate
\be
\begin{aligned}
f(s_0,t)
&\;\ge\;
\frac{f(s_0,t_0)}{\left(2+s_0(t_0-M)^2\right)^{c}}\frac{(t-M)^2}{(t_0-M)^2},
\end{aligned}
\ee
for $t\in[M,t_0]$ or $t\in[t_0,M]$ respectively as $M<t_0$ or $M>t_0$.
\end{proposition}
{\bf Remark}. Due to translation invariance, it may typically be convenient to assume $t_0=0$.

\noindent\begin{figure}[h!]
\centering
\caption{\it Graphic for the barrier in Proposition \ref{PropPolynomialTBounds}, drawn using $(s_0,t_0)=(1,0)$ and $f_0=f(s_0,t_0)=2$.
 }
\noindent\begin{subfigure}[b]{0.4\textwidth}
\includegraphics[scale=0.65]{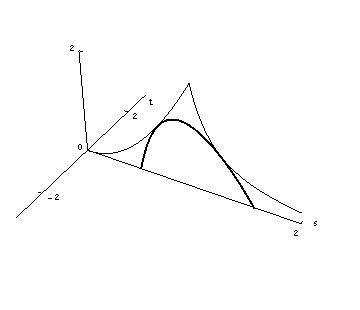}
\label{FigTBarrier}
\caption{Thin curve is the polynomial lower bound on $s\mapsto{}f(s,t_0)$ from Corollary \ref{CorPolyAndExpBound}.
Thick curve is the ``optimal'' curve described in equation (\ref{IneqInequalityBounds}).}
\end{subfigure}
\quad\quad\quad
\quad\quad\quad
\begin{subfigure}[b]{0.4\textwidth}
\includegraphics[scale=0.65]{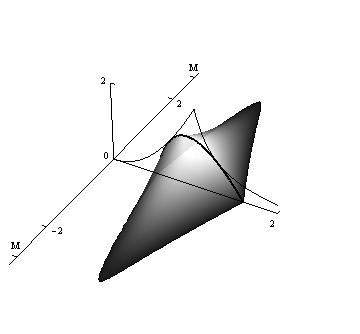}
\caption{Computed lower barrier $\psi$ shown for both $M<t_0$ and $M>t_0$.
The depiction here is for $t_0=0$; see equation (\ref{IneqBoundingFFromBelowByBarrier})}
\end{subfigure}
\end{figure}
{\it Proof}.
Given $(s_0,t_0)$, the previous lemma implies that the function
\be
H_{t_0}(s)\;=\;\begin{cases}
f(s_0,t_0)\left(\frac{s}{s_0}\right)^c \quad & s\in[0,s_0) \\
f(s_0,t_0)\left(\frac{s}{s_0}\right)^{-c} & s\in[s_0,\infty)
\end{cases}
\ee
is a lower barrier: $f(s,t_0)>H_{t_0}(s)$.
The function
\be
\psi(s,t)\;=\;Ds_0^{-1}\left(-s_0s^{-1}-s_0^{-1}s+2+s_0(t-M)^2\right)
\ee
satisfies $s^3\psi_{ss}+\psi_{tt}=0$, and we will find $D=D(M,s_0,t_0)$ so that on $\{t=t_0\}$ we have $\varphi(s,t_0)<H_{t_0}(s)$.
Then restricting to the precompact domain
\be
\Omega\;\triangleq\;\left\{\varphi\,>\,0\right\}\,\cap\,\{t\,<\,t_0\},
\ee
this implies that $\psi<f$ on $\Omega$.
To simplify appearances a bit, we'll use substitutions $\sigma=\frac{s}{s_0}$, $\tau=\sqrt{1+s_0(t-M)^2/2}$, and $\tau_0=\sqrt{1+s_0(t_0-M)^2/2}$, and then choose the largest $D>0$ so that
\be
\begin{aligned}
&f_0\sigma^c\,-\,Ds_0^{-1}\left(-\sigma^{-1}-\sigma+2\tau_0^2\right)\;\ge\;0, \quad 0<\sigma<1 \\
&f_0\sigma^{-c}\,-\,Ds_0^{-1}\left(-\sigma^{-1}-\sigma+2\tau_0^2\right)\;\ge\;0, \quad 1<\sigma<\infty \\
\end{aligned} \label{IneqInequalityBounds}
\ee
where $f_0\triangleq{f}(s_0,t_0)$ and $\tau_0^2=1+s_0(t_0-M)^2/2$.
Noticing that these two equations are equivalent after the transformation $\sigma\mapsto\sigma^{-1}$, we need only consider one of them, say the first one.
Given any $D>0$ the function $G(\sigma)=f_0\sigma^c\,-\,Ds_0^{-1}\left(-\sigma^{-1}-\sigma+2\tau_0^2\right)$ has a single extremum, which is a global minimum.
The optimal choice for $D$ would then solve simultaneously $G'(\sigma)=0$, $G(\sigma)=0$.
We find that
\be
\begin{aligned}
&\sigma\;=\;\frac{c}{c-1}\left(\tau_0^2\;-\;\sqrt{\tau_0^4-\frac{c^2-1}{c^2}}\right), \\
&Ds_0^{-1}\;=\;\frac{(c+1)f_0}{2}\;\frac{\sigma^{c}}{\tau_0^2-\sigma}
\;=\;\frac12(c^2-1)f_0\,\frac{\left(\frac{c}{c-1}\right)^c\left(\tau_0^2-\sqrt{\tau_0^4-\frac{c^2-1}{c^2}}\right)^c}{-\tau_0^2\;+\;c\sqrt{\tau_0^4-\frac{c^2-1}{c^2}}}.
\end{aligned} \label{EqnEstForDs0Inv}
\ee
Therefore
\be
\begin{aligned}
\psi(s,t)
&\;=\;
\frac12(c^2-1)f_0\,\frac{\left(\frac{c}{c-1}\right)^c\left(\tau_0^2-\sqrt{\tau_0^4-\frac{c^2-1}{c^2}}\right)^c}{-\tau_0^2\;+\;c\sqrt{\tau_0^4-\frac{c^2-1}{c^2}}}\left(-\left(\frac{s}{s_0}\right)^{-1}\,-\,\left(\frac{s}{s_0}\right)\,+\,2\tau^2\right)
\end{aligned}
\ee
is a lower barrier: $\psi(s,t)<f(s,t)$ on $\Omega$.
Restricting to $t=t_0$ we get our polynomial lower decay bound on $f$:
\be
\begin{aligned}
f(s_0,t)&\ge\;\psi(s_0,t)
\;=\;
\frac12(c^2-1)f_0\,\frac{\left(\frac{c}{c-1}\right)^c\left(\tau_0^2-\sqrt{\tau_0^4-\frac{c^2-1}{c^2}}\right)^c}{-\tau_0^2\;+\;c\sqrt{\tau_0^4-\frac{c^2-1}{c^2}}}\,\cdot\,s_0(t-M)^2, \\
& where \quad \tau_0^2=1+\frac12s_0(t_0-M)^2.
\end{aligned} \label{IneqBoundingFFromBelowByBarrier}
\ee
It is possible to estimate the expression for $Ds_0^{-1}$ in (\ref{EqnEstForDs0Inv}) by something simpler; for instance
\be
Ds_0^{-1}\;\ge\;\frac{f_0}{2(2\tau_0^2)^c\left(\tau_0^2-1\right)}
\ee
so we have the simplified estimates
\be
\begin{aligned}
&f(s,t)\;\ge\;\frac{f_0}{2(2\tau_0^2)^c\left(\tau_0^2-1\right)}\left(-\left(\frac{s_0}{s}\right)\,-\,\frac{s}{s_0}\,+\,2\,+\,s_0(t-M)^2\right) \\
&f(s_0,t)\;\ge\;\frac{f_0}{(2\,+\,s_0(t_0-M)^2)^c}\frac{(t-M)^2}{(t_0-M)^2}
\end{aligned}
\ee
\qed

{\bf Remark}. If one works even harder at choosing an optimal barrier in Proposition \ref{PropPolynomialTBounds}, one may obtain $f(s_0,t)<Cf(s_0,t_0)\left(t/t_0\right)^{-2}$.
However this will not be necessary for us.
See the conjecture at the end of this subsection.

\begin{corollary}[Almost-monotonicity in $s$] \label{CorAlmostIncreasingFunction}
Assume $f\in{C}^\infty(H^2)$ is a non-negative solution of (\ref{EqnOrigModifiedPDE}).
There is a number $\alpha\in(0,1)$ with the following property.
Fixing $(s_0,t_0)$, if $s>s_0$ then we have $f(s,t_0)>\alpha{f}(s_0,t_0)$.
\end{corollary}

\noindent\begin{figure}[h!]
\centering
\caption{\it Graphic for the barrier in Corollary \ref{CorAlmostIncreasingFunction}, drawn using $(s_0,t_0)=(1,0)$ and $\alpha{f}_0=2$.
 }
\noindent\begin{subfigure}[b]{0.4\textwidth}
\includegraphics[scale=0.65]{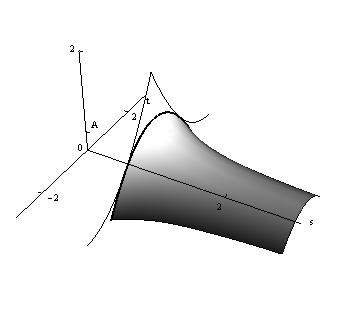}
\caption{The thin curve is the quadratic lower bound on $s\mapsto{}f(s,t_0)$ from Proposition \ref{PropPolynomialTBounds}.
The darker curve underneath is the ``optimal'' parabola determined by (\ref{EqnsOptimalDForSBarrier}).
The barrier $\psi(s,t)$ is shown in the limit where $\eta=0$.}
\label{FigSBarrier}
\end{subfigure}
\quad\quad\quad
\quad\quad\quad
\begin{subfigure}[b]{0.4\textwidth}
\includegraphics[scale=0.6]{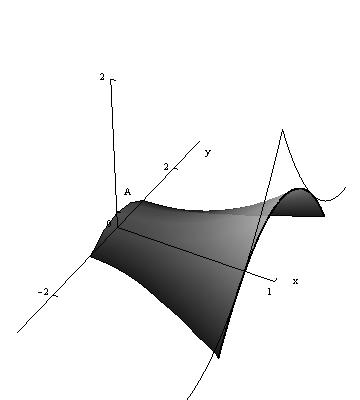}
\caption{This is the same barrier $\psi$, expressed in the $(x,y)$-coordinate system. \\
\\
}
\label{FigXBarrier}
\end{subfigure}
\end{figure}

{\it Proof}.
Again the strategy is to place a test function underneath $f$.
We use
\be
\psi(s,t)\;=\;D\left(s^{-1}\,-\,(t-t_0)^2\right)\,-\,\eta^2(s-s_0)\,+\,A
\ee
for positive $D$.
Similar to the previous lemma, note that the domain $\Omega=\left\{\psi>0\right\}\cap\{s>s_0\}$ is precompact, and $\psi=0$ on $\partial\Omega\setminus\{s=s_0\}$, so to check that $\psi<f$ on $\Omega$, it is enough to check that $\psi\le{f}$ on $\{s=s_0\}$.
Setting $s=s_0$ the previous lemma gives
\be
\begin{aligned}
&f(s_0,t)\;\ge\;\alpha\,\cdot\,f_0\,\cdot\,\frac{(t-M)^2}{(t_0-M)^2}, \quad \text{where} \quad
\alpha\;=\;\frac{\frac12(c+1)s_0(t_0-M)^2\,\sigma_0^c}{\tau_0-\sigma_0},  \\
& \quad\quad
\sigma_0=\frac{c}{c-1}\left(\tau_0^2-\sqrt{\tau_0^4-\frac{c^2-1}{c^2}}\right), \quad
\tau_0^2\;=\;1+\frac12s_0(t_0-M)^2.
\end{aligned}
\ee
We pick some $N$ (to be chosen below), and find the optimal $D$ so that $N-D(t-t_0)^2\le\alpha{f}_0\left(\frac{t-M}{t_0-M}\right)^2\le{f}(s_0,t)$.
To do so, we again use the first derivative trick, and simultaneously solve
\be
\begin{aligned}
H(t)&\;=\;\alpha{f}_0\left(\frac{t-M}{t-t_0}\right)^2\,-\,N\,+\,D(t-t_0)^2\;=\;0 \\
H'(t)&\;=\;2\alpha{f}_0\frac{t-M}{(t-t_0)^2}\,+\,2D(t-t_0) \;=\;0
\end{aligned} \label{EqnsOptimalDForSBarrier}
\ee
for $t$ and $D$.
We find $D=\frac{N\alpha{f}_0}{(\alpha{f}_0-N)(t_0-M)^2}$.
Therefore the optimal barrier is
\be
\begin{aligned}
&\psi(s,t)\;=\;N\left(1\,-\,\frac{\alpha{f}_0}{(\alpha{f}_0-N)s_0(t_0-M)^2}\right) \\
& \quad\quad\quad\quad\,+\,\frac{N\alpha{f}_0}{(\alpha{f}_0-N)s_0(t_0-M)^2}\left(s_0s^{-1}\,-\,s_0(t-t_0)^2\right) \;+\;\eta^2(s-s_0).
\end{aligned}
\ee
Since the $D$ given above is the solution to the system (\ref{EqnsOptimalDForSBarrier}), we have
\be
\begin{aligned}
&\psi(s_0,t)\;=\;N\,-\,\frac{N\alpha{f}_0(t-t_0)^2}{(\alpha{f}_0-N)(t_0-M)^2} \;\le\; \alpha{f}_0\left(\frac{t-M}{t_0-M}\right)^2\;\le\;f(s_0,t).
\end{aligned}
\ee
Now set $N=(1-\mu)\alpha{f}_0$, to obtain
\be
\begin{aligned}
&\psi(s,t)\;=\;\frac{1-\mu}{\mu}\alpha{f}_0\frac{\mu{}s_0(t_0-M)^2-1}{s_0(t_0-M)^2} \\
&\quad\quad\quad\quad\,+\,\frac{1-\mu}{\mu}\frac{\alpha{f}_0}{s_0(t_0-M)^2}\left(s_0s^{-1}\,-\,s_0(t-t_0)^2\right) \;+\;\eta^2(s-s_0).
\end{aligned}
\ee
Choose $M=M(s_0,t_0)$ so that $s_0(t_0-M)^2=4$.
Then $\tau_0^2=3$, $\sigma_0=\frac{c}{c-1}(3-\sqrt{8+c^{-2}})$, and $\alpha=\frac{(c+1)\sigma_0^c}{3-\sigma_0}$.
Finally choose $\mu=\frac12$.
Then
\be
\psi(s,t)=\frac18\alpha{}f_0\,+\,\frac18\alpha{}f_0\left(\left(\frac{s}{s_0}\right)^{-1}-\left(\frac{t-t_0}{t-M}\right)^2\right)\,+\,\eta^2(s-s_0) \label{EqnBarrierS}
\ee
Now we have that $\psi<f$ on $\partial\Omega$, so $\psi<f$ on $\Omega$, and this is independent of $\eta^2$.
Sending $\eta\rightarrow0$ and setting $t=t_0$, we obtain $f(s,t_0)>\frac{\alpha}{8}f(s_0,t_0)$.
\qed

\begin{corollary}[Polynomial Upper bounds for $t$] \label{CorPolyGrowthInT}
There exists some $c<\infty$ so that if $f\in{C}^\infty(H^2)$ is a non-negative solution of (\ref{EqnNewModifiedPDE}), then
\be
\begin{aligned}
&f(s,t_0)\;\le\;2^{c+2}(1+st^2)^cf(s,t)
\end{aligned}
\ee
for any $t,t_0\in\mathbb{R}$.
In the $(x,y)$ system, this is $f(x,y_0)\;\le\;C\left(\frac{x^2+(y-y_0)^2}{x^2}\right)^c\,f(x,y).$
\end{corollary}
{\it Proof}.
First assume $t_0>t$ in the conclusion of Corollary \ref{PropPolynomialTBounds}, and set $M=t-(t_0-t)=2t-t_0$ so $t_0-M=2(t_0-t)$ and $t-M=t-t_0$.
We get
\be
\begin{aligned}
&f(s_0,t)\;\ge\;\frac{f(s_0,t_0)}{(1+s_0(t_0-t)^2)^c}\frac{1}{2^{c+2}}.
\end{aligned}
\ee
\qed

\begin{corollary}[Global minimum occurs at $s=0$] \label{CorBoundedAtSEqZero}
If $f\ge0$ solves (\ref{EqnNewModifiedPDE}) on $\{s\ge0\}$, then $f$ is continuous at all boundary points $\{s=0\}$, $f|_{s=0}$ is constant, and the value of $f$ at $s=0$ is a strict minimum.
\end{corollary}
{\it Proof}.
For $(s_0,t_0)$, again choose $M=M(s_0,t_0)$ so $s_0(t_0-M)^2=4$ as above, but now assume $t\in\left[t_0-\sqrt{\frac{2}{s_0}},t_0+\sqrt{\frac{2}{s_0}}\right]$ so $|t|\le(M-t_0)/\sqrt{2}$.
Then using the barrier (\ref{EqnBarrierS}) with these values of $t$ we have
\be
\begin{aligned}
\psi(s,t)
&\;\ge\;\frac{1}{16}\alpha{}f_0\,+\,\eta^2(s-s_0) \\
\end{aligned}
\ee
Again we can send $\eta\rightarrow0$, so therefore $f(s,t)>\frac{1}{16}\alpha{f}_0$ on $(s,t)\in\left[s_0,\infty\right)\times\left[t_0-\sqrt{\frac{2}{s_0}},t_0+\sqrt{\frac{2}{s_0}}\right]$ (where $f_0=f(s_0,t_0)$).
This is independent of $s_0$.
So let $(s_i,t_i)$ be a sequence of points where $s_i\searrow0$ and $t_i\rightarrow{}t_\infty$.
Define numbers $f_i=f(s_i,t_i)$.
First, it is impossible that $\sup_i{}f_i=\infty$, because $f(s,t_i)>\frac18\alpha{f}(s_i,t_i)$ for $s>s_i$ would then force $f(s,t_i)\rightarrow\infty$ for any fixed $s$, which is impossible.

Define the values of $f$ on $\{s=0\}$ by simply setting $f(0,t)={\lim\,\inf}_{(s',t')\rightarrow(0,t)}f(s',t')$.
Also, by subtracting a constant from $f$ if necessary, assume $\inf{}f(s,t)=0$ (where the infimum is taken over the whole half-plane).

Now passing to a subsequence of $(s_i,t_i)$, we can assume $f_i$ converges to some $f_\infty<\infty$, where the subsequence can be chosen so that $f_\infty={\lim\,\sup}_{i\rightarrow\infty}f(s_i,t_i)$.
Then (using the extended definition of $f$) we have that $f(0,t)\ge\frac{1}{16}\alpha{f}_\infty$.
Then from the previous result also $f(s,t)\ge\frac{1}{8}\alpha\left(\frac{1}{16}\alpha{f}_\infty\right)$.
But $\inf\,{f}=0$, so necessarily $f_\infty=0$.
Since $f\ge0$ and $f_\infty$ was chosen as a ``$\lim\inf$,'' we have that if $(s_j,t_j)$ is {\it any} other sequence that converges to $(0,t)$, this forces $\lim_jf(s_j,t_j)$ to exist, and to equal zero.
But $t$ is arbitrary, so $f$ is continuous and constant at $s=0$, as claimed.
\qed

{\bf Remark}. The foundational result is Proposition \ref{PropIntGradientBounds}, which cannot hold for supersolutions; see Example 1.
Nevertheless we conjecture that Corollary \ref{CorPolyAndExpBound} holds for supersolutions, whether bounded at $\{x=0\}$ or not, and also that $c$ can be taken to be $2$.

\subsection{Fundamental Solution Methods} \label{SubSectionFundSolMethods}

If $f$ is a solution function and $x_0>0$, then the function $y\mapsto{f}(x_0,y)$ has polynomial growth in $y$ by Corollary \ref{CorPolyGrowthInT}.
The main idea of this section is that the polynomial growth bound allows the use of Fourier transform methods, and allows us to construct fundamental solutions on strips $\{0<x<x_0\}$ as well as half-planes $\{x_0<x\}$.
By the standard method of convolutions of boundary data with with fundamental solutions, we are able to then show that necessarily $y\mapsto{f}(0,y)$ also has polynomial growth.
Then we show that this is impossible.

For the moment we require the apriori assumption that $y\mapsto{f}(0,y)$ already has exponential growth: this technical assumption is necessary in the proof of \ref{LemmaExpBoundsImplyFourierOnSmallX}, where we show that, on strips $\{0<x<\epsilon\}$ then a solution $f$ is actually equal to the usual convolution with a fundamental solution.
This exponential assumption is removed in \S\ref{SubSectionGeneralCaseHTwo}.

A fundamental solution for (\ref{EqnOrigUnmodifiedPDE}) can be written in terms of elementary functions, namely
\be
G(x,y)\;=\;\frac{\frac12x^2}{(x^2+y^2)^{3/2}}. \label{EqnFundSolForOrigEqn}
\ee
This solves $G_{xx}-x^{-1}G_x+G_{yy}=0$, and, restricted to $\{x=0\}$, $y\mapsto{}G(0,y)$ is the Dirac delta (this is justified in Example 2 below).
If $\psi(y)$ are boundary values on $\{x=0\}$, then
\be
\varphi(x,y)\;=\;\int_{-\infty}^\infty G(x,\,y-t)\,\psi(t)\,dt
\ee
solves $\varphi_{xx}-x^{-1}\varphi_x+\varphi_{yy}=0$ with $\varphi(0,y)=\psi(y)$.

By contrast, in light of Proposition \ref{PropUnspecifiabilityOfBoundaryValues} there can be no fundamental solution for the equation $f_{xx}+3x^{-1}f_x+f_{yy}=0$ on $\{x\ge0\}$.
We have the following fundamental solutions on certain right half-planes and certain strips:
\be
\begin{aligned}
&\mathcal{G}^\epsilon(x,y)\;=\;\frac{\epsilon}{x}\int_0^\infty\frac{K_1(\omega{}x)}{K_1(\omega\epsilon)}\,\cos(\omega{}y)\;d\omega, \quad \text{on}\;\{x\ge\epsilon\}, \\ 
&\mathcal{G}_\epsilon(x,y)\;=\;\frac{\epsilon}{x}\int_0^\infty\frac{I_1(\omega{x})}{I_1(\omega\epsilon)}\,\cos(\omega{}y)\,d\omega, \quad \text{on}\;\{0\le{x}\le\epsilon\},
\end{aligned} \label{EqnFundSol}
\ee
where $K_1$ and $I_1$ are the usual modified Bessel functions of the first and second kinds.

Define one-variable functions $\psi_\epsilon$ by
\be
\psi_\epsilon(y)=f(\epsilon,y)
\ee
For simplicity, from here on we assume $f(x,y)=f(x,-y)$, as we could always replace $f(x,y)$ with $f(x,y)+f(x,-y)$.
This allows us to use a Fourier cosine representation for $\psi_\epsilon$, which the author likes better than the exponential representation.
\begin{lemma} \label{LemmaBasiFourierProperties}
The Fourier transform $\widehat{\psi_\epsilon}$ of $\psi_\epsilon$ exists in the sense of distributions, and $\widehat{\psi_\epsilon}(\omega)$ is smooth except at $\omega=0$.
\end{lemma}
{\it Proof}.
Since $\psi_\epsilon$ is an even function, it has a Fourier cosine transform
\be
\widehat{\psi_\epsilon}(\omega)\;=\;\int_{-\infty}^\infty \psi_\epsilon(t)\cos(\omega{t})\,dt.
\ee
The function $\psi_\epsilon(y)$ has polynomial growth at worst by Corollary \ref{CorPolyGrowthInT}, so this transform exists in the distributional sense.
Because $\psi_\epsilon$ is non-negative and therefore has no large oscillations, the transform must be smooth except possibly at $\omega=0$.
\qed

\begin{lemma} \label{LemmaFourierOnBigX}
We have that $f(x,y)=\left(\mathcal{G}^\epsilon*_y\psi_\epsilon\right)(x,y)$ on the half-plane $\{x\ge\epsilon\}$.
\end{lemma}
{\it Proof}.
The assertion is that if we set
\be
F^\epsilon(x,y)\;\triangleq\;\int_{-\infty}^\infty\mathcal{G}^\epsilon(x,y-t)\,\psi_\epsilon(t)\,dt, \label{EqnFundRightHalfPlaneDef}
\ee
then actually $f(x,y)=F^\epsilon(x,y)$ on $\{x\ge\epsilon\}$.
To prove this, we first look at the strip $\{\epsilon<x<\epsilon'\}$ and prove that fundamental solution methods let us recover the solution exactly, and then we let $\epsilon'\rightarrow\infty$ and see that we are left with precisely (\ref{EqnFundRightHalfPlaneDef}).

For certain values of $\alpha_\omega$, $\beta_\omega$ to be determined below, define
\be
\mathcal{G}^\epsilon_{\epsilon'}(x,y)
\;=\; \frac1x \int_0^\infty \left(\alpha_\omega{K}_1(\omega{x})\,+\,\beta_\omega{I}_1(\omega{x})\right)\,d\omega
\ee
Then on whatever the domain of convergence happens to be, we have that $\mathcal{G}^\epsilon_{\epsilon'}$ indeed solves (\ref{EqnOrigModifiedPDE}).
Choosing $\alpha_\omega=\frac{\epsilon{}I_1(\omega\epsilon')}{K_1(\omega\epsilon)I_1(\omega\epsilon')-K_1(\omega\epsilon')I_1(\omega\epsilon)}$ and $\beta_\omega=\frac{-\epsilon'{}K_1(\omega\epsilon')}{K_1(\omega\epsilon)I_1(\omega\epsilon')-K_1(\omega\epsilon')I_1(\omega\epsilon)}$ we have that 
\be
\mathcal{G}^\epsilon_{\epsilon'}(\epsilon,y)\;=\;\delta_y, \quad\quad
\mathcal{G}^\epsilon_{\epsilon'}(\epsilon',y)\;=\;0.
\ee
Similarly, if we set
\be
\widetilde{\mathcal{G}}^\epsilon_{\epsilon'}(x,y)
\;=\; \frac1x \int_0^\infty \left(\tilde\alpha_\omega{K}_1(\omega{x})\,+\,\tilde\beta_\omega{I}_1(\omega{x})\right)\,d\omega
\ee
with $\tilde\alpha_\omega=\frac{-\epsilon{}I_1(\omega\epsilon)}{K_1(\omega\epsilon)I_1(\omega\epsilon')-K_1(\omega\epsilon')I_1(\omega\epsilon)}$ and $\tilde\beta_\omega=\frac{\epsilon'{}K_1(\omega\epsilon)}{K_1(\omega\epsilon)I_1(\omega\epsilon')-K_1(\omega\epsilon')I_1(\omega\epsilon)}$
we obtain
\be
\widetilde{\mathcal{G}}^\epsilon_{\epsilon'}(\epsilon,y)\;=\;0, \quad\quad
\widetilde{\mathcal{G}}^\epsilon_{\epsilon'}(\epsilon',y)\;=\;\delta_y.
\ee
Then define
\be
F^\epsilon_{\epsilon'}(x,y)\;=\;
\int_{-\infty}^\infty\mathcal{G}^\epsilon_{\epsilon'}(x,y-t)\psi_\epsilon(t)\,dt
\;+\;\int_{-\infty}^\infty\overline{\mathcal{G}}^\epsilon_{\epsilon'}(x,y-t)\psi_{\epsilon'}(t)\,dt. \label{DefOfFUpper}
\ee
Because of the polynomial growth of $\psi_\epsilon$, $\psi_{\epsilon'}$, the function $F^\epsilon_{\epsilon'}$ is well-defined on the strip $\{\epsilon<x<\epsilon'\}$, and agrees with $f(x,y)$ on the boundary $\{x=\epsilon\}\cup\{x=\epsilon'\}$.
Standard Fourier theory implies that $F^\epsilon_{\epsilon'}(x,y)$ also has polynomial growth in $y$ for each fixed $x\in(\epsilon,\epsilon')$.

Next we show that $F^\epsilon_{\epsilon'}=f$ on the interior of the strip $\{\epsilon<x<\epsilon'\}$.
Note that $F^\epsilon_{\epsilon'}-f$ has zero boundary conditions and at worst polynomial growth; we'll use a barrier argument to show it is exactly zero.
Consider the function
\be
\eta(x,y)\;=\;\frac{1}{x}J_1(x/(2\epsilon'\alpha))\,\cosh(y/(2\epsilon'\alpha)).
\ee
where $\alpha$ is the second zero of $J_1(x)$.
As noted in Table \ref{TableModEigenfunctions}, $\eta$ solves the PDE (\ref{EqnOrigModifiedPDE}).
Further, $\eta(x,y)>0$ for $x\in[0,2\epsilon')$, and indeed we have an exponential lower growth estimate:
\be
\eta(x,y)\;\ge\;\frac{1}{\epsilon'}J_1(1/(2\alpha))\,\cosh(y/(2\epsilon'\alpha)).
\ee
Due to the polynomial growth of $F^\epsilon_{\epsilon'}-f$, we see that give any $\epsilon$, $\epsilon'$ and any positive constant $C$ we have that $C\eta(x,y)>F^\epsilon_{\epsilon'}(x,y)-f(x,y)$ for sufficiently large $y$.
But then the maximum principle says indeed $C\eta(x,y)>F^\epsilon_{\epsilon'}(x,y)-f(x,y)$ on the entire strip.
But $C$ is arbitrary, we can let $C\rightarrow0$ to obtain $F^\epsilon_{\epsilon'}(x,y)-f(x,y)\le0$.
Similarly using $-C\eta(x,y)$ we obtain $F^\epsilon_{\epsilon'}(x,y)-f(x,y)\ge0$, so indeed we have proven that $F^\epsilon_{\epsilon'}(x,y)=f(x,y)$ on the strip.

Next we fix $\epsilon$ and let $\epsilon'\rightarrow0$.
Considering the coefficients $\alpha_\omega$, $\beta_\omega$, $\tilde\alpha_\omega$, $\tilde\beta_\omega$, we have $\alpha_\omega\rightarrow\frac{\epsilon}{K_1(\omega\epsilon)}$, while $\beta_\omega$, $\tilde{\alpha}_{\omega}$, and $\tilde{\beta}_{\omega}$ all converge to zero.
Convergence of the second integral in (\ref{DefOfFUpper}) is not an issue, as the polynomial growth bounds on $\psi_{\epsilon'}(y)$ only improve as $\epsilon'\rightarrow\infty$.
This proves the second equality of
\be
f(x,y)=\lim_{\epsilon'\rightarrow\infty}F^\epsilon_{\epsilon'}(x,y)\;=\;F^\epsilon(x,y)
\ee
and finishes the lemma.
\qed

\begin{lemma} \label{LemmaExpBoundsImplyFourierOnSmallX}
Assume $y\mapsto{}f(0,y)$ has exponential growth bounds, meaning there is some $C$ where $f(0,y)\le{e}^{C|y|}$.
Set $\psi_\epsilon(y)=f(\epsilon,y)$.
If $\epsilon=\epsilon(C)>0$ is sufficiently small, then actually $f(x,y)=(\mathcal{G}_\epsilon*_y\psi_\epsilon)(x,y)$ on $\{0\le{x}\le\epsilon\}$.
Further, $y\mapsto{}f(0,y)$ is actually polynomially bounded: $f(0,y)\le{C}_1+{C}_2|y|^c$ for some universal constants $c$, $C_1$, and $C_2$.
\end{lemma}
{\it Proof}.
We would like to imitate the proof above, setting
\be
F^\epsilon_{\epsilon'}(x,y)\;=\;
\int_{-\infty}^\infty\mathcal{G}^\epsilon_{\epsilon'}(x,y-t)\psi_\epsilon(t)\,dt
\;+\;\int_{-\infty}^\infty\widetilde{\mathcal{G}}^\epsilon_{\epsilon'}(x,y-t)\psi_{\epsilon'}(t)\,dt \label{DefOfFUpper2}
\ee
and this time letting $\epsilon\rightarrow0$.
However the polynomial growth bounds from Corollary \ref{CorPolyGrowthInT} actually deteriorate as $\epsilon\rightarrow0$ (whereas the bounds improve as $\epsilon'\rightarrow\infty$), so even though the coefficients converge: $\alpha_\omega\rightarrow0$, $\beta_\omega\rightarrow0$, $\tilde\alpha_\omega\rightarrow0$, $\tilde\beta_\omega\rightarrow\frac{\epsilon'}{I_1(\omega\epsilon')}$, it is not clear that the first integral in (\ref{DefOfFUpper2}) actually vanishes in the limit.

So we argue differently.
Simply define
\be
\begin{aligned}
F_\epsilon(x,y)
&\;\triangleq\;\int_{-\infty}^\infty \mathcal{G}_\epsilon(x,\,y-t)\,\psi_\epsilon(t)\,dt \\
&\;=\;\frac{\epsilon}{x}\int_{0}^\infty \frac{I_1(\omega x)}{I_1(\omega \epsilon)}\left(\int_{-\infty}^\infty \cos(\omega(y-t))\,\psi_\epsilon(t)\,dt\right)\,d\omega \\
&\;=\;\frac{\epsilon}{x}\int_{0}^\infty \frac{I_1(\omega x)}{I_1(\omega \epsilon)}\cos(\omega{y})\,\widehat{\psi_\epsilon}(\omega)\,d\omega
\end{aligned}
\ee
Then $F_\epsilon(x,y)$ has polynomial growth bounds everywhere, including at $x=0$.
To see this, note that
\be
\begin{aligned}
F_\epsilon(0,y)
&\;=\;\epsilon\int_{0}^\infty \frac{\omega}{2I_1(\omega \epsilon)}\cos(\omega{y})\,\widehat{\psi_\epsilon}(\omega)\,d\omega
\end{aligned}
\ee
which is certainly convergent.
To explain why, notice this is simply the Fourier inverse of $\frac{\epsilon\omega}{2I_1(\omega \epsilon)}\widehat{\psi_\epsilon}(\omega)$, and notice that indeed $\frac{\omega\epsilon}{2I_1(\omega \epsilon)}\widehat{\psi_\epsilon}(\omega)-\widehat{\psi_\epsilon}(\omega)$ is actual a smooth function of $\omega$.
This is seen by noticing that $\lim_{\omega\rightarrow0}\frac{\omega\epsilon}{2I_1(\omega \epsilon)}=1$, and that by Lemma \ref{LemmaBasiFourierProperties} the non-smooth part of $\widehat{\psi_\epsilon}(\omega)$ occurs only at $\omega=0$.
Then, because only the non-smooth part of $\widehat{\psi_\epsilon}$ controls the growth of $F_\epsilon(0,y)$, we see that $F_\epsilon(0,y)$ and $F_\epsilon(\epsilon,y)$ have the same polynomial growth (eg. a Dirac $\delta$-function implies $O(1)$ growth, the $k^{th}$ derivative $\delta^{(k)}$ of a delta implies polynomial growth of order $k$, and so on).

Now consider the function $f(x,y)-F_\epsilon(x,y)$.
This is zero on $\{x=\epsilon\}$, and by assumption has, at worst, exponential growth in $y$ of order $C$ at $\{x=0\}$.
But then if $\epsilon$ is small enough, the function
\be
\eta(x,y)\;=\;A\frac{1}{x}J_1(2Cx)\,\cosh(2Cy)
\ee
exists and is positive on $\{0\le{x}\le\epsilon\}$, has exponential growth in $y$ of order $2C$ at $\{x=0\}$, and solves the PDE (\ref{EqnOrigModifiedPDE}).
Thus, no matter what $A>0$ might be, we have that $\eta(x,y)>f(x,y)-F_\epsilon(x,y)$ for $y$ sufficiently large, and therefore $\eta(x,y)>f(x,y)-F_\epsilon(x,y)$ everywhere on $\{0\le{x}\le\epsilon\}$.
Sending $A\rightarrow0$ gives $f(x,y)\le{F}_\epsilon(x,y)$.
Similarly using $-\eta(x,y)$ we obtain the opposite inequality, and so $f(x,y)=F_\epsilon(x,y)$ on $\{0\le{x}\le\epsilon\}$.
\qed

\begin{proposition} \label{ThmExpBoundsImplyConstancy}
Assume $f\in{C}^0(\overline{H^2})\cap{C}^\infty(H^2)$ is a non-negative solution to (\ref{EqnOrigModifiedPDE}) on the right half-plane, and assume $y\mapsto{}f(0,y)$ has exponential growth bounds.
Then $f(x,y)$ is constant.
\end{proposition}
{\it Proof}.
By subtracting a constant, we can assume $\inf_{x>0}{f}(x,y)=0$.
Lemma \ref{LemmaExpBoundsImplyFourierOnSmallX} implies polynomial bounds on the boundary: $f(0,y)<C_1+C_2y^c$.
By Lemma \ref{LemmaFourierOnBigX} we have that $f(x,y)$ on $\{x\ge\epsilon\}$ is
\be
\begin{aligned}
f(x,y)
&\;=\;\frac{1}{x}\int_{-\infty}^\infty\int_0^\infty \frac{K_1(\omega{x})}{K_1(\omega\epsilon)}\,\cos(\omega(y-t))\psi_\epsilon(t)\,d\omega\, dt \\
&\;=\;\frac{1}{x}\int_0^\infty \frac{K_1(\omega{x})}{K_1(\omega\epsilon)}\widehat{\psi_\epsilon}(\omega)\cos(\omega{y})\,d\omega
\end{aligned}
\ee
Sending $\epsilon\rightarrow0$ forces $f(x,y)\equiv0$.
To see this, note that $x\mapsto\frac{{K}_1(\omega{x})}{x}$ exists and is finite everywhere including $\{x=0\}$, and that $\frac{1}{K_1(\omega\epsilon)}\rightarrow0$ as $\epsilon\rightarrow0$.
Using the uniform polynomial bounds on $\psi_\epsilon(y)$ as $\epsilon\rightarrow0$, we have that $\lim_{\epsilon\rightarrow0}\widehat{\psi}_\epsilon(x)$ converges to a distribution in the distributional sense.
Therefore the function $(x,\omega)\mapsto\lim_{\epsilon\rightarrow0}\frac{K_1(\omega{x})}{{x}K_1(\omega\epsilon)}\widehat{\psi_\epsilon}(\omega)$ converges to 0 in the distributional sense.
\qed

{\bf Remark.} It may be helpful to explain the proof of Proposition \ref{ThmExpBoundsImplyConstancy} in natural language.
We found earlier that, on the half-planes $\{x\ge\epsilon\}$, $f$ equals the convolution of $f(\epsilon,y)$ with the fundamental solutions $\mathcal{G}^\epsilon(x,y)$ produced above.
We also shows that $y\mapsto{f}(0,y)$ indeed must have {\it polynomial} bounds, and therefore we can send $\epsilon\rightarrow0$ without any trouble.
The underlying idea of the proof is that the fundamental solution $\mathcal{G}^\epsilon$ actually converges to $0$ everywhere on $\{x>0\}$.
In the distributional sense, one has actually $\lim_{\epsilon\rightarrow0}\mathcal{G}^\epsilon(x,y)=\delta(x,y)$.

\noindent\begin{figure}[h!]
\centering
\caption{\it The thinning of the fundamental solution on half-planes.}
\noindent\begin{subfigure}[b]{0.4\textwidth}
\includegraphics[scale=0.65]{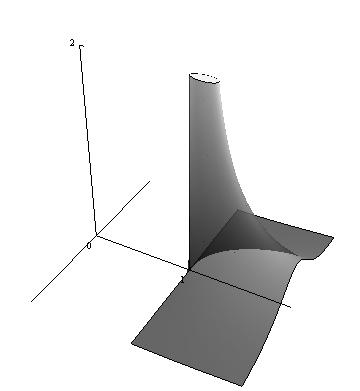}
\caption{Fundamental solution on $\{x\ge1\}$.}
\label{ThinningFig}
\end{subfigure}
\quad\quad\quad
\quad\quad\quad
\begin{subfigure}[b]{0.4\textwidth}
\includegraphics[scale=0.65]{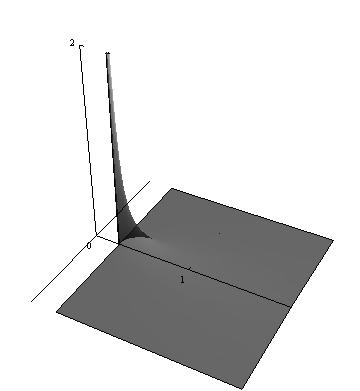}
\caption{Fundamental solution on $\{x\ge\frac14\}$.}
\label{HiFig}
\end{subfigure}
\end{figure}

\subsection{Constancy of solutions in the General Case of Well-Defined Boundary Values} \label{SubSectionGeneralCaseHTwo}

Here we remove the exponential growth condition on $y\mapsto{f}(0,y)$.
To do this, we first show that if $f\in{C}^0(\overline{H^2})\cap{C}^\infty(H^2)$, then $f\in{C}^{\infty}(\overline{H^2})$, and then we use a blow-up argument to obtain a non-constant, non-negative function with exponential growth bounds on the boundary, leading to a contradiction.\footnote{This shows a fundamental difference between equations (\ref{EqnOrigUnmodifiedPDE}) and (\ref{EqnOrigModifiedPDE}); see Example 4.}

\begin{proposition}[Regularity at the boundary] \label{PropRegularityAtBoundary}
Assume $\Omega$ is a pre-compact domain with $\overline{\Omega}\cap\{x=0\}$ non-empty.
Assume $f$ satisfies (\ref{EqnOrigModifiedPDE}) on $\Omega$, and $f\in{C}^0(\overline{\Omega})\cap{C}^\infty(\Omega)$.
Then $f$ is $C^\infty$ on $\Omega\cup\left(\partial\Omega\cap\{x=0\}\right)$.
\end{proposition}
{\it Proof}.
We have $f\in{C}^\infty(\Omega)$.
If $p$ is a point in the interior of $\partial\Omega\cap\{x=0\}$, we can find an open quadrilateral $Q\subset\Omega$ with $p\in\partial{Q}$.
By translation invariance with respect to $y$ and by simultaneous-scale invariance in $x$ and $y$, we may assume $Q=\{0<x<1\}\cap\{-1<y<1\}$.

Now we can use a series representations for the solution $f$ on $Q$.
We get
\be
\begin{aligned}
f(x,y)
&\;=\;A(x, y)\,+\,B(x,y)\,+\,C(x,y)
\end{aligned}
\ee
where $A$, $B$, and $C$ are given by
\be
A(x,y)\;=\;\sum_{n=0}^{\infty} a_n\frac{1}{x}\frac{I_1(\pi{n}x)}{I_1(\pi{n})}\cos(\pi{n}y)\,+\,\sum_{n=0}^{\infty}b_n\frac{1}{x}\frac{I_1(\pi{n}x)}{I_1(\pi{n})}\sin(\pi{n}y)
\ee
where $a_n$, $b_n$ are the usual Fourier cosine and sine coefficients for $y\mapsto{f}(1,y)$;
\be
B(x,y)\;=\;\sum_{n=0}^{\infty}c_n\frac{1}{x\alpha_n}J_1(\lambda_nx)\frac{\sinh(\lambda_n(1+y))}{2\cosh(\lambda_n)\sinh(\lambda_n)}
\ee
where $\lambda_n$ is the $n^{th}$ zero of $J_1$, $\alpha_n=-\frac12J_0(\lambda_n)J_2(\lambda_n)$ is a standard normalizing constant, and the $c_n$ are the Bessel series coefficients for the function $x\mapsto{f}(x,1)-A(x,1)$; and
\be
C(x,y)\;=\;\sum_{n=0}^{\infty}d_n\frac{1}{x\alpha_n}J_1(\lambda_nx)\frac{\sinh(\lambda_n(1-y))}{2\cosh(\lambda_n)\sinh(\lambda_n)}
\ee
where $\lambda_n$ is the $n^{th}$ zero of $J_1$, $\alpha_n$ is the normalizing constant from above, and the $d_n$ are the Bessel series coefficients for the function $x\mapsto{f}(x,-1)-A(x,-1)$.

Using that as $x\rightarrow0$ we have $x^{-1}I_1(\pi{n}x)\rightarrow\frac12\pi{n}$ and $x^{-1}J_1(\lambda_n{x})\rightarrow\frac12\lambda_n$, we have
\be
\begin{aligned}
A(0,y)&\;=\;\sum_{n=0}^{\infty} \frac{\pi{n}a_n/2}{I_1(\pi{n})}\cos(\pi{n}y)\,+\,\sum_{n=0}^{\infty}b_n\frac{\pi{n}b_n/2}{I_1(\pi{n})}\sin(\pi{n}y) \\
B(0,y)&\;=\;\sum_{n=0}^{\infty}\frac{\lambda_nc_n}{2\alpha_n}\frac{\sinh(\lambda_n(1+y))}{2\cosh(\lambda_n)\sinh(\lambda_n)} \\
C(0,y)&\;=\;\sum_{n=0}^{\infty}\frac{\lambda_nd_n}{2\alpha_n}\frac{\sinh(\lambda_n(1-y))}{2\cosh(\lambda_n)\sinh(\lambda_n)}
\end{aligned} \label{EqnsSeriesForABC}
\ee
Of course the sequences $\{a_n\}$, $\{b_n\}$, $\{c_n\}$, $\{d_n\}$ are square-summable (the Plancherel theorem), and since $I_1(n\pi)$ is exponentially decreasing with $n$, we easily see that the series for $A(0,y)$ and for every derivative $(\partial/\partial{y})^kA(0,y)$ is summable.
For $B$ and $C$, note that $\alpha_n\approx\sqrt{\lambda_n}$, and that for each $y\in(-1,1)$ the coefficients $\frac{\sinh(\lambda_n(1+y))}{2\cosh(\lambda_n)\sinh(\lambda_n)}$, $\frac{\sinh(\lambda_n(1-y))}{2\cosh(\lambda_n)\sinh(\lambda_n)}$ are exponentially decreasing with respect to $n$.
Therefore for $y\in(-1,1)$ each series in (\ref{EqnsSeriesForABC}) converges absolutely.
Taking derivatives with respect to $y$ in any of these series only creates coefficients that are polynomial in $\lambda_n$ (which grows like $n$), and so the series in (\ref{EqnsSeriesForABC}) all remain absolutely convergent no matter how many derivatives with respect to $y$ are taken.

The fact that $f(x,y)$ is identically equal to $A(x,y)+B(x,y)+C(x,y)$ on $Q$ follows from the fact that both expressions are finite and equal to one another on $\partial{Q}\setminus\{x=0\}$ (by construction), and then by Proposition \ref{PropUnspecifiabilityOfBoundaryValues}.
Thus $f$ is $C^\infty$ on $\overline{Q}$, including at $\{x=0\}$.
\qed

{\bf Remark}. Proposition \ref{PropRegularityAtBoundary} helps justify the colloquialism that points on the ``singular boundary,'' namely $\partial\Omega\cap\partial{H}^2$, are ``interior'' boundary points.
Not only are the values of $\varphi$ uniquely determined there by its values on $\partial\Omega\setminus\partial{H}^2$, but $C^\infty$ regularity holds.

{\bf Remark}.
Our regularity proof in Proposition \ref{PropRegularityAtBoundary} was a simple spin-off of the Fourier analysis explored above, but other proofs exist.
In the literature, the generalized Heston equation
\be
x\left(\sigma^2\varphi_{xx}+2\rho\sigma\varphi_{xy}+\varphi_{yy}\right)\,+\,\left(2c_0-q-x\right)\varphi_y\,+\,2\kappa\left(\theta-x\right)\varphi_x+c_0\varphi=0 \label{EqnHeston}
\ee
has seen substantial study in recent years, although it appears that our Liouville-type theorem was not proved, and the methods are substantially different.
Our equation $x\left(\varphi_{xx}+\varphi_{yy}\right)+\nu\varphi_x=0$ does not quite have the form (\ref{EqnHeston}), but does have the same behavior at the singular boundary.
See \cite{Feh} for existence/uniqueness of the Dirichlet problem, and \cite{DH}, \cite{FP1}, \cite{FP2} (and references therein) for $C^\infty$ regularity on non-singular boundary components.
In \cite{DH} the parabolic version of our equation $\varphi_t=x(\varphi_{xx}+\varphi_{yy})+\nu\varphi_x$, specifically, was studied for $\nu>0$.
Notice that our conclusion in Proposition \ref{PropRegularityAtBoundary} specifically fails to deal with the ``corner points'' of the domain.
This is dealt with at length in \cite{FP1}.

\begin{lemma}[Reduction to the exponential case] \label{LemmaReductionToExponential}
Assume $f\ge0$ is $C^0(\overline{H^2})\cap{C}^\infty(H^2)$ solves $f_{xx}+f_{yy}+3x^{-1}f_x=0$ on the right half-plane.
Then there is another function $\tilde{f}>0$ on $H^2$ with $\tilde{f}\in{C}^0(\overline{H^2})\cap{C}^\infty(H^2)$ so that the function $y\mapsto\tilde{f}(0,y)$ has at worst exponential growth.
If $f$ is non-constant, then $\tilde{f}$ is non-constant.
\end{lemma}
{\it Proof}. 
So assume there is such an $f\ge0$ that is not constant.
Proposition \ref{PropRegularityAtBoundary} actually guarantees $f\in{C}^\infty(\overline{H^2})$, and Corollary \ref{CorAlmostIncreasingFunction} shows that actually $f>0$ on $\overline{H^2}$.
We can thus consider the function $y\mapsto\left.\frac{d}{dy}\right|_{(0,y)}\log{f}$.
If this is bounded, then obviously $y\mapsto{}f(0,y)$ has exponential growth and putting $\tilde{f}=f$, we are done.

Otherwise, there exists some sequence $\{y_i\}$ such that $\lim_{i\rightarrow\infty}\left|\left.\frac{d}{dy}\right|_{(0,y_i)}\log{f}\right|=\infty$.
For each $i$, we can re-choose the point $y_i$ so that $\left|\frac{d}{dy}\log{f}\right|$ is ``almost largest'' in an appropriate neighborhood.
Specifically, set suppose $M_i=\left|\left.\frac{d}{dy}\right|_{(0,y_i)}\log{f}\right|$, and suppose there is some $y\in(y_i-2^iM_i,y_i+2^iM_i)$ with $\left|\left.\frac{d}{dy}\right|_{(0,y)}\log{f}\right|\;>\;2\left|\left.\frac{d}{dy}\right|_{(0,y_i)}\log{f}\right|$, then re-choose $y_i$ to be this new $y$.
Continuing, we eventually obtain a sequence $y_i$ so that
\begin{itemize}
\item[{\it{i}})] $\lim_{i\rightarrow\infty} \left|\left.\frac{d}{dy}\right|_{(0,y)}\log{f}\right|\rightarrow\infty$
\item[{\it{ii}})] $\left|\left.\frac{d}{dy}\right|_{(0,y)}\log{f}\right|\le{4}\left|\left.\frac{d}{dy}\right|_{(0,y_i)}\log{f}\right|$ on $(y_i-2^iM_i,\,y_i+2^iM_i)$, where $M_i=\left|\left.\frac{d}{dy}\right|_{(0,y)}\log{f}\right|$.
\end{itemize}
Finally, consider the sequence of functions
\be
f_i(x,y)\;=\;f(x/M_i,\,(y-y_i)/M_i) \; / \; f(0,y_i).
\ee
These have the properties
\begin{itemize}
\item[{\it{i}})] Each $f_i$ is non-negative and solves the PDE on $H^2$.
\item[{\it{ii}})] $f_i(0,0)=1$ and $\left|\left.\frac{d}{dy}\right|_{(0,0)}\log\,f_i\right|=1$.
\item[{\it{iii}})] We have $\left|\left.\frac{d}{dy}\right|_{(0,y)}\log{f}_i\right|\le{4}$ on the interval $(-2^i,\,2^i)$
\end{itemize}
Combining ({\it{ii}}) with Corollary \ref{CorAlmostIncreasingFunction} shows that $f_i$ is uniformly bounded at least at the point $(1,0)$.
Corollary \ref{CorPolyAndExpBound} implies $f_i$ does not converge to $\infty$ anywhere, and therefore $f_i$ converges to a $C^\infty(\overline{H^2})$ solution $\tilde{f}$.
From ({\it{iii}}) we have that $\tilde{f}(0,y)$ is at worst exponential.
Of course ({\it{ii}}) implies $\tilde{f}$ is non-constant.
Proposition \ref{ThmExpBoundsImplyConstancy} now provides the contradiction.
\qed

Finally we write down our primary analytical theorems.
\begin{theorem} \label{ThmFinalConstancy}
If $f\in{C}^0(\overline{H^2})\cap{C}^\infty(H^2)$, $f\ge0$, and $f$ solves $f_{xx}+f_{yy}+3x^{-1}f_x=0$, then $f$ is constant.
\end{theorem}
{\it Proof}.
By Lemma \ref{LemmaNewEqnBondaryFacts} and Proposition \ref{PropRegularityAtBoundary} we have that $f>0$ and $f\in{C}^\infty(\overline{H^2})$.
By Lemma \ref{LemmaReductionToExponential}, we may assume that $y\mapsto{f}(0,y)$ has exponential growth at worst.
Then Proposition \ref{ThmExpBoundsImplyConstancy} implies $f$ is constant.
\qed

\begin{theorem} \label{ThmFinalProportionalityToXSquared}
Assume $\varphi\in{C}^2(\overline{H^2})\cap{C}^\infty(H^2)$, $\varphi\ge0$, $\varphi=0$ at $\{x=0\}$, and $\varphi$ solves $\varphi_{xx}+\varphi_{yy}-x^{-1}\varphi_x=0$.
Then $\varphi(x,y)=A\,x^2$ for some constant $A\ge0$.
\end{theorem}
{\it Proof}.
We have that $f(x,y)=\varphi(x,y)\cdot{x}^{-2}$ solves $f_{xx}+f_{yy}+3x^{-1}f_x=0$.
By Lemma \ref{LemmaNewEqnBondaryFacts} $f$ is actually $C^0(\overline{H^2})\cap{C}^\infty(H^2)$.
Therefore Theorem \ref{ThmFinalConstancy} implies $f$ is constant, and so $\varphi(x,y)=A\cdot{x}^2$ for some constant $A\ge0$.
\qed

Here we apply our analytic results to the geometric situation.
\begin{theorem}[Potentials on half-plane polytopes] \label{ThmHalfPlaneFlat}
Assume the polytope $(\Sigma^2,g_\Sigma)$ is a closed half-plane and obeys (A)-(F).
Possibly after affine recombination of $\varphi^1$, $\varphi^2$, we have
\be
\varphi^1\;=\;\frac12x^2, \quad \varphi^2=y.
\ee
\end{theorem}
{\it Proof}.
After a possible affine re-combination of $\{\varphi^1,\varphi^2\}$, the image of the moment map $\Phi=(\varphi^1,\varphi^2)$ is precisely the set $H^2=\{\varphi^1\ge0\}$.
Then $x=\sqrt{\mathcal{V}}$ is zero on the boundary, and Proposition \ref{PropXISOneOneOnto} guarantees the complex coordinate $z=x+\sqrt{-1}y$ is 1-1, onto, $C^\infty$, and has $C^\infty$ inverse.
In real terms, the map $(x,y)$ is actually a coordinate system on $\Sigma^2$, sending the half-plane $\{\phi^1\ge0\}$ to the half-plane $\{x\ge0\}$.
Because the map $(x,y):\Sigma^2\rightarrow\{x\ge0\}$ has $C^\infty$ inverse, the functions $\varphi^1$, $\varphi^2$ are $C^0(\overline{H^2})\cap{}C^\infty(H^2)$ where $H^2=\{x>0\}$.

By Proposition \ref{PropCTwoAndLipschitz}, we have $\varphi^1,\varphi^2\in{}C^2(\overline{H^2})\cap{C}^\infty(H^2)$.
Since $\varphi^1(0,y)=0$, Theorem \ref{ThmFinalProportionalityToXSquared} guarantees $\varphi^1=\frac{\alpha}{2}x^2$.

As for $\varphi^2$, by ({\it{ix}}) of \S\ref{SubSecPolyTopeEdges} we have $\varphi^2=y+C_2(y)x^2+C_3(y)x^3+\dots$.
Thus the function $\varphi^2-y$ is zero on $\{x=0\}$ so that Theorem \ref{ThmFinalProportionalityToXSquared} gives $\varphi^2-y=\frac{\beta}{2}x^2$.
After possibly recombination of $\varphi^1$, $\varphi^2$, we therefore have $\varphi^1=\frac12x^2$, $\varphi^2=y$.

We have the transition matrix $A=(\frac{\partial\varphi^i}{\partial{x}^j})=\left(\begin{array}{cc} x & 0 \\ 0 & 1 \end{array}\right)$, so $\frac{\det(A)}{x}=1$ and from (\ref{EqnMetricConstruction}) and (\ref{EqnCoordSigmaScalarComp}) we have $g_\Sigma=dx^2+dy^2$ and $K_\Sigma=0$.
In momentum coordinates, from (\ref{EqnMetricConstruction}) we have $g_{\Sigma}=\frac{1}{4\varphi^1}(d\varphi^1)^2+(d\varphi^2)^2=\left(d\sqrt{\varphi^1}\right)^2+(d\varphi^2)^2$.

Using (\ref{EqnsGJOmegaM}) this gives the metric $g$ on $M^4$:
\be
g=\left(d\sqrt{\varphi^1}\right)^2+4\varphi^1\left(d\theta^1\right)^2\,+\,\left(d\varphi^2\right)^2\,+\,\left(d\theta^2\right)^2.
\ee
Substituting, say, $\varphi^1=r^2$, we easily see this is the flat metric on $\mathbb{R}^2\times\mathbb{R}^2$, where the Killing field on the first factor is a rotational field and on the second factor is a translation field.
\qed

\begin{corollary}[Half-plane polytopes are flat] \label{CorHalfPlaneFlat}
Assume $(M^4,J,\omega,\mathcal{X}_1,\mathcal{X}_2)$ has polytope $(\Sigma^2,g_\Sigma)$ where $\Sigma^2$ is a closed half-plane, and assume $s=0$ on $M^4$.
Then, possibly after affine recombination of $\varphi^1$, $\varphi^2$, we have
\be
\varphi^1\;=\;\frac12x^2, \quad \varphi^2=y.
\ee
Further, the metric $g$ on $M^4$ is flat.
\end{corollary}
{\it Proof}.
The K\"ahler reduction satisfies the hypotheses of Theorem \ref{ThmHalfPlaneFlat}.
\qed

\subsection{Examples} \label{SubSectionAnalyticExamples}

In this, our first of two ``examples'' section, we focus on analytical examples chosen to demonstrate the themes and the limitations of our analytical lemmata.
See \S\ref{SubSectionGeometricExamples} for geometrically-themed examples.

\noindent{\bf Example 1}. {\it Superharmonic solutions.}

The function
\be
\varphi(x,y)\;=\;\sqrt{2}-\sqrt{1-(x^2+y^2)+\sqrt{(1-(x^2+y^2))^2+4y^2}}
\ee
satisfies $\varphi(0,y)=0$, and solves $\varphi_{xx}-x^{-1}\varphi_x+\varphi_{yy}=0$ away from a singular set.
The point $(1,0)$ can be seen to be a branch point, and the singular ray is a branch cut.
The function $f(x,y)=\varphi(x,y)\cdot{x}^{-2}$ is smooth at $\{x=0\}$, and has the same singular locus.
\noindent\begin{figure}[h!]
\centering
\caption{\it Graphs for Example 1: Supersolutions.}
\noindent\begin{subfigure}[b]{0.4\textwidth}
\includegraphics[scale=0.4]{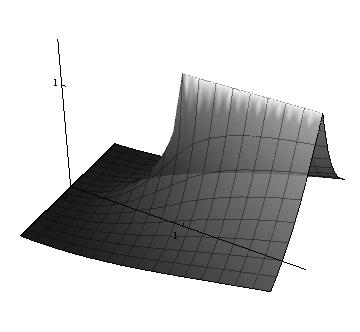}
\label{FigBranchCut}
\caption{Graph of $\varphi$, solving $\varphi_{xx}+\varphi_{yy}-x^{-1}\varphi_x\le0$ a.e., showing a branch cut.}
\end{subfigure}
\quad\quad\quad
\quad\quad\quad
\begin{subfigure}[b]{0.4\textwidth}
\includegraphics[scale=0.4]{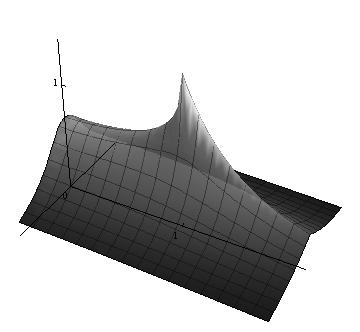}
\caption{Graph of $f=\varphi\cdot{x}^{-2}$ solving $f_{xx}+f_{yy}+3x^{-1}f_x\le0$ a.e., showing an internal salient.}
\end{subfigure}
\end{figure}

Neither $\varphi$ nor $f$ is Lipschitz, but are Holder continuous with Holder exponent $\frac12$.
The gradient does not exists along the ray $\{(x,0)\,\big|\,x\ge1\}$ so Proposition \ref{PropIntGradientBounds} is obviously meaningless, but we note that the conclusions of Corollary \ref{CorPolyAndExpBound} and Corollary \ref{CorPolyGrowthInT} actually hold with $c=2$.

In the general superharmonic case, provided some version of Corollary \ref{CorPolyAndExpBound} actually holds, then the proofs of Corollary \ref{CorAlmostIncreasingFunction}, Corollary \ref{CorPolyGrowthInT}, and Corollary \ref{CorBoundedAtSEqZero} actually go through.
In this example, we see the conclusions of all of these results indeed hold.

\noindent{\bf Example 2}. {\it Solutions that are Step Functions and $\delta$-functions on $\{x=0\}$}.

Consider the functions
\be
\begin{aligned}
&\varphi(x,y)\;=\;\frac12\left(1+\frac{y}{\sqrt{x^2+y^2}}\right) \\
&\psi(x,y)\;=\;\frac12\frac{x^2}{(x^2+y^2)^{\frac32}}.
\end{aligned} \label{EqnUnitStepAndDiracDelta}
\ee
These both solve the PDE (\ref{EqnOrigUnmodifiedPDE}) for $\nu=-1$ and are non-negative.
Note that $\psi=\frac{\partial\varphi}{\partial{y}}$.
Also, we have $y\mapsto\varphi(0,y)$ is the unit step function, and therefore $y\mapsto\psi(0,y)$ is the unit Dirac delta function.
This justifies the assertion, from the beginning of \S\ref{SubSectionFundSolMethods}, that $G(x,y)=\psi(x,y)$ is a fundamental solution on the right half-plane.

\noindent\begin{figure}[h!]
\centering
\caption{\it Graphs for Example 2: Steps and Deltas on $\{x=0\}$.}
\noindent\begin{subfigure}[b]{0.4\textwidth}
\includegraphics[scale=0.5]{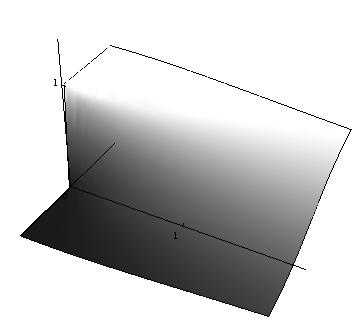}
\caption{Graph of non-negative solution $\varphi$ where $y\mapsto\varphi(0,y)$ is a step.}
\label{FigStepExample}
\end{subfigure}
\quad\quad\quad
\quad\quad\quad
\begin{subfigure}[b]{0.4\textwidth}
\includegraphics[scale=0.5]{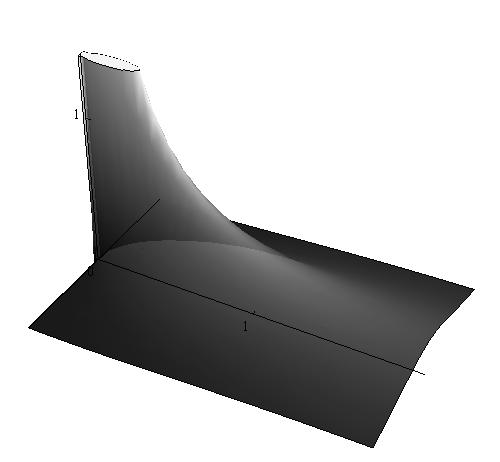}
\caption{Graph of non-negative solution $\psi$ where $y\mapsto\psi(0,y)$ is a Dirac-$\delta$.}
\end{subfigure}
\end{figure}

Setting $f=\varphi\cdot{x}^{-2}$ and $g=\psi\cdot{x}^{-2}$, we have non-negative functions on the half-plane that satisfy (\ref{EqnOrigModifiedPDE}), but which are not $C^0(\overline{H^2})$.\
In the case of $f$, we have a solution with finite boundary values on $x<0$ and infinite values on $x\ge0$.
In the case of $g$, we have a solution with finite boundary values on $x\ne0$, and $g(0,0)=\infty$.

\noindent{\bf Example 3}. {\it Solution that is a pulse on $\{x=0\}$}.

Solutions are translation-invariant in $y$, so using translations of (\ref{EqnUnitStepAndDiracDelta}) we construct
\be
\begin{aligned}
&\varphi(x,y)\;=\;\frac12\left(
\frac{y+1/2}{\sqrt{x^2+(y+1/2)^2}}
-\frac{y-1/2}{\sqrt{x^2+(y-1/2)^2}}
\right) \\
\end{aligned}
\ee
which solves $\varphi_{xx}-x^{-1}\varphi_{x}+\varphi_{yy}=0$.
Restricted to the boundary $\{x=0\}$, this is a unit pulse.
The associated function $f=\varphi\cdot{x}^{-2}=(x^2+y^2)^{-3/2}$ solves $f_{xx}+f_{yy}+3x^{-1}f_x=0$.

\noindent\begin{figure}[h!]
\centering
\caption{\it Graphs for Example 3: A pulse, and a solution of (\ref{EqnOrigModifiedPDE}) not determined by boundary values.}
\noindent\begin{subfigure}[b]{0.4\textwidth}
\includegraphics[scale=0.5]{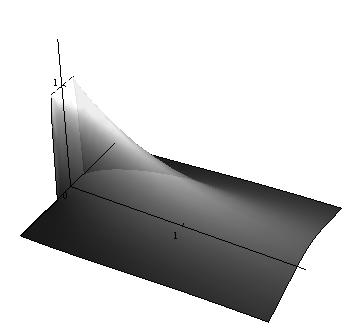}
\label{FigStep}
\caption{Graph of non-negative $\varphi$ where $y\mapsto\varphi(0,y)$ is a pulse function. \\}
\end{subfigure}
\quad\quad\quad
\quad\quad\quad
\begin{subfigure}[b]{0.4\textwidth}
\includegraphics[scale=0.5]{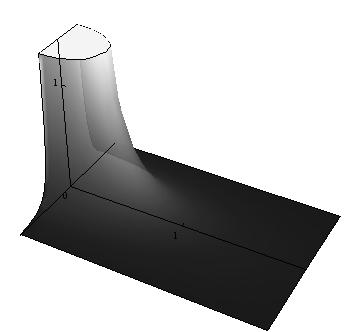}
\caption{Graph of non-negative, non-constant $f$ with finite and infinite values on the boundary.}
\end{subfigure}
\end{figure}
\indent Then $f$ has finite values on $x\notin[-\frac12,\frac12]$.
If one considers, say, the quadrilateral $Q=\{0<x<1\}\cap\{-1<y<1\}$, then we see that $f$ has finite values on $\partial{Q}$, and so using the Fourier and Bessel series from (\ref{EqnsSeriesForABC}) we can construct an $L^\infty$ solution $\tilde{f}$ to the PDE that has $\tilde{f}=f$ on $\partial{Q}$.
In particular, the function $f-\tilde{f}$ is zero on $\partial{Q}$ but is non-constant.

This illustrates the necessity of the hypothesis of boundedness in Lemma \ref{PropUnspecifiabilityOfBoundaryValues}.

\noindent{\bf Example 4}. {\it A solution that is only $C^{1,\alpha}(\overline{H^2})\cap{}C^\infty(H^2)$}.

Consider the function
\be
\varphi\;=\;{y}\sqrt{x^2+y^2}\,+\,{x^2}\log\left(\frac{y+\sqrt{x^2+y^2}}{x}\right),
\ee
which solves $\varphi_{xx}+\varphi_{yy}-x^{-1}\varphi_x=0$.
Note that $\varphi\in{}C^{1,\alpha}(\overline{H^2})\cap{C}^\infty(H^2)$ for all $\alpha\in[0,1)$ but $\varphi\notin{}C^{1,1}(\overline{H^2})$.

The fact that $y\mapsto\varphi(0,y)$ is piecewise quadratic rather than piecewise-linear provides an important counterpoint to the crucial assertion in Proposition \ref{PropCTwoAndLipschitz} where piecewise linearity is asserted when the solution $\varphi$ is a momentum function in a geometric situation.

\noindent\begin{figure}[h!]
\centering
\caption{\it Graph for Example 4: Piecewise quadratic growth at $\{x=0\}$.}
\includegraphics[scale=0.5]{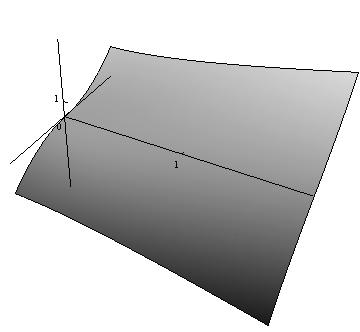}
\label{FigQuadraticBound}
\end{figure}
\indent This example shows the particular need for Propostion \ref{PropCTwoAndLipschitz} where we show that the solutions that arise from our geometrical situations are in fact $C^2$ up to the boundary, except at corners where they are Lipschitz.

Finally, this example show that although the PDE (\ref{EqnOrigUnmodifiedPDE}) has fundamental solutions (\ref{EqnFundSolForOrigEqn}) and has the ability to specify boundary values, but it is not hypoelliptic at the boundary $\{x=0\}$.
We have already seen that (\ref{EqnOrigModifiedPDE}) is not hypoelliptic at the boundary.

\section{The Quarter-Plane and General Polytopes} \label{SectionQuarterPlaneAndGeneral}

\subsection{Polytope has one corner}

First we consider the case that the polytope has just one corner.
Take the example
\be
\begin{aligned}
\varphi^1&\;=\; a y\,+\,b \sqrt{x^2\,+\,y^2}\,+\,\frac{\alpha}{2}x^2, \\
\varphi^2&\;=\; c y\,+\,d \sqrt{x^2\,+\,y^2} \,+\, \frac{\beta}{2}x^2.
\end{aligned}
\ee
Setting $\Phi=(\varphi^1,\varphi^2)^T$ we have Jacobian and Jacobian determinant
\be
\begin{aligned}
&D\Phi
\;=\;\left(\begin{array}{ccc}
bx\left(x^2+y^2\right)^{-\frac12}+\alpha{}x & & a+by\left(x^2+y^2\right)^{-\frac12} \\
dx\left(x^2+y^2\right)^{-\frac12}+\beta{}x & & c+dy\left(x^2+y^2\right)^{-\frac12}
\end{array}\right), \\
&\det\left(D\Phi\right)
\;=\;-x\left(x^2+y^2\right)^{-\frac12}\left((ad-cb)+\left(a\beta-c\alpha\right)\left(x^2+y^2\right)^{\frac12}+y\left(b\beta-\alpha{d}\right)\right).
\end{aligned}
\ee
It is impossible that $\det(D\Phi)$ be zero away from $\{x=0\}$ or undefined away from $(x,y)=(0,0)$, or else $d\varphi^1$ and $d\varphi^2$ would be colinear which would imply the polytope has some additional edge or corner somewhere.
Therefore $r=\sqrt{x^2+y^2}$ and $y$ must obey the linear inequalities
\be
\begin{aligned}
&\det(M)\,+\,\det(N)r\,+\,\det(K)y\;\ne\;0, \quad\quad 0\le|y|\le{r} \label{SystemForRY}
\end{aligned}
\ee
where
\be
\begin{aligned}
&M=\left(\begin{array}{cc} a & b \\ c & d \end{array}\right), \quad 
N=\left(\begin{array}{cc} a & \alpha \\ c & \beta \end{array}\right), \quad 
K=\left(\begin{array}{cc} b & \alpha \\ d & \beta \end{array}\right).
\end{aligned} \label{EqnsMatricesDefs}
\ee
In order for (\ref{SystemForRY}) to hold, it is necessary and sufficient that $\left|\frac{\det(K)}{\det(N)}\right|\le1$ and $\frac{\det(K)}{\det(M)}\ge0$ and $\det(M)\ne0$.

\begin{lemma} \label{LemmaAnalyticOneVertexLemma}
Consider the map $\Phi:\overline{H^2}\rightarrow\mathbb{R}^2$ where $\Phi=(\varphi^1,\varphi^2)^T$, and each $\varphi^i$ satisfies $(\varphi^i)_{xx}-x^{-1}(\varphi^i)_x+(\varphi^i)_{yy}=0$.
Assume $\Phi$ is $C^2$ except at $(0,0)$, and that $\Phi$ is Lipschitz everywhere.
Then either the image is either the half-plane and (possibly after $GL(2,\mathbb{R})$ recombination  of $\varphi^1$, $\varphi^2$) we have
\be
\begin{aligned}
&\varphi^1\;=\;x^2 \\
&\varphi^2\;=\;y
\end{aligned}
\ee
or else (possibly after a $GL(2,\mathbb{R})$ recombination) the polytope is the first quadrant with
\be
\begin{aligned}
\varphi^1&\;=\; \frac{1}{\sqrt{2}}\left(-y\,+\,\sqrt{x^2\,+\,y^2}\right)\,+\,\frac{\alpha}{2}x^2, \\
\varphi^2&\;=\; \frac{1}{\sqrt{2}}\left(y\,+\,\sqrt{x^2\,+\,y^2}\right)\,+\,\frac{\beta}{2}x^2
\end{aligned} \label{EqnsOneVertexPotentials}
\ee
for some $\alpha$, $\beta$ where $\alpha,\beta\ge0$.
\end{lemma}
{\it Proof}.
As indicated set $\Phi(x,y)=(\varphi^1(x,y),\varphi^2(x,y))^T$.
At any smooth point on $\{x=0\}$, the fact that $(\varphi^i)_{xx}-x^{-1}(\varphi^i)_x+(\varphi^i)_{yy}=0$ implies that $(\varphi^i)_x=0$ to first order.
But then $\lim_{x\rightarrow0}x^{-1}(\varphi^i)_x=(\varphi^i)_{xx}$, so the first two terms cancel and we have $(\varphi^i)_{yy}=0$ at $\{x=0\}$.

Therefore, away from $(0,0)$, the map $y\mapsto\Phi(0,y)$ is actually linear.
Including $(0,0)$, this map is piece-wise linear, and consists of two rays joined at one point 
\footnote{For an example of how this can fail if the ``Lipschitz'' hypothesis is removed, see Example 2.
For an example of how this can fail if the ``$C^2$ except at $(0,0)$'' hypothesis is removed, see Example 4.}.
After possibly an affine recombination of $\varphi^1$, $\varphi^2$, we can therefore assume that $y\mapsto\Phi(0,y)$ is the map $y\mapsto\left(\frac{1}{\sqrt{2}}\left(-y+|y|\right),\,\frac{1}{\sqrt{2}}\left(y+|y|\right)\right)$.

Consider the function $\psi^1=\varphi^1-\frac{1}{\sqrt{2}}\left(-y+\sqrt{x^2+y^2}\right)$.
On $\{x=0\}$ we now have that $\psi^1(0,y)=0$.
We wish to use Theorem \ref{ThmFinalProportionalityToXSquared}, but we must show that $\psi^1\ge0$.
First consider $\tilde\psi^1=\psi^1+x^2$, which is still zero on $\{x=0\}$.
This may have negative components, but the negative part is within the wedge below the line $y=1-x$ (in fact, it is below the parabola $y=\frac12(1-x^2)$.
The function $\eta_\epsilon=-\epsilon(x^4-4(y-1)^2x^2)$ is zero on $x=0$ and on $y=1-\frac12x$; therefore $\eta_\epsilon<0$ on $y=1-x$.
Further, the growth of $\eta_\epsilon$ is quadratic, whereas the growth of the negative part $\tilde\psi^1_{-}$ is linear.

Thus for any $\epsilon>0$ we have $\eta_\epsilon<-\tilde\psi^1_{-}$ for any $\epsilon$.
Sending $\epsilon\rightarrow0$, we see that $\tilde\psi^1_{-}=0$.
Thus $\tilde\psi^1\ge0$.
Now we can use Theorem \ref{ThmFinalProportionalityToXSquared} to see there exists some $\alpha\in\mathbb{R}$ so that $\tilde\psi^1(x,y)=\frac{\tilde\alpha}{2}x^2$ and therefore necessarily $\varphi^1(x,y)=\frac{1}{\sqrt{2}}\left(-y+\sqrt{x^2+y^2}\right)+\frac{\alpha}{2}x^2$ for some constant $\alpha=\tilde\alpha-1$.
Similarly we have $\varphi^2=\frac{1}{\sqrt{2}}\left(y+\sqrt{x^2+y^2}\right)+\frac{\beta}{2}x^2$.
The fact that we must have both $\alpha,\beta\ge0$ follows simply from the considerations coming immediately after (\ref{EqnsMatricesDefs}).
\qed

In \cite{AS} it was noticed that there is a 2-parameter family of metrics on any open polytope $\Sigma^2$ without parallel lines that produces scalar flat metrics on $(M^4,J,\omega)$.
The next theorem shows that up to homothety these are {\it precisely all} metrics when the polytope has a single corner.
This is the geometric version of Theorem \ref{LemmaAnalyticOneVertexLemma}.
\begin{theorem}[cf. Theorem \ref{TheoremQuarterPlane}] \label{ThmGeometricSingleVertex}
Assume $(\Sigma,g_\Sigma)$ satisfies (A)-(F)---for instance it may be the reduction of some scalar flat $(M^4,J,\omega,\mathcal{X}_1,\mathcal{X}_2)$---and assume the associated metric polytope $(\Sigma^2,g_\Sigma)$ is closed and has a single corner.
Then (up to homothethy and $SL(2,\mathbb{R})$ re-combination) the momentum functions  $\varphi^1$, $\varphi^2$ necessarily have the form
\be
\begin{aligned}
\varphi^1&\;=\; \frac{1}{\sqrt{2}}\left(-y\,+\,\sqrt{x^2\,+\,y^2}\right)\,+\,\frac{\alpha}{2}x^2, \\
\varphi^2&\;=\; \frac{1}{\sqrt{2}}\left(y\,+\,\sqrt{x^2\,+\,y^2}\right) \,+\, \frac{\beta}{2}x^2
\end{aligned} \label{EqnMomentumFunctionsInTheorem}
\ee
where $\alpha,\beta\ge0$.
Therefore metric $g_\Sigma$ belongs precisely to a 2-parameter family of possibilities, parametrized by $\alpha$ and $\beta$.
\end{theorem}
{\it Proof}.
The ``closure'' assumption is that the polytope contains both of its rays: geometrically, neither ray can be ``infinitely far away''.
By Proposition \ref{PropCTwoAndLipschitz}, the potential functions $\varphi^1$, $\varphi^2$, when considered as functions of $(x,y)$ and as solutions to the PDE (\ref{EqnOrigUnmodifiedPDE}), are indeed $C^2$ on edges and are Lipschitz everywhere.
But Lemma \ref{LemmaAnalyticOneVertexLemma} then implies $\varphi^1$, $\varphi^2$ (after possible $GL(2,\mathbb{R})$-recombination) necessarily have the form given above.
\qed

\begin{corollary} \label{CorOneVertexMetricClassification}
If the K\"ahler manifold $(M^4,J,\omega,\mathcal{X}_1,\mathcal{X}_2)$ is scalar-flat and has a closed polytope with a single corner, then the metric is necessarily in the 2-parameter family given by the momentum functions (\ref{EqnMomentumFunctionsInTheorem}) and the recipe in \S\ref{SubSectionReconstructionOfTheMetric}.
\end{corollary}
{\it Proof}. Combine Theorem \ref{ThmGeometricSingleVertex} with the description in \S\ref{SubSectionReconstructionOfTheMetric}. \qed

\subsection{Open polytopes with more than one corner}

In the case where the polytope is unbounded but has more than two edges, the situation is a bit more complex.
Homothetic $GL(2,\mathbb{R})$ transformations of $\varphi^1$, $\varphi^2$ are able to stabilize precisely two vertices, but unlike the one-vertex case, the speed of parametrization of the two ray edges cannot also be fixed with such transformations.

Define the ``outline'' of the polytope to be the image of the map $y\mapsto\left(\varphi^1(0,y),\varphi^2(0,y)\right)^T$.
In the case that the polytope contains all of its boundaries, then Proposition \ref{PropCTwoAndLipschitz} asserts that this map is Lipschitz and piecewise linear.
After fixing two vertices in place via $GL(2,\mathbb{R})$-transformation, the speed at which each linear segment is traversed is an independent variable.
Therefore each of the $N$ edges contributes a degree of freedom to the space of polytope metrics.
Two additional degrees of freedom express the possible addition of $x^2$-terms to each momentum function.

\begin{theorem}[Generic Polytopes] \label{ThmDegreesOfFreedomGeneralPolytope}
Assume $(\Sigma^2,g_\Sigma)$ is a closed polytope that satisfies (A)-(F) and has $n\ge3$ edges.
Then the pair $(\varphi^1,\varphi^2)$, as a function of $(x,y)$, belongs to an $(n+2)$-parameter family of possibilities.
As a result, the metric $g_\Sigma$ belongs to an $(n+2)$-parameter family of possible polytope metrics.
\end{theorem}
{\it Proof}.
The proof is to give a recipe for reconstructing the metric on the polytope from the behavior of $(\varphi^1,\varphi^2)$ on the edges, and noting the degrees of freedom that appear.
As noted above, the outline map $y\mapsto\Phi(0,y)=\left(\begin{array}{c}\varphi^1(0,y) \\ \varphi^2(0,y)\end{array}\right)$ is Lipschitz and piecewise linear, so can be expressed as sums of terms of the form
\be
\Psi_i(y)=
\left(\begin{array}{c}
\alpha_i(y-y_i)+\beta_i|y-y_i|+c_i \\
\gamma_i(y-y_0)+\delta_i|y-y_0|+d_i
\end{array}\right),  \quad \Psi(0,y)=\sum_{i=1}^{n-1}\Psi_i(y).
\ee
Since $y$ is required to satisfy only $dy=-J_\Sigma{d}x$, we may add a constant to $y$ to make sure that $y_1=0$.
Given an outline (the image of $y\mapsto\Phi(0,y)$) along with the fact that the outline map is piecewise linear, the only freedom we have is choosing the speed of the parametrizations of each segment; this is $n$ degrees of freedom.
With the first corner occurring at $(\varphi^1,\varphi^2)=(0,1)$ we have $c_1=0$, $d_1=1$, and from the speed of each segment, we can determine all of the rest of the constants.

Then we simply replace $|y-y_i|$ with $\sqrt{x^2+(y-y_i)^2}$ to obtain
\be
\widetilde{\Phi}_i(x,y)\;=\;
\left(\begin{array}{cc}
\widetilde{\varphi}_i^1(x,y) \\
\widetilde{\varphi}_i^2(x,y)
\end{array}\right)
\;=\;
\left(\begin{array}{cc}
\alpha_i(y-y_i)\,+\,\beta_i\sqrt{x^2+(y-y_i)^2}\,+\,c_i \\
\gamma_i(y-y_i)\,+\,\delta_i\sqrt{x^2+(y-y_i)^2}\,+\,d_i \\
\end{array}\right).
\ee
We have that each $\widetilde{\varphi}_i^1$, $\widetilde{\varphi}_i^2$ satisfies the differential equation (\ref{EqnOrigUnmodifiedPDE}).
Since $\sum_i\widetilde{\Phi}_i(0,y)=\Phi(0,y)$ by construction, by Theorem \ref{ThmFinalProportionalityToXSquared} there must be constants $\alpha$, $\beta$ such that $\varphi^1(x,y)=\frac{\alpha}{2}\sum_i\varphi_i^1$ and $\varphi^2=\frac{\beta}{2}+\sum_i\varphi^2_i$.
Selection of $\alpha$, $\beta$ give the final two degrees of freedom.
\qed

\begin{corollary} \label{CorDegreesOfFreedomGeneralScalarFlat}
If the scalar-flat K\"ahler manifold $(M^4,J,\omega,\mathcal{X}_1,\mathcal{X}_2)$ has a closed polytope with $n\ge3$ edges, then the momentum map $\Phi=(\varphi^1,\varphi^2)^T$ must be a member of an $(n+2)$-parameter family of possibilities.
Thus the metric $g$ on $M^4$ lies within an $(n+2)$-parameter family of possible scalar flat metrics.
\end{corollary}
{\it Proof}. Combine Theorem \ref{ThmDegreesOfFreedomGeneralPolytope} with the description in \S\ref{SubSectionReconstructionOfTheMetric}.  \qed

{\bf Remark}.
To express Theorem \ref{CorDegreesOfFreedomGeneralScalarFlat} more intrinsically, note that if the polytope of $(M^4,J,\omega,\mathcal{X}_1,\mathcal{X}_2)$ has $n$ many edges, then $b_2(M^4)=n-2$.
Thus if $b_2(M^4)>0$ then there are $b_2(M^4)+4$ degrees of freedom, up to homothety, in choosing scalar-flat metrics on complete K\"ahler 4-manifolds with two commuting holomorphic symmetries.

{\bf Remark}.
For the polytope of 3 edges (where $b_2(M^4)=1$), the recipe of Theorem \ref{CorDegreesOfFreedomGeneralScalarFlat} is carried out in Example 6 below.

{\bf Remark}.
Although we do not discuss the case that the polytope is not closed (that it has ``edges'' that are infinitely far away), it is not difficult to imagine how they could be built up, namely by adding together potentials of the form found in Example 2, along with the potentials used in this section.
However we run into difficulties in the proof of uniqueness.
The issue is that, after subtracting copies of such potentials, it is still possible that the resulting functions are not zero everywhere on the boundary.
If one attempts to multiply such potentials $\varphi$ by $x^{-2}$ to get $f=x^{-2}\cdot\varphi$, there may be singularities on the boundary.
Thus the analysis of Section \ref{SectionAnalysisHalfPlane} fails, as it relies on $f\in{C}^0(\overline{H^2})\cap{C}^\infty(H^2)$.

\subsection{Examples} \label{SubSectionGeometricExamples}

{\bf Example 5}. {\it Construction of all doubly-invariant scalar-flat metrics on $\mathbb{C}^2$}.

These are precisely the Taub-NUTs, the achiral (or twisted) Taub-NUTs, and the flat metric.
All of these are ALF metrics, except the flat metric which is ALE.
According to our theorem, up to linear combinations of $\varphi^1$, $\varphi^2$, all quarter-plane examples have the form given in (\ref{EqnsOneVertexPotentials}), and using (\ref{EqnMetricConstruction}) we have polytope metric
\be
\begin{aligned}
&g_{\Sigma}\;=\;\frac{1\,+\,\frac{\alpha+\beta}{\sqrt{2}}\sqrt{x^2+y^2}+\frac{\alpha-\beta}{\sqrt{2}}y}{\sqrt{x^2+y^2}}\,\left(dx\otimes{d}x+dy\otimes{d}y\right).
\end{aligned}
\ee
Setting $M=\frac{\alpha+\beta}{\sqrt{2}}$ and $k=\frac{\alpha-\beta}{\alpha+\beta}$ (so we are free to choose any $M>0$ and any $k\in[-1,1]$) and changing to polar coordinates, we have
\be
\begin{aligned}
&g_{\Sigma}\;=\;\frac{1+Mr\left(1+k\sin\theta\right)}{r}\,\left(dr^2\,+\,r^2d\theta^2\right), \quad 0\le\theta\le\frac\pi2, \\
&K_{\Sigma}\;=\;-M\frac{1\,-\,Mk\,r\left(k+\sin\theta\right)}{\left(1+Mr\left(1+k\sin\theta\right)\right)^3}.
\end{aligned}
\ee
When $k=0$, these are the Taub-NUT metrics.
Also using (\ref{EqnMetricConstruction}), we can write down the metric in $(\varphi^1,\varphi^2)$-coordinates, which allows us to easily write down the metric on $M^4$ as well.
However for generic $\alpha,\beta>0$, the epression for $G=x(det(A))AA^T$ is rather complex.

\noindent {\bf Example 6}. {\it Construction of all doubly-invariant scalar-flat metrics on $T\mathbb{C}P^1$}.

These include the Eguchi-Hanson metric.
These metrics are represented by polytopes with three edges: two rays and an adjoining segment.
Under $GL(2,\mathbb{R})$ recombinations of $\varphi^1,\varphi^2$, this polytope is the standard polytope with vertices at $(0,1)$ and $(1,0)$.

\noindent\begin{figure}[h!]
\centering
\caption{\it Plot for Example 6: Polytope outline for toric metrics on $T\mathbb{C}P^2$, with direction of parametrization indicated.}
\includegraphics[scale=0.45]{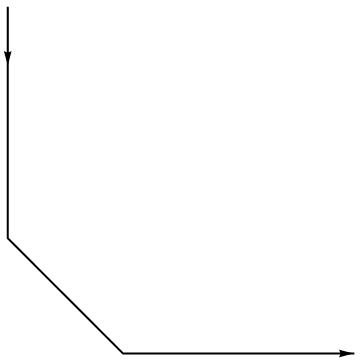}
\label{FigClippedPolytope}
\end{figure}
We have five remaining degrees of freedom: three degrees for the speed of each edge, and two degrees of freedom arising from the possible addition of $x^2$ terms to $\varphi^1$, $\varphi^2$.
Thus up to homothety all such metrics are given by
\be
\begin{aligned}
\varphi^1&=
\frac{v_2}{2\sqrt{2}}\left(y+\sqrt{x^2+y^2}\right)
+\left(\frac{v_3}{2}-\frac{v_2}{2\sqrt{2}}\right)\left(\left(y-\frac{\sqrt{2}}{v_2}\right)+\sqrt{x^2+\left(y-\frac{\sqrt{2}}{v_2}\right)^2}\right)
+\frac{\alpha}{2}x^2 \\
\varphi^2&=
1+\frac{v_1}{2}\left(-y+\sqrt{x^2+y^2}\right)
-\frac{v_2}{2\sqrt{2}}\left(y+\sqrt{x^2+y^2}\right) \\
&\quad\quad+\frac{v_2}{2\sqrt{2}}\left(\left(y-\frac{\sqrt{2}}{v_2}\right)+\sqrt{x^2+\left(y-\frac{\sqrt{2}}{v_2}\right)^2}\right)
+\frac{\beta}{2}x^2
\end{aligned}
\ee
The numbers $v_1,v_2,v_3>0$ encode the speeds at which $y\mapsto\Phi(0,y)$ maps onto the three linear pieces of the outline.
We remark that there are certainly other choices of parametrization.

\noindent{\bf Example 7}. {\it A polytope with an ``edge at infinity.''}

We begin with the 4-manifold $M^4=\mathbb{S}^2\times\mathbb{H}^2$, where $\mathbb{H}^2$ is the pseudosphere:
\be
g\;=\;dr^1\otimes{}dr^1+\sin^2(r^1)d\theta^1\otimes{d}\theta^1 \,+\,dr^2\otimes{}dr^2+e^{2r^2}d\theta^2\otimes{}d\theta^2.
\ee
The moment functions are $\varphi^1=-\cos(r^1)$ and $\varphi^2=e^{r^2}$.
Because ${r}^1\in[0,\pi]$ and $r^2\in(-\infty,\infty)$, we have that the range of $(\varphi^1,\varphi^2)$ is precisely $[0,1]\times(0,\infty)$.

\noindent\begin{figure}[h!]
\centering
\caption{\it Plots for Example 7: A polytope with an ``edge at infinity.''}
\noindent\begin{subfigure}[b]{0.4\textwidth}
\includegraphics[scale=0.65]{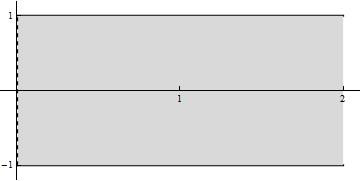}
\label{FigBoundaryAtInfinity}
\caption{The dashed segment is not part of the image of the moment map, and is infinitely far away as measured in $g_\Sigma$. \\}
\end{subfigure}
\quad\quad\quad
\quad\quad\quad
\begin{subfigure}[b]{0.4\textwidth}
\includegraphics[scale=0.5]{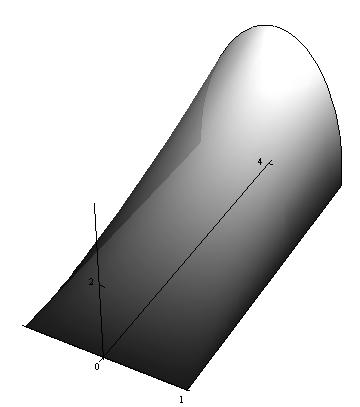}
\caption{Graph of $\sqrt{V}=x$ as a function of $\varphi^1$, $\varphi^2$.}
\end{subfigure}
\end{figure}
We have the metrics
\be
\begin{aligned}
g&\;=\;
\frac{1}{1-\left(\varphi^1\right)^2}\left(d\varphi^1\right)^2\,+\,\left(1-\left(\varphi^1\right)^2\right)\left(d\theta^1\right)^2
+\frac{1}{\left(\varphi^2\right)^2}\left(d\varphi^2\right)^2\,+\,\left(\varphi^2\right)^2\left(d\theta^2\right)^2, \\
g_\Sigma&\;=\;
\frac{1}{1-\left(\varphi^1\right)^2}\left(d\varphi^1\right)^2
+\frac{1}{\left(\varphi^2\right)^2}\left(d\varphi^2\right)^2
\end{aligned}
\ee
This produces
\be
\begin{aligned}
&x\;=\;\sqrt{\mathcal{V}}=|\nabla\varphi^1||\nabla\varphi^2|=\sqrt{1-\left(\varphi^1\right)^2}\cdot\varphi^2, \quad {\rm and} \quad
y\;=\;\varphi^1\varphi^2, \\
&\varphi^1\;=\;\frac{y}{x^2+y^2}, \quad {\rm and} \quad
\varphi^2\;=\;\sqrt{x^2+y^2}.
\end{aligned}
\ee
Figure \ref{FigStepExample} in Example 2 is actually the graph of $\varphi^1$ as a function of $x$, $y$.

\noindent {\bf Example 8}. {\it A Polytope with a disconnected outline}.

This also furnishes an example of a polytope $\Sigma^2$ so that $x=\sqrt{\mathcal{V}}$ with $y$ gives a complex function $z=x+\sqrt{-1}y$ with a critical point, so $z$ is not one-to-one; see Proposition \ref{PropXISOneOneOnto}.

Using the scalar flat metric on $\mathbb{S}^2\times\mathbb{H}^2$, where this time the hyperbolic piece has a ``two-ended trumpet'' metric, we have
\be
g\;=\;dr^1\otimes{}dr^1+\sin^2(r^1)d\theta^1\otimes{d}\theta^1 \,+\,dr^2\otimes{}dr^2+\cosh^2(r^2)d\theta^2\otimes{}d\theta^2.
\ee
This gives $\varphi^1=-\cos(r^1)$ and $\varphi^2=\sinh(r^2)$, with range $\Sigma^2=[-1,1]\times\mathbb{R}$.
The polytope metric
\be
g\;=\;\frac{1}{1-\left(\varphi^1\right)^2}\,d\varphi^1\otimes{d}\varphi^1\,+\,\frac{1}{1+\left(\varphi^2\right)^2}\,d\varphi^2\otimes{d}\varphi^2
\ee
and we have
\be
\begin{aligned}
&x\;=\;\sqrt{\left(1-\left(\varphi^1\right)^2\right)\left(1+\left(\varphi^2\right)^2\right)}, \quad y\;=\;\varphi^1\varphi^2, \\
&\varphi^1\;=\;\frac{1}{\sqrt{2}}\sqrt{1-\left(x^2+y^2\right)+\sqrt{\left(1-\left(x^2+y^2\right)^2\right)^2\,+\,4y^2}} \\
&\varphi^2\;=\;\frac{1}{\sqrt{2}}\sqrt{-1+\left(x^2+y^2\right)+\sqrt{\left(1-\left(x^2+y^2\right)^2\right)^2\,+\,4y^2}}
\end{aligned}
\ee
\noindent\begin{figure}[h!]
\centering
\caption{\it Plots for Example 8: A polytope with disconnected outline.}
\noindent\begin{subfigure}[b]{0.42\textwidth}
\includegraphics[scale=0.65]{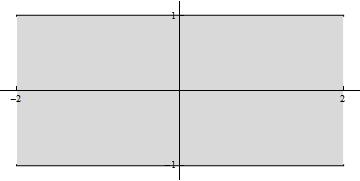}
\label{FigDisconnectedBoundary}
\caption{The polytope is the strip $-1\le\varphi^1\le1$. \\}
\end{subfigure}
\quad\quad\quad
\quad\quad\quad
\begin{subfigure}[b]{0.4\textwidth}
\includegraphics[scale=0.55]{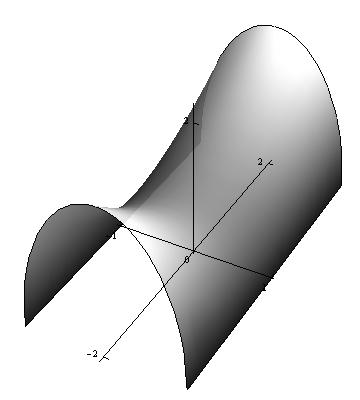}
\caption{Graph of $x=\sqrt{\mathcal{V}}$ as a function of $\varphi^1$, $\varphi^2$.}
\end{subfigure}
\end{figure}

Clearly $x$ and therefore $z=x+\sqrt{-1}y$ have a critical point.
This compliments Proposition \ref{PropXISOneOneOnto}, as the polytope has disconnected edges.

\end{document}